\documentclass[11pt]{article}

\usepackage{amsfonts,amsthm}
\usepackage{amsmath,mathtools}   
\usepackage{latexsym}
\usepackage{amssymb}
\usepackage{mathrsfs}
\usepackage{pifont}
\usepackage{tikz,ulem}
\usepackage[all]{xy}
\usepackage{hyperref}

\newtheoremstyle{newrem}{3pt}{3pt}{}{}
{\bfseries}{.}{.5em}{}

\newtheorem{theo}{Theorem}[section]
\newtheorem{thm}{Theorem}[section]

\newtheorem*{theo*}{Theorem}
\newtheorem{lemm}[theo]{Lemma}
\newtheorem{prop}[theo]{Proposition}

\theoremstyle{newrem}
\newtheorem{remax}[theo]{Remark}
\newenvironment{rema}
  {\pushQED{\qed}\remax}
  {\popQED\endremax}
\theoremstyle{definition}

\newtheorem*{term*}{Notation/Terminology}

\newcommand{\bZ}{\mathbb{Z}}

\newcommand{\cA}{\mathcal{A}}
\newcommand{\cB}{\mathcal{B}}

\newcommand{\cG}{\mathcal{G}}

\newcommand{\cK}{\mathcal{K}}

\newcommand{\cN}{\mathcal{N}}
\newcommand{\cP}{\mathcal{P}}

\newcommand{\cT}{\mathcal{T}}

\newcommand{\cW}{\mathcal{W}}

\newcommand{\bJ}{\boldsymbol{J}}

          \def\fe{{\mathfrak e}}
\def\fg{{\mathfrak g}}     \def\fh{{\mathfrak h}}     
          
     \def\fr{{\mathfrak r}}     \def\fs{{\mathfrak s}}
\def\ft{{\mathfrak t}}     \def\fu{{\mathfrak u}}
            \def\fii{{\mathfrak i}}

     \def\fS{{\mathfrak S}}     

\newcommand{\II}{{\mathbb I}}
\newcommand{\ZZ}{{\mathbb Z}}

    \def\sN{\mathsf{N}}

\newcommand\atopn[2]{\genfrac{}{}{0pt}{}{#1}{#2}}

\newcommand\ind[4]{{\scriptsize%
\!\!\begin{array}{cc}#2\\ #1\end{array}%
\!\!\!\!\begin{array}{cc}#4\\ #3\end{array}%
}}

\newcommand{\ket}[1]{|#1\rangle}
\newcommand{\bra}[1]{\langle #1|}

\def\sgn{\mathop{\rm sgn}\nolimits}

\usepackage{geometry}
\newcommand{\doi}[1]{\href{https://doi.org/#1}{doi:#1}}
\numberwithin{equation}{section}

\begin{document}

\title{\bf Representations  of the rank two Racah algebra\\
 and orthogonal multivariate polynomials}
\author{
Nicolas Cramp\'e\textsuperscript{$1$}\footnote{E-mail: crampe1977@gmail.com}~,
Luc Frappat\textsuperscript{$2$}\footnote{E-mail: luc.frappat@lapth.cnrs.fr}~,
Eric Ragoucy\textsuperscript{$2$}\footnote{E-mail: eric.ragoucy@lapth.cnrs.fr}~
\\[.9em]
\textsuperscript{$1$}
\small Institut Denis-Poisson CNRS/UMR 7013 - Universit\'e de Tours - Universit\'e
d'Orl\'eans,\\
\small~Parc de Grandmont, 37200 Tours, France.\\[.5em]
\textsuperscript{$2$}
\small Laboratoire d'Annecy-le-Vieux de Physique Th\'eorique LAPTh,\\
\small~Universit\'e Savoie Mont Blanc, CNRS, F-74000 Annecy,
 France.\\[.9em]
}
\date{}
\maketitle

\bigskip\bigskip 

\begin{center}
\begin{minipage}{14cm}
\begin{center}
{\bf Abstract}\\
\end{center}

The algebraic structure of the rank two Racah algebra is studied in detail. We provide an automorphism group of this algebra, which is isomorphic to the permutation group of five elements. This group can be geometrically interpreted  as the symmetry of a folded icosidodecahedron. It allows us to study a class of equivalent irreducible representations of this Racah algebra. They can be chosen symmetric so that their transition matrices are orthogonal. We show that their entries can be expressed in terms of Racah polynomials. 
This construction gives an alternative proof of the recurrence, difference and orthogonal relations satisfied by the Tratnik polynomials, as well as their expressions as a product of two monovariate Racah 
polynomials. Our construction provides a generalization of these bivariate polynomials together with their properties. 
\end{minipage}
\end{center}

\medskip

\begin{center}
\begin{minipage}{13cm}

\textbf{MSC2020 database:} 33C80; 33C45; 16G60
\end{minipage}
\end{center}

\clearpage
\tableofcontents
\newpage
\section{Introduction}\label{sec:intro}

The bispectral properties of the orthogonal polynomials of the Askey scheme \cite{Koek} can be described 
by algebras defined by quadratic relations \cite{Zhedanov}. When the 
Racah polynomial, which is at the top of the Askey scheme, is considered, one obtains an algebra, called the Racah algebra (of rank one). 
The $6j$-symbols, which are the overlap coefficients between 
bases with different recoupling of three $su(2)$
representations, can be expressed in terms of Racah polynomials \cite{Wilson,Wil}.
The Racah algebra provides an algebraic explanation of this connection \cite{LL,GZ,GVZ}. Indeed, the elements of the centralizer of the diagonal action of $su(2)$ in three copies of $su(2)$ 
satisfy the Racah algebra. Therefore, the $6j$-symbols are intimately connected to the representation theory of the Racah algebra.
This justifies the mathematical study of its representations \cite{HB,Huang}. 
This algebra also appears  in several other contexts. The connection between this algebra and the centralizer has been investigated in detail in \cite{CPV} where 
some quotient of the Racah algebra has been identified as the Temperley--Lieb or Brauer algebras. 
It can be embedded in different other algebras: $U(su(2))$ \cite{GZ2,Koo,BH2,CSV}, additive DAHA of type $(C_1^\vee , C_1 )$ \cite{Huang} or Bannai--Ito algebra \cite{GVZ3}.
The Racah algebra, as the algebras associated to the other polynomials of the Askey scheme, allows one to characterize the Leonard pair and 
plays an important role in algebraic combinatorics \cite{Ter,GWH}.
Finally, in the context of  superintegrable models, some Hamiltonians on the $2$-sphere possess the Racah algebra as a symmetry \cite{Post1,KMP1,GVZ2}.

The symmetry of the Hamiltonian of superintegrable models on $d$-sphere ($d> 2$) leads to identify the structure of the Racah algebra for higher rank.
In \cite{KMP}, the model on the $3$-sphere was examined and the Tratnik bivariate Racah polynomials \cite{Tra} have been identified as the bases of the symmetry algebra.
The generalization to a $d$-sphere has been made in \cite{DGVV1} where relations for the higher rank Racah algebra have been computed (see also \cite{LMZ}).
The rank two Racah algebra also showed up in the recoupling problem of four copies of $su(1,1)$ algebra in \cite{Post}. Then, the centralizer of the diagonal action
on $n$ copies has been studied in \cite{DIVV}. 
The exact isomorphism between this centralizer and a quotient of the higher rank Racah algebra, called special Racah algebra, has been proven in \cite{CGPV}.
The connection with multivariate orthogonal polynomials \cite{GI} has been also initiated in \cite{BV}.
\\ 

In this paper, we consider the rank two Racah algebra $R(4)$, which is the first generalization of the usual Racah algebra $R(3)$. Their definitions are recalled in details in section \ref{sec:defR4}.
We study a specific quotient of $R(4)$, called special Racah algebra $sR(4)$ defined in \cite{CGPV}, obtained by setting to zero certain Casimir elements of $R(4)$.
 An automorphism group of $sR(4)$ is identified in section \ref{sect:sym-icosi}, which is isomorphic to the permutation group of five elements $\fS_5$.
 It is also identified with the symmetry group of the folded icosidodecahedron, when its vertices are identified with abelian subalgebras of $sR(4)$, see figure \ref{fig}.

We study representations of $sR(4)$ in section \ref{sec:rep}. We focus on representations such that one of the abelian subalgebras acts diagonally. 
The vectors of the corresponding representation space are labeled by two integers, and we show that the other generators of $sR(4)$ act tridiagonally on these vectors, 
in the sense that the actions of these generators on a vector shift only by $\pm 1$ these two labels.
This representation depends on the values of the central elements of $sR(4)$. A finite-dimensional representation can be extracted by imposing 
a constraint on these values. We show that this representation is irreducible and can be chosen real symmetric when the parameters satisfy a set of inequalities.

In section \ref{sec:equiv}, the automorphisms of $sR(4)$ allow us to build a representation associated to each element of the permutation group $\fS_5$. 
We demonstrate that choosing accurately their parameters, these representations become in fact equivalent. 
The transition matrices between these equivalent representations are computed and their entries are identified with Racah polynomials.
This identification together with relations among automorphisms of $sR(4)$ allow us to obtain some relations on monovariate Racah polynomials, which turn out 
to be equivalent to the Racah and Biedenharn--Elliott relations (section \ref{sec:pmr}).

Finally, in section \ref{sec:bivar},  the transition matrices associated to paths on the folded icosidodecahedron are computed. 
For length two paths, the transition matrix entries are expressed as the product of two monovariate Racah polynomials, which is identified with the bivariate Tratnik polynomials. As a by-product, our construction of $sR(4)$ representations 
allows us to recover the unitarity, recurrence and difference identities obeyed by these objects.
Similarly, for length three paths, the transition matrix entries can be computed and appear to be a sum  of products of three monovariate Racah polynomials.
This provides a generalization of the bivariate Tratnik polynomials, which has not been studied yet.
Once again, our construction provides  unitarity, recurrence and difference identities obeyed by these objects. 
We conclude the paper with various open problems in section \ref{sec:conclu} and gather some technical formulas in appendices.

\section{Special Racah algebra $sR(4)$\label{sec:sR4}}

\subsection{Definitions and Casimir elements of the $R(4)$ algebra\label{sec:defR4}}

The Racah algebra $R(4)$ is generated by the elements $C_I$ with $I\subseteq \{1,2,3,4\}$. 
Let us emphasize that the order in the subset $I$ is irrelevant and, for example, $C_{12}=C_{21}$ (to lighten the notations, 
we write only the elements of the subsets instead of the subset itself. For example, $C_{12}$ stands for $C_{\{1,2\}}$). 
By convention $C_\emptyset=0$ and for two disjoint subsets $I,J$, one sets $C_{IJ}=C_{I \cup J}$.
The defining relations of the Racah algebra $R(4)$ are
\begin{align}
& [C_I,C_J]=0 \quad \text{for} \quad I \cap J =\emptyset\ \ \text{or}\ \ I\subseteq J \,,
\label{eq:commCICJ} \\
& C_I=\frac{1}{2}\sum_{i,j\in I \atop i\neq j} C_{ij} -( |I|-2 ) \sum_{i\in I} C_i\,,
\label{eq:defCI}
\end{align}
for three nonempty disjoint subsets $I,J,K\subseteq \{1,2,3,4\}$, 
\begin{align}
&[C_{IJ},C_{JK}]=[C_{JK},C_{IK}]=[C_{IK},C_{IJ}]\,, \label{eq:comm2} \\
&\big[C_{JK},[C_{IJ},C_{JK}]\big]=2 C_{IK}C_{JK}-2C_{JK}C_{IJ} + 2(C_J-C_K)(C_{IJK} -C_I)\,,
\label{eq:comm3} 
\end{align}
and, for $i, j, k,\ell \in  \{1,2,3,4\}$ pairwise distinct,
\begin{eqnarray}
 \frac{1}{2}\big[C_{k\ell}, [C_{ij},C_{jk}] \big]&=&C_{ik}C_{j\ell}-C_{i\ell}C_{jk} +(C_i+C_\ell)C_{jk}+(C_j+C_k)C_{i\ell}\nonumber \\
 &&-(C_j+C_\ell)C_{ik}-(C_i+C_k)C_{j\ell}+ (C_i-C_j)(C_\ell-C_k)\,, 
 \label{eq:comm3bis} \\
 \frac{1}{2}\big[ [C_{ij},C_{jk}] , [C_{jk},C_{k\ell}] \big]&=&- \big( [C_{ij},C_{j\ell}]+ [C_{ik},C_{k\ell}] \big) (C_{jk}-C_j-C_k) \,. \label{eq:comm4} 
\end{eqnarray}

The above presentation of $R(4)$ is not unique and different other possibilities exist in the literature.
Instead of $C_{ij}$, another set of generators $P_{ij}=C_{ij}-C_i-C_j$ (for $i\neq j$) and $P_{ii}=2C_i$ has been used previously 
to define the Racah algebra $R(4)$ \cite{DGVV1,DIVV,CGPV}.
In this paper, we prefer to define the $R(4)$ algebra using all the elements $C_I$ with $I\subseteq \{1,2,3,4\}$ since it simplifies the  computations.

One can show that the relations  \eqref{eq:comm4} can be replaced by (for pairwise  distinct indices $i,j,k,\ell\in\{1,2,3,4\}$):
\begin{equation}\label{eq:comm4bis}
2C_j\,[C_{ik},C_{k\ell}] +(C_{ij}-C_i-C_j)[C_{k\ell},C_{jk}]  +(C_{j\ell}-C_j-C_\ell) [C_{jk},C_{ik}] 
-(C_{jk}-C_j-C_k)[C_{ij},C_{j\ell}] 
=0\,.
 \end{equation}
Such types of relations, expressed in the $P_{ij}$ basis, appeared already in \cite{CGPV}.
Note that thanks to \eqref{eq:comm2}, relation 
\eqref{eq:comm4bis} is invariant under any permutation of the indices $(i,k,\ell)$, one gets only $4$ such relations.

Using relations \eqref{eq:defCI} (see also \eqref{eq:Ctrou}), we can define the algebra $R(4)$  using only the elements $C_{j}$, $j=1,2,3,4$, $C_{1234}$,
$C_{12}$, $C_{23}$, $C_{34}$, $C_{123}$ and $C_{234}$. 
For completeness, we provide the equivalent set of relations of $R(4)$ using only these elements in appendix \ref{app:algR4}.
They will be used for the construction of representations.\\

From relations \eqref{eq:commCICJ} satisfied in the Racah algebra $R(4)$, one can show that $C_i$ ($i\in \{1,2,3,4\}$) and $C_{1234}$ are central.  
There exist other central elements given by
\begin{eqnarray} \label{eq:w123}
 w_{I,J,K}&=&\frac{1}{4}\Big( [C_{IJ},C_{JK}]\Big)^2 -\frac{1}{2}\{C_{IJ}^2,C_{JK}\} -\frac{1}{2}\{C_{IJ},C_{JK}^2\}  + C_{IJ}^2+C_{JK}^2 +\{C_{IJ},C_{JK}\} \nonumber  \\
 &&+\frac{1}{2}(C_I+C_J+C_K+C_{IJK})( \{C_{IJ},C_{JK} \} -2C_{IJ}-2C_{JK} )\nonumber\\
&& -(C_I-C_{IJK})(C_J-C_K) C_{IJ}-(C_I-C_{J})(C_{IJK}-C_K) C_{JK}\nonumber\\
&&+(C_IC_K-C_{IJK}C_J)(C_{IJK}-C_I+C_J-C_K)+(C_I+C_K)(C_{IJK}+C_J)\,,
\end{eqnarray}
where $I,J,K\subseteq \{1,2,3,4\}$ are three nonempty disjoint subsets and $\{.,.\}$ is the anticommutator.
It should  be mentioned that  the  element $w_{1,2,3}$ is  essentially the known Casimir element of the Racah algebra $R(3)$ given in \cite{GZ,GVZ}. 
It is quite a surprise that this element is still a central element in the Racah algebra $R(4)$. For example, it is not trivial that it commutes with the generator $C_{34}$ \cite{CGPV}.
We can show that $w_{I,J,K}$ is symmetric under the exchange of the sets $I$, $J$ and $K$. 
In addition, one gets
\begin{equation}\label{eq:x1234}
 2x_{1234}:=w_{12,3,4}-w_{1,3,4}-w_{2,3,4}=w_{1,23,4}-w_{1,2,4}-w_{1,3,4}=w_{1,2,34}-w_{1,2,3}-w_{1,2,4}\,.
\end{equation}
Therefore, there are only $5$ independent Casimir elements $w_{1,2,3}$, $w_{1,2,4}$, $w_{1,3,4}$, $w_{2,3,4}$ and $x_{1234}$.
These 5 elements have been identified in \cite{CGPV}.

The special Racah algebra $sR(4)$ is the quotient of the Racah algebra $R(4)$ by the relations
\begin{equation}\label{eq:specialR4}
 w_{1,2,3}=0\ ,\qquad w_{1,2,4}=0\ ,\qquad w_{1,3,4}=0\ ,\qquad w_{2,3,4}=0\ , \qquad x_{1234}=0\,.
\end{equation}
The terminology ``special'' comes from \cite{CFGPRV}, where a similar quotient of the Askey--Wilson 
algebra has been defined and denoted as the “special Askey--Wilson algebra” (this was inspired by the nomenclature of Lie groups).

For any three distinct subsets $I,J,K\subseteq \{1,2,3,4\}$, the subalgebras of $sR(4)$ generated by $C_{I}$, $C_{J}$, $C_{K}$, $C_{IJ}$, $C_{IK}$, $C_{JK}$ and $C_{IJK}$
are isomorphic and describe the so-called special Racah algebra $sR(3)$.

\subsection{A web of abelian subalgebras}

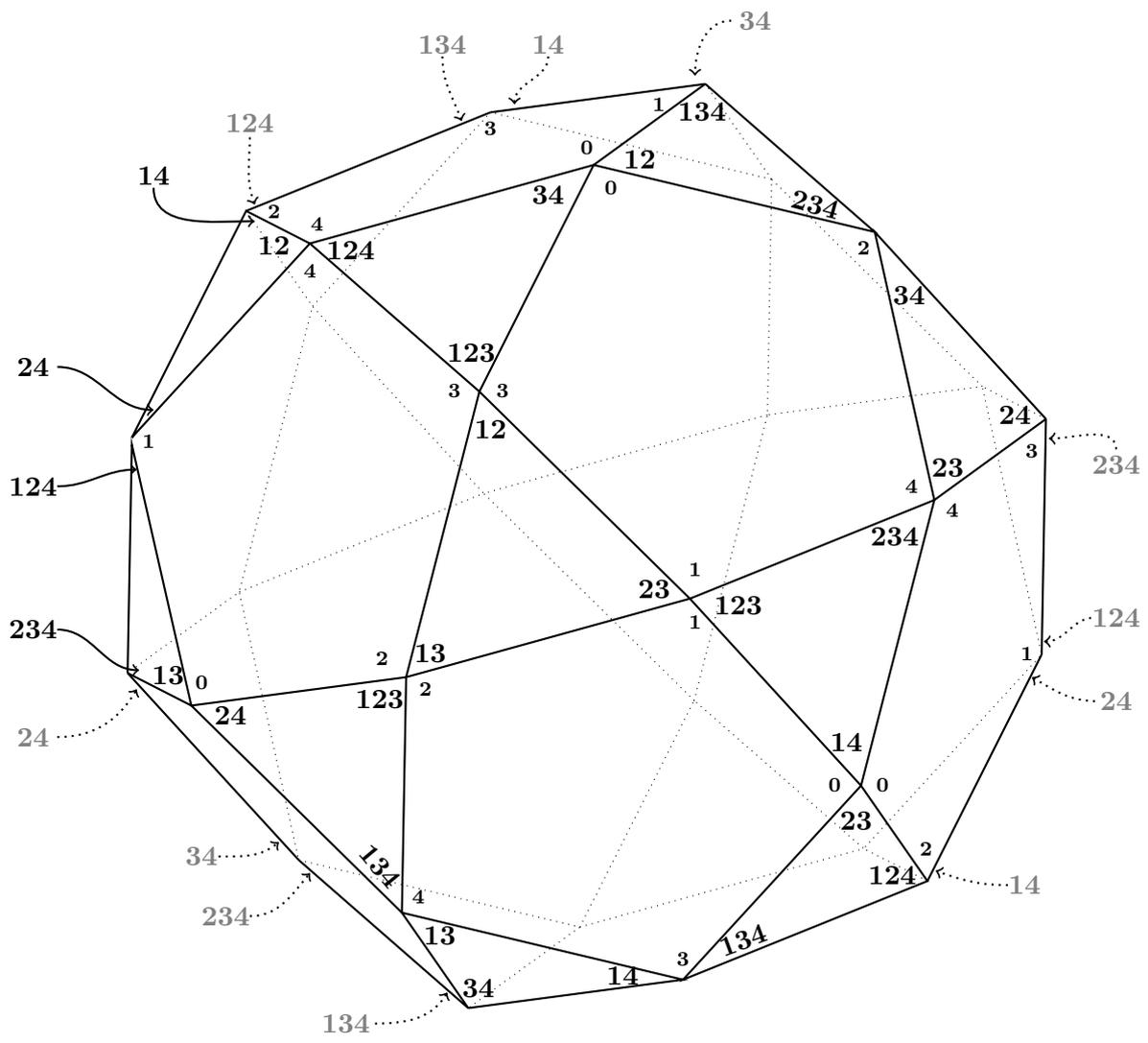
\begin{figure}[htbp]
\begin{center}
\begin{tikzpicture}[scale=3.3]

\draw [thick] (-1.908111787, -.5319769439)--(-1.890410045, .4242268192)-- (-1.643196449, -.6682790668) -- (-1.908111787, -.5319769439);

\node at (-0.48, 0.81) {$\mathbf{123}$};
\node at (-0.4, 0.49) {$\mathbf{12}$};
\node at (-0.35, 0.65) {\scriptsize$\mathbf{ 3}$};
\node at (-0.55, 0.65) {\scriptsize$\mathbf{ 3}$};

\node at (-0.65, -0.45) {$\mathbf{13}$};
\node at (-0.86, -0.64) {$\mathbf{123}$};
\node at (-0.67, -0.6) {\scriptsize$\mathbf{ 2}$};
\node at (-0.85, -0.47) {\scriptsize$\mathbf{ 2}$};

\node at (0.28, -0.18) {$\mathbf{23}$};
\node at (0.63, -0.25) {$\mathbf{123}$};
\node at (0.45, -0.1)  {\scriptsize$\mathbf{ 1}$};
\node at (0.45, -0.32) {\scriptsize$\mathbf{ 1}$};

\node at (0.22, 1.62) {$\mathbf{12}$};
\node at (-0.16, 1.47) {$\mathbf{34}$};
\node at (0.1, 1.5) {\scriptsize$\mathbf{ 0}$};
\node at (-0., 1.67) {\scriptsize$\mathbf{ 0}$};

\node at (-1.3, 1.26) {$\mathbf{12}$};
\node at (-0.98, 1.24) {$\mathbf{124}$};
\node at (-1.12, 1.35) {\scriptsize$\mathbf{ 4}$};
\node at (-1.15, 1.15) {\scriptsize$\mathbf{ 4}$};
\node at (-1.3, 1.4) {\scriptsize$\mathbf{ 2}$};

\node at (1.5, 0.33)  {$\mathbf{23}$}; 
\node at (1.28, 0.04) {$\mathbf{234}$};
\node at (1.52, 0.15) {\scriptsize$\mathbf{ 4}$}; 
\node at (1.35, 0.25) {\scriptsize$\mathbf{ 4}$};

\node at (1.08, -0.82) {$\mathbf{14}$};
\node at (1.12, -1.15) {$\mathbf{23}$};
\node at (1.03, -1) {\scriptsize$\mathbf{ 0}$};
\node at (1.23, -1.)  {\scriptsize$\mathbf{ 0}$};

\node at (0.15, -1.8) {$\mathbf{14}$}; 
\node at (0.65, -1.65) [rotate=20] {$\mathbf{134}$};
\node at (0.4, -1.73) {\scriptsize$\mathbf{ 3}$}; 

\node at (1.27, -1.38) {$\mathbf{124}$};
\node at (1.41, -1.27) {\scriptsize$\mathbf{ 2}$};

\node at (-0.61, -1.63) {$\mathbf{13}$}; 
\node at (-0.86, -1.34) [rotate=-50 ] {$\mathbf{134}$};
\node at (-0.7, -1.47) {\scriptsize$\mathbf{ 4}$}; 

\node at (1.83, -0.45) {\scriptsize$\mathbf{ 1}$};

\draw [dotted,thick,<-] (1.92,0.45) to [out=0,in=90] (2.2,0.4); 
\node at (2.2, 0.34) [color=black!50] {$\mathbf{234}$};

\node at (1.78, 0.55) {$\mathbf{24}$};
\node at (1.85, 0.4) {\scriptsize$\mathbf{ 3}$};

\node at (-0.45, -1.85) {$\mathbf{34}$};

\node at (-1.74, -0.54) {$\mathbf{13}$};
\node at (-1.48, -0.71) {$\mathbf{24}$};
\node at (-1.6, -0.57) {\scriptsize$\mathbf{ 0}$};

\node at (1.34, 1.05) {$\mathbf{34}$}; 
\node at (0.95, 1.43) [rotate=-15]   {$\mathbf{234}$};
\node at (1.15, 1.25) {\scriptsize$\mathbf{ 2}$}; 

\node at (0.48, 1.82) {$\mathbf{134}$};
\node at (0.3, 1.85)  {\scriptsize$\mathbf{ 1}$};
\node at (-0.4, 1.75) {\scriptsize$\mathbf{ 3}$};

\node at (-1.8, 1.55) {$\mathbf{14}$};
\draw [thick,->] (-1.8,1.5) to [out=-90,in=-180] (-1.38,1.36); 

\node at (-2.3, .75)  {$\mathbf{24}$};
\draw [thick,->] (-2.2, .75) to [out=0,in=180] (-1.8, .57); 
\node at (-1.82, .44) {\scriptsize$\mathbf{ 1}$};

\node at (-2.3, .25) {$\mathbf{124}$};
\draw [thick,->] (-2.2, .25) to [out=0,in=-180] (-1.87, .32); 

\node at (-2.3, -0.35)  {$\mathbf{234}$};
\draw [thick,->] (-2.2,- 0.35) to [out=0,in=-180] (-1.86, -.52); 

\node at (-1.4,1.77) [color=black!50] {$\mathbf{124}$};
\draw [thick,dotted,<-] (-1.38, 1.43) to [out=110,in=-90] (-1.4,1.72); 

\node at (-2.3, -0.8) [color=black!50]  {$\mathbf{24}$};
\draw [thick,dotted,->] (-2.2,- 0.8) to [out=0,in=-120] (-1.87, -.6); 

\node at (0.7, 2.2) [color=black!50]  {$\mathbf{34}$};
\draw [thick,dotted,->] (0.6,2.2) to [out=-180,in=90] (0.45, 1.97); 

\node at (-1., -2.) [color=black!50] {$\mathbf{134}$};
\draw [thick,dotted,->] (-0.88,-2) to [out=0,in=-120] (-0.58, -1.87); 

\node at (1.82, -1.42) [color=black!50]  {$\mathbf{14}$};
\draw [thick,dotted,->] (1.75,-1.42) to [out=180,in=-20] (1.45, -1.36); 

\node at (2.2, -0.65) [color=black!50]  {$\mathbf{24}$};
\draw [thick,dotted,<-] (1.85,- 0.55) to [out=-30,in=180] (2.12, -.65); 

\node at (2.2, -0.3) [color=black!50]  {$\mathbf{124}$};
\draw [thick,dotted,<-] (1.9,- 0.4) to  [out=0,in=180] (2.11, -.3); 

\node at (-1.6, -1.3) [color=black!50]  {$\mathbf{34}$};
\draw [thick,dotted,->] (-1.53,- 1.3) to [out=0,in=220] (-1.28, -1.24); 

\node at (-1.5, -1.55) [color=black!50]  {$\mathbf{234}$};
\draw [thick,dotted,->] (-1.4,- 1.55) to [out=0,in=230]  (-1.15, -1.37); 

\node at (-.6,2.1) [color=black!50] {$\mathbf{134}$};
\draw [thick,dotted,->] (-.6, 2.05) to  [out=-90,in=110] (-.52,1.78); 

\node at (-.16,2.1) [color=black!50] {$\mathbf{14}$};
\draw [thick,dotted,->] (-.16, 2.05) to [out=-90,in=90]  (-.3,1.85); 

\draw [dotted]  (-1.908111787, -.5319769439)--(-1.196979681, -1.314983596)--(-1.444193277, -.1924777098) -- (-1.908111787, -.5319769439);

\draw [thick] (-1.908111787, -.5319769439)--(-1.196979681, -1.314983596);

\draw [dotted] (-1.196979681, -1.314983596)--(-.4925605311, -1.935210443)--(-0.0286420207, -1.595711209) --(-1.196979681, -1.314983596);

\draw [thick] (-1.196979681, -1.314983596)--(-.4925605311, -1.935210443);

\draw  [thick] (-.4925605311, -1.935210443)--(-.7683376602, -1.535525063)--(.3999999999, -1.816252676) -- (-.4925605311, -1.935210443);

\draw  [thick] (-.7683376602, -1.535525063)--(-1.643196449, -.6682790668)--(-.7506359181, -.5493212998) --  (-.7683376602, -1.535525063);

\draw  [dotted] (.7683376602, 1.535525063)--(.4925605311, 1.935210443)--(-.3999999999, 1.816252676) --(.7683376602, 1.535525063);

\draw [thick] (.4925605311, 1.935210443)--(-.3999999999, 1.816252676) ;

\draw  [dotted] (1.643196449, .6682790668)--(.7683376602, 1.535525063)--(.7506359181, .5493212998)--(1.643196449, .6682790668);

\draw [dotted] (1.908111787, .5319769439)--(1.643196449, .6682790668)--(1.890410045, -.4542268192)-- (1.908111787, .5319769439);

\draw [thick] (1.890410045, -.4542268192)-- (1.908111787, .5319769439);

\draw [ thick]  (1.908111787, .5319769439)--(1.444193277, .1924777098)--(1.196979681, 1.314983596) -- (1.908111787, .5319769439);

\draw [thick]  (.4925605311, 1.935210443)--(1.196979681, 1.314983596)--(0.0286420207, 1.595711209) -- (.4925605311, 1.935210443);

\draw [dotted]  (-1.139774127, 1.003548119)--(-.3999999999, 1.816252676)--(-1.415551256, 1.403233499) -- (-1.139774127, 1.003548119);

\draw [thick] (-.3999999999, 1.816252676)--(-1.415551256, 1.403233499);

\draw [dotted] (-1.444193277, -.1924777098)--(-.4286420207, .2205414672)--(-1.139774127, 1.003548119)-- (-1.444193277, -.1924777098);

\draw [dotted]  (.7506359181, .5493212998)--(-.4286420207, .2205414672)--(.4462167681, -.6467045291) --(.7506359181, .5493212998);

\draw  [dotted] (-0.286420207e-1, -1.595711209)--(1.150635918, -1.266931376)--(.4462167681, -.6467045291) --(-0.286420207e-1, -1.595711209);

\draw  [dotted] (1.415551256, -1.403233499)--(1.890410045, -.4542268192)--(1.150635918, -1.266931376) -- (1.415551256, -1.403233499);

\draw [thick] (1.415551256, -1.403233499)--(1.890410045, -.4542268192);

\draw  [thick] (1.139774127, -1.003548119)--(1.415551256, -1.403233499)--(.3999999999, -1.816252676) -- (1.139774127, -1.003548119);

\draw  [thick] (1.444193277, .1924777098)--(1.139774127, -1.003548119)--(.4286420207, -.2205414672) --  (1.444193277, .1924777098);

\draw [thick]  (.4286420207, -.2205414672)--(-.7506359181, -.5493212998)--(-.4462167681, .6467045291) -- (.4286420207, -.2205414672);

\draw [ thick]  (0.286420207e-1, 1.595711209)--(-.4462167681, .6467045291)--(-1.150635918, 1.266931376) -- (0.286420207e-1, 1.595711209);

\draw [ thick] (-1.890410045, .4542268192) -- (-1.415551256, 1.403233499)--(-1.150635918, 1.266931376) --(-1.890410045, .4542268192);

\end{tikzpicture}

\caption{The connection graph for $sR(4)$ and the icosidodecahedron. The graph is a folded icosidodecahedron, which is the quotient of the 
 icosidodecahedron by its central symmetry, \textit{i.e.} the vertices and the edges which are related by the central symmetry are identified.
The labels of the same vertex and displayed inside the triangle correspond to the indices of 
the generators belonging to the same abelian subalgebras.
For clarity, we draw the labels in the backside of the picture in gray, with dotted arrows. Moreover, most of them
are not displayed: they can be deduced using the central symmetry. The labels inside the pentagonal faces correspond to central elements.
  \label{fig}}
\end{center}
\end{figure}

To study the representations of the special Racah algebra, it is useful to consider its abelian subalgebras.
It is easy to see that couples of the form $(C_{ij},C_{ijk})$ or $(C_{ij},C_{k\ell})$ (for $i,j,k,\ell$ pairwise distinct) are abelian subalgebras of $sR(4)$.
Hence, both elements of these abelian subalgebras are simultaneously diagonalizable, which fixes a basis of the representation in the case where their joint spectrum is nondegenerate. 
We are then interested in computing the connection coefficients between two different bases. This study is highly facilitated by introducing the associated 
connection graph displayed in figure \ref{fig} and corresponds to a folded icosidodecahedron, which is the quotient of the 
 icosidodecahedron by its central symmetry (\textit{i.e.} the vertices and the edges related by the symmetry with respect to the center of the icosidodecahedron are identified).
The vertices of this graph are given by the abelian subalgebras and the edges link 
two abelian algebras differing by only one generator.
The web of the abelian subalgebras given in figure \ref{fig} generalizes the result of \cite{DGVV1, DIVV} 
where only the subalgebras of the type $(C_{ij},C_{ijk})$ have been considered.
In this latter case, the connection graph reduces to a truncated tetrahedron.

The labels inside a chosen triangle correspond to the noncentral generators of a subalgebra $sR(3)$.
The labels of its three vertices that are outside the triangle
are the central elements of this $sR(3)$ subalgebra.
The labels in the pentagons correspond to the central elements $C_j$ (with 0 associated to $C_{1234}$). 
They are written in such a way that they allow to fully reconstruct isomorphisms of $sR(4)$, see  section \ref{sect:sym-icosi} below.

\subsection{Symmetries of the icosidodecahedron and isomorphisms of algebra\label{sect:sym-icosi}}

An important point to remark is that the previous graph is a regular graph, meaning that all vertices are equivalent. 
Indeed, each vertex is surrounded by two triangles and two pentagons. 
In addition, we remark that the automorphism group of this graph is the permutation group of five elements $\fS_5$ (it has been verified 
using a formal mathematical software).
This rises the question of the link between symmetries of the connection graph and isomorphisms of $sR(4)$.

\begin{prop}\label{eq:propiso}
The linear maps given by:
 \begin{alignat}{2}
 \label{eq:geners}
 \mathfrak{s}:   &\quad C_{12},C_{13},C_{14},C_{23},C_{24},C_{34} &\quad& \mapsto \quad  C_{123},C_{124},C_{134},C_{34},C_{24},C_{23} \\
  &\quad C_{123},C_{124},C_{134},C_{234} &\quad& \mapsto \quad C_{12},C_{13},C_{14},C_{234} \nonumber \\
  &\quad C_1,C_2,C_3,C_4,C_{1234}    &\quad&\mapsto \quad  C_{1234}, C_4, C_3, C_2,C_1\,, \nonumber\\[6pt]
 \label{eq:genert}
  \mathfrak{t}:  &\quad  C_{12},C_{13},C_{14},C_{23},C_{24},C_{34}&\quad&\mapsto\quad C_{13},C_{23},C_{34},C_{12},C_{14},C_{24} \\
   &\quad C_{123},C_{124},C_{134},C_{234} &\quad&\mapsto\quad C_{123},C_{134},C_{234},C_{124}\nonumber \\
  &\quad C_1,C_2,C_3,C_4,C_{1234}   &\quad&\mapsto \quad C_{3}, C_1, C_2, C_4,C_{1234}\,, \nonumber\\[6pt]
 \label{eq:genern}
  \mathfrak{i}:  &\quad  C_{12},C_{13},C_{14},C_{23},C_{24},C_{34}&\quad&\mapsto\quad C_{12},C_{13},C_{234},C_{23},C_{134},C_{124} \\
   &\quad C_{123},C_{124},C_{134},C_{234} &\quad&\mapsto\quad C_{123},C_{34},C_{24},C_{14}\nonumber \\
  &\quad C_1,C_2,C_3,C_4,C_{1234}   &\quad&\mapsto \quad C_{1}, C_2, C_3, C_{1234},C_4\,, \nonumber
  \end{alignat}
are isomorphisms of $sR(4)$. They obey the relations 
\begin{equation}\label{eq:A5}
\begin{aligned}
&  \mathfrak{s}^2=\fe\ ,\quad \mathfrak{t}^3=\fe\ ,\quad \fii^2=\fe\,,\\
&(\mathfrak{s}\,\mathfrak{t})^5=\fe\, ,\quad (\mathfrak{s}\,\mathfrak{i})^4=\fe\, ,\quad (\mathfrak{s}\,\ft\,\mathfrak{i})^6=\fe\, ,\quad 
\ft\,\fii=\fii\,\ft \,,\quad (\mathfrak{i}\,\mathfrak{s}\,\ft\,\fs)^2 = \fe\,,
\end{aligned}
\end{equation} 
where $\fe$ is the identity map. They generate the symmetric group $\fS_5$. 
In addition, $\ft$ and $\fs$ generate the alternating group $\cA_5$. 
\end{prop}
\proof To prove that $\fii$, $\mathfrak{s}$ and  $\mathfrak{t}$ are algebra homomorphisms, one first checks by direct computations
that their action leave the relations of $sR(4)$ unchanged. 
In addition, they verify the relations 
\eqref{eq:A5}. The first row of these relations show that they are invertible, so that they are bijections.
They show that $\ft$ and $\fs$ generate the alternating group $\cA_5$ since $\mathfrak{s}^2=\fe$ ,$\mathfrak{t}^3=\fe$ and $(\mathfrak{s}\,\mathfrak{t})^5=\fe$  are the usual 
defining relations of $\cA_5$. 
To prove that we have the symmetric group $\fS_5$, we introduce
\begin{equation}\label{eq:S5}
\fh_{1}= \fs\,\fii\,\fs\,, \quad \fh_{2}= \fs\,\fii\,\fs\,\ft\,, \quad \fh_{3}= \fii\,\fs\,\ft\,\fs\,, \quad \fh_{4}= \fii\,,
\end{equation} 
which can be inverted as 
\begin{equation}\label{eq:S5inv}
 \fii=\fh_{4}\,, \quad \ft= \fh_{1}\,\fh_{2}\,, \quad \fs=  \big(\fh_{2}\,\fh_{1}\, \fh_{3}\,\fh_{4}\big)^2\, \fh_{2}\, \fh_{1}\,.
\end{equation} 
It is easy to see using relations \eqref{eq:A5} that $\fh_j$, $j=1,2,3,4$ obey the defining relations of $\fS_5$:
\begin{equation}\label{eq:S5bis}
\fh_j^2=\fe\,,\ j=1,2,3,4 \quad \mbox{and}\quad
\begin{cases}
\fh_{j}\,\fh_{k}\,\fh_{j}=\fh_{k}\,\fh_{j}\,\fh_{k}\,\,, \quad |j-k|=1\,,\\
\fh_{j}\,\fh_{k}=\fh_{k}\,\fh_{j}\,, \quad |j-k|\neq1\,.
\end{cases}
\end{equation} 
We checked that the group generated by the generators $\fii$, $\ft$ and $\fs$ has dimension at least 120. This proves that it is indeed 
$\fS_5$.
\endproof

The actions of $\mathfrak{s}$ and $\mathfrak{t}$ \eqref{eq:geners} and \eqref{eq:genert} on the generators $C_{ij}$ and $C_{ijk}$ can be obtained directly
from figure \ref{fig}:
the map $\mathfrak{s}$ corresponds to the reflection with respect to the plane containing the center of the icosidodecahedron and 
the vertices $(C_{12},C_{123})$, $(C_{24},C_{234})$;  
the map $\mathfrak{t}$ is the rotation of $2\pi/3$ around the axis passing through the center 
of the icosidodecahedron and the center of the triangle $(C_{12},C_{123})$, $(C_{13},C_{123})$, $(C_{23},C_{123})$.
 Thanks to this identification, we can obtain the previous isomorphisms generating $\cA_5$ directly from the symmetries of the figure \ref{fig}.
Acting with a symmetry on figure \ref{fig} and reading how the labels are transformed, one gets an isomorphism of $sR(4)$.
For example, the rotation of the figure with an angle $2\pi/5$ around the axis passing through the center 
of the icosidodecahedron and the center of the pentagon $(C_{12},C_{123})$, $(C_{23},C_{123})$, $(C_{234},C_{23})$, $(C_{34},C_{234})$, $(C_{12},C_{34})$
is a symmetry corresponding to the isomorphism
\begin{alignat}{2}
 \label{eq:gener}
 \mathfrak{r}:   &\quad C_{12},C_{13},C_{14},C_{23},C_{24},C_{34} &\quad& \mapsto \quad  C_{34},C_{124},C_{13},C_{123},C_{14},C_{234} \\
  &\quad C_{123},C_{124},C_{134},C_{234} &\quad& \mapsto \quad C_{12},C_{134},C_{24},C_{23} \nonumber \\
  &\quad C_1,C_2,C_3,C_4,C_{1234}    &\quad&\mapsto \quad  C_3, C_4, C_{1234}, C_1,C_2\,. \nonumber
\end{alignat}
\\
 The map $\fii$ is much harder to visualize using the icosidodecahedron image: for example, it  sends the pentagon
 $(C_{12},C_{123})$, $(C_{13},C_{123})$, $(C_{13},C_{24})$, $(C_{124},C_{34})$, $(C_{12},C_{124})$ to the pentagon
 $(C_{12},C_{123})$, $(C_{13},C_{123})$, $(C_{134},C_{13})$, $(C_{134},C_{34})$, $(C_{12},C_{34})$, which is indeed a pentagon of the folded 
 icosidodecahedron (but not of the icosidodecahedron). The  map is illustrated in the following figure:
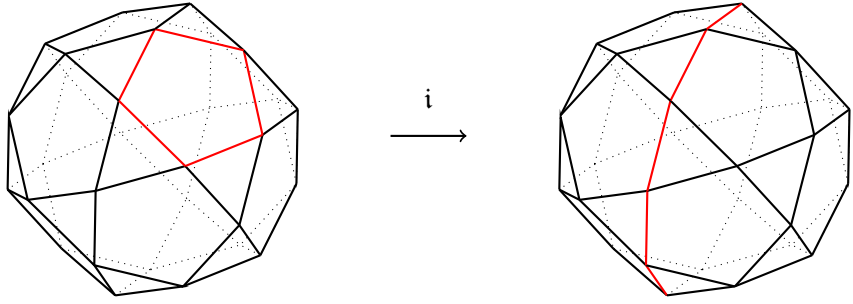
\begin{center}
\begin{tabular}{ccc}
\qquad
\begin{tikzpicture}[scale=1]
\draw [thick] (-1.908111787, -.5319769439)--(-1.890410045, .4242268192)-- (-1.643196449, -.6682790668) -- (-1.908111787, -.5319769439);

\draw [dotted]  (-1.908111787, -.5319769439)--(-1.196979681, -1.314983596)--(-1.444193277, -.1924777098) -- (-1.908111787, -.5319769439);

\draw [thick] (-1.908111787, -.5319769439)--(-1.196979681, -1.314983596);

\draw [dotted] (-1.196979681, -1.314983596)--(-.4925605311, -1.935210443)--(-0.0286420207, -1.595711209) --(-1.196979681, -1.314983596);

\draw [thick] (-1.196979681, -1.314983596)--(-.4925605311, -1.935210443);

\draw  [thick] (-.4925605311, -1.935210443)--(-.7683376602, -1.535525063)--(.3999999999, -1.816252676) -- (-.4925605311, -1.935210443);

\draw  [thick] (-.7683376602, -1.535525063)--(-1.643196449, -.6682790668)--(-.7506359181, -.5493212998) --  (-.7683376602, -1.535525063);

\draw  [dotted] (.7683376602, 1.535525063)--(.4925605311, 1.935210443)--(-.3999999999, 1.816252676) --(.7683376602, 1.535525063);

\draw [thick] (.4925605311, 1.935210443)--(-.3999999999, 1.816252676) ;

\draw  [dotted] (1.643196449, .6682790668)--(.7683376602, 1.535525063)--(.7506359181, .5493212998)--(1.643196449, .6682790668);

\draw [dotted] (1.908111787, .5319769439)--(1.643196449, .6682790668)--(1.890410045, -.4542268192)-- (1.908111787, .5319769439);

\draw [thick] (1.890410045, -.4542268192)-- (1.908111787, .5319769439);

\draw [ thick]  (1.908111787, .5319769439)--(1.444193277, .1924777098);
\draw [ thick, color=red]  (1.444193277, .1924777098)--(1.196979681, 1.314983596);
\draw [ thick]  (1.196979681, 1.314983596) -- (1.908111787, .5319769439);

\draw [thick]  (.4925605311, 1.935210443)--(1.196979681, 1.314983596);
\draw [thick, color=red]  (1.196979681, 1.314983596)--(0.0286420207, 1.595711209);
\draw [thick]  (0.0286420207, 1.595711209) -- (.4925605311, 1.935210443);

\draw [dotted]  (-1.139774127, 1.003548119)--(-.3999999999, 1.816252676)--(-1.415551256, 1.403233499) -- (-1.139774127, 1.003548119);

\draw [thick] (-.3999999999, 1.816252676)--(-1.415551256, 1.403233499);

\draw [dotted] (-1.444193277, -.1924777098)--(-.4286420207, .2205414672)--(-1.139774127, 1.003548119)-- (-1.444193277, -.1924777098);

\draw [dotted]  (.7506359181, .5493212998)--(-.4286420207, .2205414672)--(.4462167681, -.6467045291) --(.7506359181, .5493212998);

\draw  [dotted] (-0.286420207e-1, -1.595711209)--(1.150635918, -1.266931376)--(.4462167681, -.6467045291) --(-0.286420207e-1, -1.595711209);

\draw  [dotted] (1.415551256, -1.403233499)--(1.890410045, -.4542268192)--(1.150635918, -1.266931376) -- (1.415551256, -1.403233499);

\draw [thick] (1.415551256, -1.403233499)--(1.890410045, -.4542268192);

\draw  [thick] (1.139774127, -1.003548119)--(1.415551256, -1.403233499);
\draw  [thick] (1.415551256, -1.403233499)--(.3999999999, -1.816252676);
\draw  [thick] (.3999999999, -1.816252676) -- (1.139774127, -1.003548119);

\draw  [thick] (1.444193277, .1924777098)--(1.139774127, -1.003548119);
\draw  [thick] (1.139774127, -1.003548119)--(.4286420207, -.2205414672) ;
\draw  [thick, color=red] (.4286420207, -.2205414672) --  (1.444193277, .1924777098);

\draw [thick]  (.4286420207, -.2205414672)--(-.7506359181, -.5493212998);
\draw [thick] (-.7506359181, -.5493212998)--(-.4462167681, .6467045291) ;
\draw [thick, color=red]  (-.4462167681, .6467045291) -- (.4286420207, -.2205414672);

\draw [ thick, color=red]  (0.286420207e-1, 1.595711209)--(-.4462167681, .6467045291);
\draw [ thick]  (-.4462167681, .6467045291)--(-1.150635918, 1.266931376);
\draw [ thick]  (-1.150635918, 1.266931376) -- (0.286420207e-1, 1.595711209);

\draw [ thick] (-1.890410045, .4542268192) -- (-1.415551256, 1.403233499);
\draw [ thick] (-1.415551256, 1.403233499)--(-1.150635918, 1.266931376);
\draw [ thick] (-1.150635918, 1.266931376) --(-1.890410045, .4542268192);

\end{tikzpicture}\qquad
&\qquad
\begin{tikzpicture}[scale=1]
\draw [thick, ->] (-0.5, 2)--(0.5, 2);
\node at (0,0) {};
\node at (0,2.5) {$\fii$};
\end{tikzpicture}\qquad
&\qquad
\begin{tikzpicture}[scale=1]
\draw [thick] (-1.908111787, -.5319769439)--(-1.890410045, .4242268192)-- (-1.643196449, -.6682790668) -- (-1.908111787, -.5319769439);

\draw [dotted]  (-1.908111787, -.5319769439)--(-1.196979681, -1.314983596)--(-1.444193277, -.1924777098) -- (-1.908111787, -.5319769439);

\draw [thick] (-1.908111787, -.5319769439)--(-1.196979681, -1.314983596);

\draw [dotted] (-1.196979681, -1.314983596)--(-.4925605311, -1.935210443)--(-0.0286420207, -1.595711209) --(-1.196979681, -1.314983596);

\draw [thick] (-1.196979681, -1.314983596)--(-.4925605311, -1.935210443);

\draw  [thick, color=red] (-.4925605311, -1.935210443)--(-.7683376602, -1.535525063);
\draw  [thick] (-.7683376602, -1.535525063)--(.3999999999, -1.816252676) ;
\draw  [thick] (.3999999999, -1.816252676) -- (-.4925605311, -1.935210443);

\draw  [thick] (-.7683376602, -1.535525063)--(-1.643196449, -.6682790668);
\draw  [thick] (-1.643196449, -.6682790668)--(-.7506359181, -.5493212998);
\draw  [thick, color=red] (-.7506359181, -.5493212998) --  (-.7683376602, -1.535525063);

\draw  [dotted] (.7683376602, 1.535525063)--(.4925605311, 1.935210443)--(-.3999999999, 1.816252676) --(.7683376602, 1.535525063);

\draw [thick] (.4925605311, 1.935210443)--(-.3999999999, 1.816252676) ;

\draw  [dotted] (1.643196449, .6682790668)--(.7683376602, 1.535525063)--(.7506359181, .5493212998)--(1.643196449, .6682790668);

\draw [dotted] (1.908111787, .5319769439)--(1.643196449, .6682790668)--(1.890410045, -.4542268192)-- (1.908111787, .5319769439);

\draw [thick] (1.890410045, -.4542268192)-- (1.908111787, .5319769439);

\draw [ thick]  (1.908111787, .5319769439)--(1.444193277, .1924777098);
\draw [ thick]  (1.444193277, .1924777098)--(1.196979681, 1.314983596);
\draw [ thick]  (1.196979681, 1.314983596) -- (1.908111787, .5319769439);

\draw [thick]  (.4925605311, 1.935210443)--(1.196979681, 1.314983596);
\draw [thick]  (1.196979681, 1.314983596)--(0.0286420207, 1.595711209) ;
\draw [thick, color=red]  (0.0286420207, 1.595711209) -- (.4925605311, 1.935210443);

\draw [dotted]  (-1.139774127, 1.003548119)--(-.3999999999, 1.816252676)--(-1.415551256, 1.403233499) -- (-1.139774127, 1.003548119);

\draw [thick] (-.3999999999, 1.816252676)--(-1.415551256, 1.403233499);

\draw [dotted] (-1.444193277, -.1924777098)--(-.4286420207, .2205414672)--(-1.139774127, 1.003548119)-- (-1.444193277, -.1924777098);

\draw [dotted]  (.7506359181, .5493212998)--(-.4286420207, .2205414672)--(.4462167681, -.6467045291) --(.7506359181, .5493212998);

\draw  [dotted] (-0.286420207e-1, -1.595711209)--(1.150635918, -1.266931376)--(.4462167681, -.6467045291) --(-0.286420207e-1, -1.595711209);

\draw  [dotted] (1.415551256, -1.403233499)--(1.890410045, -.4542268192)--(1.150635918, -1.266931376) -- (1.415551256, -1.403233499);

\draw [thick] (1.415551256, -1.403233499)--(1.890410045, -.4542268192);

\draw  [thick] (1.139774127, -1.003548119)--(1.415551256, -1.403233499);
\draw  [thick] (1.415551256, -1.403233499)--(.3999999999, -1.816252676);
\draw  [thick] (.3999999999, -1.816252676) -- (1.139774127, -1.003548119);

\draw  [thick] (1.444193277, .1924777098)--(1.139774127, -1.003548119);
\draw  [thick] (1.139774127, -1.003548119)--(.4286420207, -.2205414672) ;
\draw  [thick] (.4286420207, -.2205414672) --  (1.444193277, .1924777098);

\draw [thick]  (.4286420207, -.2205414672)--(-.7506359181, -.5493212998);
\draw [thick, color=red] (-.7506359181, -.5493212998)--(-.4462167681, .6467045291) ;
\draw [thick]  (-.4462167681, .6467045291) -- (.4286420207, -.2205414672);

\draw [ thick, color=red]  (0.286420207e-1, 1.595711209)--(-.4462167681, .6467045291);
\draw [ thick]  (-.4462167681, .6467045291)--(-1.150635918, 1.266931376);
\draw [ thick]  (-1.150635918, 1.266931376) -- (0.286420207e-1, 1.595711209);

\draw [ thick] (-1.890410045, .4542268192) -- (-1.415551256, 1.403233499);
\draw [ thick] (-1.415551256, 1.403233499)--(-1.150635918, 1.266931376);
\draw [ thick] (-1.150635918, 1.266931376) --(-1.890410045, .4542268192);

\end{tikzpicture}
\end{tabular}
\end{center}

\section{Representation theory of $sR(4)$ \label{sec:rep}}

In this section, the representations of the rank 2 special Racah algebra $sR(4)$ are computed.
We focus on the ones in which the central elements $C_i$, $i=1,2,3,4$ (resp. $C_{1234}$) take the values $\mu^{(i)}$ (resp. $\mu^{(0)}$), and
such that $C_{12}$ and $C_{123}$ are both diagonal with nondegenerate joint spectrum $\mu_n^{(12)}$ and $\mu_p^{(123)}$, $n,p\in\bZ$.
In the following, it will be  convenient to set 
\begin{equation}\label{eq:def_mu}
\mu^{(i)} = j^{(i)}(j^{(i)}+1)\,,\ i=0,1,\dots,4\,,\quad \mu_n^{(12)} = j_n^{(12)}(j_n^{(12)}+1)
\quad\text{and}\quad \mu_p^{(123)} = j_p^{(123)}(j_p^{(123)}+1).
\end{equation}
  
The basis vectors will be denoted as $\ket{n,p}$, where $n,p\in\bZ$.
The generators $C_i$, $i=1,2,3,4$, $C_{1234}$, $C_{12}$ and $C_{123}$  act as
\begin{eqnarray}
&& C_{i} \, \ket{n,p} = \mu^{(i)}\ket{n,p}\ , \quad  \quad C_{1234} \, \ket{n,p} = \mu^{(0)}\ket{n,p}\ ,\\
\label{eq:subalg12-123}
&&C_{12} \, \ket{n,p} = \mu_n^{(12)} \, \ket{n,p}\ , \quad  \quad
C_{123} \, \ket{n,p} = \mu_p^{(123)} \, \ket{n,p}\,.
\end{eqnarray}
We suppose that the basis is orthonormal and we denote the components of an element $X \in sR(4)$ as follows
\begin{equation}
[X]\ind{n}{m}{p}{q}= \bra{n,p} X \ket{m,q} \, .
\end{equation}

\subsection{Special Racah algebra $sR(3)$}

As explained previously the subalgebra of $sR(4)$ generated by $C_I$ with $I \subseteq \{1,2,3\}$ is the special Racah algebra $sR(3)$.
The representation theory of $sR(3)$ has been studied previously in \cite{GVZ2,HB}. These results are reproduced here for completeness. 

\begin{prop}
\label{prop:sR3}
The diagonal entries of $C_{12}$ are given by \eqref{eq:def_mu} with 
\begin{equation}\label{eq:j12}
j_n^{(12)}=n+a_{12}\,, 
\end{equation}
where $n\in\bZ$ and $a_{12}$ is a free parameter.
The generator $C_{23}$ takes the following tridiagonal form
\begin{align}
\label{eq:reprC23}
&C_{23}\ket{n,p} =  \psi^{+0}_{n+1,p} \ket{n+1,p} +{\psi}^{00}_{n,p}\ket{n,p}+\psi^{-0}_{n,p}\ket{n-1,p}\, , 
\end{align}
where 
\begin{eqnarray}\label{eq:Y}
&&\psi^{+0}_{n,p}\ \psi^{-0}_{n,p} = 
\mathcal{Q}_p(j_n^{(12)})\mathcal{Q}_p(-j_n^{(12)})\,,
\end{eqnarray}
and
\begin{eqnarray}\label{eq:PP}
\mathcal{Q}_p(z) &=& \frac{(z-j^{(1)}+j^{(2)})(z+j^{(1)}+j^{(2)}+1)(z-j^{(3)}+j_p^{(123)})(z-j^{(3)}-j_p^{(123)}-1)}{2z(2z-1)} \,.
\end{eqnarray}
Its diagonal entries are given by
\begin{eqnarray}
\label{eq:pinp}
&&\psi^{00}_{n,p}= -\mathcal{Q}_p(j_{n+1}^{(12)})-\mathcal{Q}_p(-j_n^{(12)})  +(j^{(2)}-j^{(3)}+1)(j^{(2)}-j^{(3)}) \,.
\end{eqnarray}
\end{prop}

\proof   
We use the presentation of appendix \ref{app:algR4}.
The commutator $[C_{23},C_{123}]=0$ leads to the relation
$\big( \mu_p^{(123)}-\mu_q^{(123)} \big) [C_{23}]\ind{n}{m}{p}{q} = 0$, from which it follows that $[C_{23}]\ind{n}{m}{p}{q} = 0$ if $p \ne q$. Then, we can set $p=q$ when considering equations 
\eqref{eq:Rac1} and  \eqref{eq:Rac2}.

Projecting \eqref{eq:Rac1}, one gets
\begin{align}
\label{eq:projRac1}
&\frac{1}{2} \big( \mu_n^{(12)}-\mu_m^{(12)} \big)^2 \,[C_{23}]\ind{n}{m}{p}{p} = \big( \mu_n^{(12)}+\mu_m^{(12)} \big) [C_{23}]\ind{n}{m}{p}{p}  \\
&\qquad + \Big( (\mu_n^{(12)})^2 - (\mu^{(1)}+\mu^{(2)}+\mu^{(3)}+\mu_p^{(123)})\mu_n^{(12)} - (\mu^{(1)}-\mu^{(2)})(\mu^{(3)}-\mu_p^{(123)}) \Big) 
\delta_{n}^{m}\,.\nonumber
\end{align}
This last equation reduces to $\frac{1}{2} \big( \mu_n^{(12)}-\mu_m^{(12)} \big)^2 [C_{23}]\ind{n}{m}{p}{p} = \big( \mu_n^{(12)}+\mu_m^{(12)} \big) [C_{23}]\ind{n}{m}{p}{p}$ if $n \ne m$. 
There exists at least an index $p$  for which $[C_{23}]\ind{n}{m}{p}{p} \ne 0$. Hence, there are only two distinct values of $m$ for fixed $n$ that satisfy the equation. 
Since for a given value of $\mu_p^{(123)}$, the $\mu_n^{(12)}$'s are all distinct by hypothesis, one can choose these two values as $m=n \pm 1$. 
The  entries $[C_{23}]\ind{n}{m}{p}{q}$ with $|n-m| > 1$ are then vanishing.
The recurrence formula obeyed by $\mu_n^{(12)}$ is obtained by considering $m=n+1$ in \eqref{eq:projRac1}:
\begin{equation}
\label{eq:recur}
\frac{1}{2} \big( \mu_n^{(12)}-\mu_{n+1}^{(12)} \big)^2 = \big( \mu_n^{(12)}+\mu_{n+1}^{(12)} \big) \,. 
\end{equation}
Moreover, the projection $\bra{n+2,p}\eqref{eq:Rac2}\ket{n,p}$ gives
\begin{equation}
\mu_{n+2}^{(12)} - 2 \mu_{n+1}^{(12)} + \mu_n^{(12)} = 2 \,,
\end{equation}
which, together with \eqref{eq:recur}, is easily solved and provides the form of the diagonal entries of $C_{12}$ given in the proposition. 

Coming back to \eqref{eq:projRac1} with $n=m$, we get 
\begin{eqnarray}
\label{eq:pinpp}
&&\psi^{00}_{n,p}= \frac{(\mu_n^{(12)}+\mu^{(2)}-\mu^{(1)})(\mu_p^{(123)}-\mu^{(3)}-\mu_n^{(12)})}{2\mu_n^{(12)}} + \mu^{(2)}+\mu^{(3)} \,,
\end{eqnarray}
which can be written as \eqref{eq:pinp} using \eqref{eq:PP}.

The projection $\bra{n,p}\eqref{eq:Rac2}\ket{n,p}$ leads to:
\begin{equation}
\label{eq:Ynp1}
\begin{split}
&\big( \mu_n^{(12)}-\mu_{n+1}^{(12)}-1 \big) Y_{n+1,p}
- \big( \mu_{n-1}^{(12)}-\mu_n^{(12)}+1 \big) Y_{n,p}= \\
&\qquad \big(\psi^{00}_{n,p}\big)^2  + (2\mu_n^{(12)}-\mu^{(1)}-\mu^{(2)}-\mu^{(3)}-\mu_p^{(123)})\, \psi^{00}_{n,p} - (\mu^{(1)}-\mu_p^{(123)})(\mu^{(3)}-\mu^{(2)}) \,,
\end{split}
\end{equation}
where $Y_{n,p}=\psi^{+0}_{n,p}\ \psi^{-0}_{n,p}$. One can check that the other components of \eqref{eq:Rac1} and \eqref{eq:Rac2}  do not give any further information.

In the special Racah algebra, $w_{123} = 0$. The evaluation of this equation on $\bra{n,p}$ and $\ket{n,p}$ leads to:
\begin{equation}
\label{eq:Ynp2}
\begin{split}
&\left( \frac{\big( \mu_n^{(12)}-\mu_{n+1}^{(12)} \big)^2}{4} + \mu_n^{(12)} - 1 \right)  Y_{n+1,p} 
+ \left( \frac{\big( \mu_n^{(12)}-\mu_{n-1}^{(12)} \big)^2}{4} + \mu_n^{(12)} - 1 \right)  Y_{n,p} \\
&= (1-\mu_{n}^{(12)}) \big(\psi^{00}_{n,p}- \mu^{(1)} - \mu^{(2)} - \mu^{(3)} - \mu_p^{(123)}\big) \,\psi^{00}_{n,p}
- (\mu_{n}^{(12)})^2\, \psi^{00}_{n,p}
+ \xi_{12}\, \mu_n^{(12)} + \xi_{23}\, \psi^{00}_{n,p} + \xi_0 \,,
\end{split}
\end{equation}
where
\begin{align}
\xi_{12} &=  (\mu^{(3)}-\mu^{(2)})(\mu^{(1)}-\mu_p^{(123)}) \,, \\
\xi_{23} &=  (\mu^{(1)}-\mu^{(2)})(\mu^{(3)}-\mu_p^{(123)}) \,, \\ 
\xi_0 &= (\mu^{(1)}\mu^{(3)}-\mu^{(2)}\mu_p^{(123)})(\mu_p^{(123)}-\mu^{(1)}+\mu^{(2)}-\mu^{(3)}) + (\mu^{(2)}+\mu^{(3)})(\mu^{(1)}+\mu^{(123)}) \,.
\end{align}
Equations \eqref{eq:Ynp1} and \eqref{eq:Ynp2} form a system of two equations with unknowns 
$ Y_{n,p}$ and $ Y_{n+1,p}$, which is solved by \eqref{eq:Y}.
\endproof

\subsection{Special Racah $sR(4)$}

We now go to the representation of the $sR(4)$ algebra. For the sake of simplicity, we first present  the following lemma which provides
constraints for $C_{34}$ and $C_{123}$.
\begin{lemm}\label{lem1}
\label{prop:sR4}
The diagonal entries of $C_{123}$ are given by \eqref{eq:def_mu} with 
\begin{equation}\label{eq:j123}
j_p^{(123)} = p + a_{123},
\end{equation}
where $p\in\bZ$ and $a_{123}$ is a free parameter. The generator $C_{34}$ must have the following tridiagonal form 
\begin{align}
\label{eq:reprC34}
&C_{34}\ket{n,p} =  \rho^{0+}_{n,p+1} \ket{n,p+1} +{\rho}^{00}_{n,p}\ket{n,p}+\rho^{0-}_{n,p}\ket{n,p-1}\, , 
\end{align}
where 
\begin{eqnarray}\label{eq:Y2}
&&\rho^{0+}_{n,p}\ \rho^{0-}_{n,p}= {\widehat{\mathcal{Q}}_n}(j_p^{(123)}) {\widehat{\mathcal{Q}}_n}(-j_p^{(123)}) \,,
\end{eqnarray}
with
\begin{equation}
{\widehat{\mathcal{Q}}_n}(z) =  \frac{(z-j^{(0)}+j^{(4)})(z+j^{(0)}+j^{(4)}+1)(z-j^{(3)}+j_n^{(12)})(z-j^{(3)}-j_n^{(12)}-1)}{2z\,(2z-1)}\,.
\end{equation}
Its diagonal entries are given by
\begin{equation}\label{eq:r0}
{\rho}^{00}_{n,p} = -  \widehat{\mathcal{Q}}_n(j_{p+1}^{(123)})  - \widehat{\mathcal{Q}}_n(-j_p^{(123)}) +(j^{(3)}-j^{(4)}-1)(j^{(3)}-j^{(4)})\,.
\end{equation}
\end{lemm}
\proof Same proof as in Proposition \ref{prop:sR3}, using now relations \eqref{eq:Rac34-123}, \eqref{eq:Rac123-34} 
and $[C_{12}\,,\,C_{34}]=0$.
\endproof

Up to now, the computations were very similar to the ones done in $sR(3)$. Indeed, in the same way the elements $C_{12}$, 
$C_{23}$, $C_{13}$ with central generators $C_{1}$, 
$C_{2}$, $C_{3}$ and $C_{123}$ generate a $sR(3)$ subalgebra, the elements $C_{123}$, $C_{34}$, $C_{124}$ generate another $sR(3)$ subalgebra containing as central generators the elements $C_{12}$, 
$C_{3}$, $C_{4}$ and $C_{1234}$. The first embedding 
allowed us to get proposition \ref{prop:sR3}, while the second one led to lemma \ref{lem1}.
The remaining calculations are more involved as shown in the following proposition.
\begin{prop}
\label{prop-sR4phi}
The entries of $C_{23}$ and $C_{34}$ are related as follows
\begin{equation}\label{eq:conpsirho}
\frac{ \psi^{+0}_{n,p-1}}{\psi^{+0}_{n,p}} \ \frac{\rho^{0+}_{n,p}}{\rho^{0+}_{n-1,p}}
=\frac{(n-p+a_{12}-a_{123}+j^{(3)}+1)(n-p+a_{12}-a_{123}-j^{(3)})}{(n-p+a_{12}-a_{123}+j^{(3)})(n-p+a_{12}-a_{123}-j^{(3)}-1)}\ .
\end{equation}

The generator $C_{234}$ takes the following form
\begin{equation}
 \label{eq:reprC234}
\begin{aligned}
C_{234}\ket{n,p} =&\  \varphi^{++}_{n+1,p+1} \ket{n+1,p+1} &+&{\varphi}^{0+}_{n,p+1}\ket{n,p+1}&+&\varphi_{n,p+1}^{-+}\ket{n-1,p+1}\\
&+ \varphi_{n+1,p}^{+0} \ket{n+1,p} &+&{\varphi}^{00}_{n,p}\ket{n,p}&+&\varphi^{-0}_{n,p}\ket{n-1,p}  \\
&+ \varphi_{n+1,p}^{+-} \ket{n+1,p-1}& +&{\varphi}^{0-}_{n,p}\ket{n,p-1}&+&\varphi^{--}_{n,p}\ket{n-1,p-1}\,. \\
\end{aligned}
\end{equation}
Its diagonal entries are
\begin{equation}
\label{eq:phidiag}
\varphi^{00}_{n,p}= \frac{\mu_n^{(12)}+\mu^{(1)}-\mu^{(2)}}{2\mu_n^{(12)}} \,\rho^{00}_{n,p}- \frac{1}{2\mu_n^{(12)}} \big( (\mu_n^{(12)}-\mu^{(1)}-\mu^{(2)}-\mu^{(0)})\mu_n^{(12)} + \mu^{(0)}(\mu^{(1)}-\mu^{(2)}) \big)
\end{equation}
while the remaining entries read 
\begin{alignat}{2}
\label{eq:phinmpp}
&\varphi^{\pm 0}_{n,p} = \frac{\mu_p^{(123)}-\mu^{(4)}+\mu^{(0)}}{2\mu_p^{(123)}} \,\psi^{\pm 0}_{n,p}
\,, &\qquad &  \varphi_{n,p}^{0\pm} = \frac{\mu_n^{(12)}+\mu^{(1)}-\mu^{(2)}}{2\mu_n^{(12)}} \,\rho^{0\pm}_{n,p} \,,
\end{alignat}
and 
\begin{alignat}{2}
\label{eq:phis}
&\varphi_{n,p}^{\pm\pm} =\frac{-\psi_{n,p}^{\pm0}\ \rho_{n-1,p}^{0\pm}}{ (n-p+a_{12}-a_{123}+j^{(3)})(n-p+a_{12}-a_{123}-j^{(3)}-1) }\,,\\
\label{eq:phis2}
&\varphi_{n,p}^{\pm\mp} =\frac{-\psi_{n,p}^{\pm 0}\ \rho_{n,p}^{0\mp}}{ (n+p+a_{12}+a_{123}+j^{(3)}+1)(n+p+a_{12}+a_{123}-j^{(3)}) }\,.
\end{alignat}
\end{prop}
\proof
Equation \eqref{eq:Rac123-234} projected on $\bra{n,p}$ and $\ket{m,q}$,  gives
\begin{align}
\label{eq:projRac123-234} 
\frac{1}{2} \big( \mu_p^{(123)}-\mu_q^{(123)} \big)^2 [C_{234}]\ind{n}{m}{p}{q} &= \big( \mu_p^{(123)}+\mu_q^{(123)} \big) [C_{234}]\ind{n}{m}{p}{q} 
- \big( \mu_p^{(123)}-\mu^{(4)}+\mu^{(0)} \big) [C_{23}]\ind{n}{m}{p}{p} \delta_p^q \nonumber \\
&+ \left( (\mu_p^{(123)}-\mu^{(1)}-\mu^{(4)}-\mu^{(0)})\mu_p^{(123)} - \mu^{(1)}(\mu^{(4)}-\mu^{(0)}) \right) \delta_n^m\,\delta_p^q \,.
\end{align}
This expression for $p\ne q$ shows that $[C_{234}]\ind{n}{m}{p}{q}=0$ when $|p-q|>1$, for all values of $m$ and $n$.
Moreover, taking $p=q$ and $m \ne n$, it proves the first expression of \eqref{eq:phinmpp}. \\
In a similar way, Eq.\;\eqref{eq:Rac12-234} leads to
\begin{align}
\label{eq:projRac12-234} 
\frac{1}{2} \big( \mu_n^{(12)}-\mu_m^{(12)} \big)^2 \,[C_{234}]\ind{n}{m}{p}{q} &= \big( \mu_n^{(12)}+\mu_m^{(12)} \big) [C_{234}]\ind{n}{m}{p}{q} - 
\big( \mu_n^{(12)}+\mu^{(1)}-\mu^{(2)} \big) [C_{34}]\ind{n}{n}{p}{q} \delta_n^m \nonumber \\
&+ \left( (\mu_n^{(12)}-\mu^{(1)}-\mu^{(2)}-\mu^{(0)})\mu_n^{(12)} + \mu^{(0)}(\mu^{(1)}-\mu^{(2)}) \right) \delta_n^m\,\delta_p^q \,.
\end{align}
Following the same lines as above, one gets $[C_{234}]\ind{n}{m}{p}{q} = 0$ if $|n-m|>1$ and the second expression of \eqref{eq:phinmpp}, as well as the diagonal entries \eqref{eq:phidiag}.
We have thus shown that $[C_{234}]\ind{n}{m}{p}{q}$ is tridiagonal both in $(n,m)$ and in $(p,q)$. \\
It remains to determine the explicit expressions of $[C_{234}]\ind{n}{m}{p}{q}$ for $m=n \pm 1$ and $q=p \pm 1$. To this aim, we consider \eqref{eq:Rac23-34} projected on $\bra{n,p}$ and $\ket{m,q}$, in the case $q \ne p$,
\begin{eqnarray} \label{eq:prC1}
&&\frac{1}{2} \sum_s\left( [C_{23}]\ind{n}{s}{p}{p}[C_{23}]\ind{s}{m}{p}{p}[C_{34}]\ind{m}{m}{p}{q} - 2[C_{23}]\ind{n}{s}{p}{p}[C_{34}]\ind{s}{s}{p}{q}[C_{23}]\ind{s}{m}{q}{q}
+ [C_{34}]\ind{n}{n}{p}{q}[C_{23}]\ind{n}{s}{q}{q}\,[C_{23}]\ind{s}{m}{q}{q}\right)\nonumber \\
&=&\!\!\!\! [C_{23}]\ind{n}{m}{p}{p}[C_{34}]\ind{m}{m}{p}{q} + [C_{34}]\ind{n}{n}{p}{q}[C_{23}]\ind{n}{m}{q}{q}
-\sum_s [C_{234}]\ind{n}{s}{p}{q}\,[C_{23}]\ind{s}{m}{q}{q}
+ (\mu^{(2)}-\mu^{(3)})\,[C_{234}]\ind{n}{m}{p}{q}.
\end{eqnarray}
The equation obtained for $m=n - 2$ and $q=p - 1$ allows to get: 
\begin{equation}\label{eq:phipp2}
\begin{aligned}
&\frac{1}{2} \Big( \psi^{+0}_{n,p}\psi^{+0}_{n-1,p}\rho^{0+}_{n-2,p}-2\psi^{+0}_{n,p}\rho^{0+}_{n-1,p}\psi^{+0}_{n-1,p-1}+\rho^{0+}_{n,p}\psi^{+0}_{n,p-1}\psi^{+0}_{n-1,p-1} \Big) 
&=& - \varphi^{++}_{n,p}\psi^{+0}_{n-1,p-1} \\
&&=&- \psi^{+0}_{n,p} \varphi^{++}_{n-1,p}\,.
\end{aligned}
\end{equation}
The second equality in the above relations is obtained starting from \eqref{eq:Rac23-34} where $C_{234}C_{23}=C_{23}C_{234}$ is used.
Similarly, for $m=n+ 2$ and $q=p - 1$ in \eqref{eq:prC1}, one gets: 
\begin{equation}\label{eq:phipm2}
\begin{aligned}
&\frac{1}{2} \Big( \psi^{-0}_{n-1,p}\psi^{-0}_{n,p}\rho^{0+}_{n,p}-2\psi^{-0}_{n-1,p}\rho^{0+}_{n-1,p}\psi^{-0}_{n,p-1}+\rho^{0+}_{n-2,p}\psi^{-0}_{n-1,p-1}\psi^{-0}_{n,p-1} \Big) 
&=& - \varphi^{-+}_{n-1,p}\psi^{-0}_{n,p-1} \\
&&=&- \psi^{-0}_{n-1,p} \varphi^{-+}_{n,p}\,.
\end{aligned}
\end{equation}
Relations \eqref{eq:phipp2} provide two different expressions for $\varphi^{++}_{n,p}$. Comparing them, one gets the following constraint
\begin{equation}\label{eq:Zpp}
 1-3Z^{++}_{n-1,p}+3Z^{++}_{n-1,p}Z^{++}_{n,p}-Z^{++}_{n-1,p}Z^{++}_{n,p}Z^{++}_{n+1,p}=0\,, \quad\text{with } \quad 
 Z^{++}_{n,p}=\frac{ \psi^{+0}_{n,p-1}\ \rho^{0+}_{n,p}}{\psi^{+0}_{n,p}\ \rho^{0+}_{n-1,p}} \ .
\end{equation}
Similarly, from \eqref{eq:phipm2}, one gets the following constraint
\begin{equation}\label{eq:Zpm}
 1-3Z^{-+}_{n-1,p}+3Z^{-+}_{n-1,p}Z^{-+}_{n,p}-Z^{-+}_{n-1,p}Z^{-+}_{n,p}Z^{-+}_{n+1,p}=0\,, \quad\text{with } \quad 
 Z^{-+}_{n,p}=\frac{\psi^{-0}_{n,p} \rho^{0+}_{n,p}}{\psi^{-0}_{n,p-1} \rho^{0+}_{n-1,p}} \ .
\end{equation}
From relation \eqref{eq:Y}, one gets 
\begin{eqnarray}\label{eq:Y3}
\frac{ \psi^{-0}_{n,p}}{ \psi^{-0}_{n,p-1} }&=&\frac{ \psi^{+0}_{n,p-1}}{ \psi^{+0}_{n,p}} 
\frac{ (n+p+a_{12}+a_{123}-j^{(3)})  (n-p+a_{12}-a_{123}-j^{(3)}-1)  }
{ (n+p+a_{12}+a_{123}-j^{(3)}-1) (n-p+a_{12}-a_{123}-j^{(3)})   }\nonumber\\
&\times& \frac{  (n-p+a_{12}-a_{123}+j^{(3)}) (n+p+a_{12}+a_{123}+j^{(3)}+1) }
{  (n-p+a_{12}-a_{123}+j^{(3)}+1) (n+p+a_{12}+a_{123}+j^{(3)}) }\,.
\end{eqnarray}
Then, we can express  $Z^{-+}_{n,p}$ in terms of $Z^{++}_{n,p}$ and relations \eqref{eq:Zpp} and \eqref{eq:Zpm} can be solved 
to get an expression of $Z^{++}_{n+1,p}$ and of $Z^{++}_{n-1,p}$ in terms of $Z_{n,p}^{++}$.
Replacing in the first relation $n$ by $n-1$, one gets two relations between $Z^{++}_{n-1,p}$ and $Z_{n,p}^{++}$.
Solving these relations, we find that
\begin{equation}
 Z^{++}_{n,p}=\frac{(n-p+a_{12}-a_{123}+j^{(3)}+1)(n-p+a_{12}-a_{123}-j^{(3)})}{(n-p+a_{12}-a_{123}+j^{(3)})(n-p+a_{12}-a_{123}-j^{(3)}-1)}\,.
\end{equation}
This proves relation \eqref{eq:conpsirho}.
Taking into account this relation, one gets from \eqref{eq:phipp2} the expression \eqref{eq:phis} for $\varphi^{++}_{n,p}$.
The values \eqref{eq:phis}-\eqref{eq:phis2} of $\varphi^{--}_{n+1,p+1}$, $\varphi^{+-}_{n,p+1}$, $\varphi^{-+}_{n+1,p}$ are computed similarly.
\endproof

Let us remark that, from relations \eqref{eq:Y} and \eqref{eq:Y2}, relation \eqref{eq:conpsirho} leads to the following relations
\begin{eqnarray}\label{eq:conpsirho3}
&&\frac{ \psi^{+0}_{n,p-1}}{\psi^{+0}_{n,p}} \ \frac{\rho^{0-}_{n-1,p}}{\rho^{0-}_{n,p}}
=\frac{(n+p+a_{12}+a_{123}+j^{(3)})(n+p+a_{12}+a_{123}-j^{(3)}-1)}{(n+p+a_{12}+a_{123}+j^{(3)}+1)(n+p+a_{12}+a_{123}-j^{(3)})}\,, \\
&&\frac{ \psi^{-0}_{n,p}}{\psi^{-0}_{n,p-1}} \ \frac{\rho^{0+}_{n,p}}{\rho^{0+}_{n-1,p}}
=\frac{(n+p+a_{12}+a_{123}+j^{(3)}+1)(n+p+a_{12}+a_{123}-j^{(3)})}{(n+p+a_{12}+a_{123}+j^{(3)})(n+p+a_{12}+a_{123}-j^{(3)}-1)} \,,\\
&&\frac{ \psi^{-0}_{n,p}}{\psi^{-0}_{n,p-1}} \ \frac{\rho^{0-}_{n-1,p}}{\rho^{0-}_{n,p}}
=\frac{(n-p+a_{12}-a_{123}+j^{(3)})(n-p+a_{12}-a_{123}-j^{(3)}-1)}{(n-p+a_{12}-a_{123}+j^{(3)}+1)(n-p+a_{12}-a_{123}-j^{(3)})}\ .
\end{eqnarray}
Note also that \eqref{eq:conpsirho} can be solved as
\begin{equation}\label{eq:unp}
\begin{aligned}
\psi^{+0}_{n,p}&
=(n-p+a_{12}-a_{123}+j^{(3)})\,\prod_{j=-\infty}^pu_{n,p-j}\ ,\\
\rho^{0+}_{n,p}&
=(n-p+a_{12}-a_{123}-j^{(3)})\,\prod_{k=-\infty}^nu_{n-k,p}\ ,
\end{aligned}
\end{equation}
where $u_{n,p}$ are  arbitrary nonzero numbers. It can be checked that these numbers can be absorbed through a 
conjugation  by a diagonal matrix. 
Note however that they play a nontrivial role when one considers particular bases for the finite-dimensional representations (see Section \ref{sec:finite}).

Up to now, we obtained necessary conditions to get a representation. The following theorem ensures that they are also sufficient. 
\begin{thm}\label{th:repinf} The diagonal matrices
\begin{equation}
 C_{12}\ket{n,p}=(n+a_{12})(n+a_{12}+1)\ket{n,p}\,, \qquad C_{123}\ket{n,p}=(p+a_{123})(p+a_{123}+1)\ket{n,p}\,,
\end{equation}
the tridiagonal matrices
\begin{align}
&C_{23}\ket{n,p} =  \psi^{+0}_{n+1,p} \ket{n+1,p} +{\psi}^{00}_{n,p}\ket{n,p}+\psi^{-0}_{n,p}\ket{n-1,p}\, ,\\
&C_{34}\ket{n,p} =  \rho^{0+}_{n,p+1} \ket{n,p+1} +{\rho}^{00}_{n,p}\ket{n,p}+\rho^{0-}_{n,p}\ket{n,p-1}\, , 
\end{align}
and the matrix
 \begin{equation}
 \label{eq:reprC234th}
\begin{aligned}
C_{234}\ket{n,p} =&  \varphi^{++}_{n+1,p+1} \ket{n+1,p+1} &+&{\varphi}^{0+}_{n,p+1}\ket{n,p+1}&+&\varphi_{n,p+1}^{-+}\ket{n-1,p+1}\\
+& \varphi_{n+1,p}^{+0} \ket{n+1,p} &+&{\varphi}^{00}_{n,p}\ket{n,p}&+&\varphi^{-0}_{n,p}\ket{n-1,p}  \\
+ & \varphi_{n+1,p}^{+-} \ket{n+1,p-1}& +&{\varphi}^{0-}_{n,p}\ket{n,p-1}&+&\varphi^{--}_{n,p}\ket{n-1,p-1}\, , \\
\end{aligned}
\end{equation}
provide a representation of $sR(4)$ if $\psi$ and $\rho$ verify  \eqref{eq:Y}, \eqref{eq:pinp}, \eqref{eq:Y2}, \eqref{eq:r0}, \eqref{eq:conpsirho} and $\varphi$ satisfies 
\eqref{eq:phidiag}-\eqref{eq:phis2}.

Hence we get representations of $sR(4)$ labeled by the multiplet $\big\{j^{(1)}, j^{(2)},j^{(3)}, j^{(4)} ,j^{(0)} ,a_{12},a_{123}\big\}$.
\end{thm}
\proof All the defining relations of $sR(4)$ displayed in appendix \ref{app:algR4} are verified by direct computation.
In addition, one can check that the represented generators indeed obey \eqref{eq:specialR4}.
\endproof

\subsection{Finite-dimensional symmetric irreducible representations \label{sec:finite}}

In the previous section, we computed infinite-dimensional representations of $sR(4)$, labeled by parameters $a_{12}$, $a_{123}$ and $j^{(i)}$ 
($i=0,1,\dots 4$). In this subsection, we provide constraints between these parameters to get finite-dimensional representations.

From now on, we suppose that the parameters $j^{(i)}$ are real and, without any loss of generality, we can choose them nonnegative ($j^{(i)}\geq0$), 
since we have the invariance $j^{(i)}\rightarrow -j^{(i)}-1$. 
\paragraph{Finite-dimensional representations.}
To get finite-dimensional representations, one has to impose that two nondiagonal entries of the tridiagonal generators vanish. 

On the one hand, it is done by fixing the value of the parameters 
$a_{12}$ and $a_{123}$.
Different choices of the parameters $a_{12}$ and $a_{123}$ lead to equivalent representations since they simply correspond to some shifts of $n$ or $p$.
One chooses
\begin{equation}\label{eq:N}
a_{12} = -j^{(1)} - j^{(2)} - 1 \,, \qquad a_{123} = - j^{(4)} - j^{(0)} -1 \,, 
\end{equation}
and one gets $j_n^{(12)}=n-j^{(1)} - j^{(2)} - 1$ and $j_p^{(123)}=p-j^{(4)}-j^{(0)}-1$.

On the other hand, let us define 
\begin{equation}
 N = j^{(1)} + j^{(2)} - j^{(3)} + j^{(4)} + j^{(0)}  \,.
 \end{equation}
In this case, the nondiagonal entries of $C_{23}$ and $C_{34}$, given by \eqref{eq:Y} and \eqref{eq:Y2}, become 
\begin{align}
\psi^{+0}_{n,p}\ \psi^{-0}_{n,p} &= \frac{n(2j^{(1)}+1-n)(2j^{(2)}+1-n)(2j^{(1)}+2j^{(2)}+2-n)(N-n-p+2j^{(3)}+2) }{(2j^{(1)}+2j^{(2)} +2-2n)^2(2j^{(1)}+2j^{(2)}+1-2n) } \nonumber \\
 &\times \frac{ (N+1-n-p)(n-p-2j^{(1)}-2j^{(2)}+2j^{(3)}+N)(p-n-N+2j^{(1)}+2j^{(2)}+1) }{(2j^{(1)}+2j^{(2)}+3-2n) }\label{eq:psif}\,, \\
\rho^{0+}_{n,p}\ \rho^{0-}_{n,p} &= \frac{p(2j^{(4)}+1-p)(2j^{(0)}+1-p)(2j^{(4)}+2j^{(0)}+2-p)(N-n-p+2j^{(3)}+2) }{(2j^{(0)}+2j^{(4)} +2-2p)^2(2j^{(0)}+2j^{(4)}+1-2p)} \nonumber \\
 &\times \frac{ (N+1-n-p)(p-n-2j^{(0)}-2j^{(4)}+2j^{(3)}+N)(n-p-N +2j^{(0)}+2j^{(4)}+1  ) }{(2j^{(0)}+2j^{(4)}+3-2p) }\,. \label{eq:rhof} 
\end{align}
If we suppose that $N$ is a positive integer ($N> 0$), we deduce that
\begin{eqnarray}\label{eq:psi0m1}
&&[C_{23}]\ind{0}{-1}{p}{p} [C_{23}]\ind{-1}{0}{p}{p} = 0 \qquad \text{and} \qquad  [C_{23}]\ind{N-p}{N-p+1}{p}{p} [C_{23}]\ind{N-p+1}{N-p}{p}{p} = 0 \,,\\
&&[C_{34}]\ind{n}{n}{0}{-1}[C_{34}]\ind{n}{n}{-1}{0}=0 \qquad \text{and} \qquad [C_{34}]\ind{n}{n}{N-n}{N-n+1} [C_{34}]\ind{n}{n}{N-n}{N-n+1} =0\,.
\end{eqnarray}
Similar relations hold for the other generators of $sR(4)$. We deduce that the space spanned by the following vectors
\begin{equation}\label{eq:FN}
 \mathcal{F}_N=\text{span}_\mathbb{C}\left\{\ket{n,p} \ | \   n,p\geq 0\,,\ \ n+p \leq N\ \right\}\,,
\end{equation}
provides a finite-dimensional representation of $sR(4)$.
The dimension of this representation is given by
\begin{equation}\label{eq:dimF}
 \text{dim}(\mathcal{F}_N)=\begin{pmatrix} N+2\\2 \end{pmatrix}\,.
\end{equation}

\paragraph{Real symmetric representations.}

To get a symmetric representation, we demand that $\psi^{+0}_{n,p}= \psi^{-0}_{n,p}$ and $\rho^{0+}_{n,p}= \rho^{0-}_{n,p}$, for $n,p\geq 0$ and $n+p \leq N$.
Then to obtain a real representation, it is necessary that the r.h.s. of \eqref{eq:psif} and \eqref{eq:rhof} be positive.
These constraints are satisfied for example if the following inequalities for $j^{(i)}$ hold:
\begin{equation}
\begin{aligned}
& j^{(1)}\geq j^{(2)} \geq j^{(4)} \,,\quad j^{(0)} \geq j^{(4)} \,,\quad  j^{(i)}\geq 0\,,\quad i=0,1,...,4\,,
\\
&j^{(1)}+j^{(2)}+j^{(0)}+j^{(4)} \geq j^{(3)}\geq j^{(1)}+j^{(2)}+j^{(0)}-j^{(4)}\,. 
\end{aligned}
\label{eq:conj}
\end{equation}
From now on, we assume that these conditions are satisfied. A discussion about the other possibilities is postponed to the conclusion of this paper.
Let us emphasize that each factor in the rational functions \eqref{eq:psif} and \eqref{eq:rhof} are positive with the above constraints \eqref{eq:conj}. 
Then, we can take 
\begin{equation}
 \psi^{+0}_{n,p}= \psi^{-0}_{n,p}=: \psi_{n,p}\,,\qquad \rho^{0+}_{n,p}= \rho^{0-}_{n,p}=: \rho_{n,p}\,,
\end{equation}
where  $\psi_{n,p}$ and  $\rho_{n,p}$ are the positive square roots of the r.h.s. of \eqref{eq:psif} and \eqref{eq:rhof} respectively.
We can check by direct computation that relation \eqref{eq:conpsirho} is satisfied in this case. This solution corresponds to a particular choice of the coefficients $u_{n,p}$ in \eqref{eq:unp}.
It follows that $C_{234}$ is also symmetric and that its entries satisfy:
\begin{equation}
\varphi^{++}_{n,p} = \varphi^{--}_{n,p} =: \varphi^{D}_{n,p}, \quad
\varphi^{+-}_{n,p} = \varphi^{-+}_{n,p} =: \varphi^{A}_{n,p}, \quad
\varphi^{0+}_{n,p} = \varphi^{0-}_{n,p} =: \varphi^{V}_{n,p}, \quad
\varphi^{+0}_{n,p} = \varphi^{-0}_{n,p} =: \varphi^{H}_{n,p}.
\end{equation}
Since the explicit expressions of the different functions $\psi_{n,p}, \rho_{n,p}, \varphi_{n,p}$ are rather cumbersome, we postpone them in the Appendix \ref{app:C}.
It follows that all the other generators in this representation are also symmetric.

\paragraph{Irreducible representation.} 

The representation is irreducible.
Indeed, as the spectrum of the abelian subalgebra $(C_{12},C_{123})$ is nondegenerate, 
one can construct projectors on any $\ket{n,p}\in \mathcal{F}_N$.
By acting with these projectors on any nonvanishing vector $v\in \mathcal{F}_N$, one can get a given vector $\ket{n,p}$. 
Then with the action of $C_{23}$ and $C_{34}$  (by remarking that $\psi_{n,p}$ and $\rho_{n,p}$ do not vanish for $n,p\geq 0$ and $n+p \leq N$ if the 
conditions \eqref{eq:conj} are fulfilled), a linear combination of $\ket{n,p}$ and $\ket{n\pm 1,p\pm 1}$ 
with nonvanishing coefficients is obtained. By acting again with the projectors, each vector $\ket{n\pm 1,p\pm 1}$ can be reached. All the vectors of $\mathcal{F}_N$ can be obtained by
iterating this procedure. Therefore the whole space $\mathcal{F}_N$ 
is generated by action of the matrices of the representation on any nonvanishing vector, which proves the irreducibility of the representation.
 
 \paragraph{Summary of the results.} The previous results and discussions lead to the following theorem:

 \begin{thm}\label{th:repfin} Let the vectors $\ket{n,p}$ belong to the vector space $\mathcal{F}_N$ given in \eqref{eq:FN}. The generators of $sR(4)$
 acting on these vectors as
\begin{align}
 &C_{12}\ket{n,p}=(n-j^{(1)}-j^{(2)}-1)(n-j^{(1)}-j^{(2)}+1)\ket{n,p}\,, \\ &C_{123}\ket{n,p}=(p-j^{(0)}-j^{(4)})(p-j^{(0)}-j^{(4)}+1)\ket{n,p}\,,\\
&C_{23}\ket{n,p} =  \psi_{n+1,p} \ket{n+1,p} +{\psi}^{00}_{n,p}\ket{n,p}+\psi_{n,p}\ket{n-1,p}\, ,\\
&C_{34}\ket{n,p} =  \rho_{n,p+1} \ket{n,p+1} +{\rho}^{00}_{n,p}\ket{n,p}+\rho_{n,p}\ket{n,p-1}\, , 
\end{align}
and 
 \begin{equation}
 \label{eq:reprC234the}
\begin{aligned}
C_{234}\ket{n,p} =&  \varphi^{D}_{n+1,p+1} \ket{n+1,p+1} &+&{\varphi}^{V}_{n,p+1}\ket{n,p+1}&+&\varphi_{n,p+1}^{A}\ket{n-1,p+1}\\
+& \varphi_{n+1,p}^{H} \ket{n+1,p} &+&{\varphi}^{00}_{n,p}\ket{n,p}&+&\varphi^{H}_{n,p}\ket{n-1,p}  \\
+ & \varphi_{n+1,p}^{A} \ket{n+1,p-1}& +&{\varphi}^{V}_{n,p}\ket{n,p-1}&+&\varphi^{D}_{n,p}\ket{n-1,p-1}\, , \\
\end{aligned}
\end{equation}
provide a finite-dimensional real symmetric irreducible representation of $sR(4)$ with  $\psi$, $\rho$ and $\varphi$ given in appendix \ref{app:C}.
 This representation is labeled by the multiplet of real parameters $\big\{j^{(1)}, j^{(2)},j^{(3)}, j^{(4)} ,j^{(0)}\big\}$ subject to the constraint
 $N=j^{(1)}+ j^{(2)}-j^{(3)}+ j^{(4)} +j^{(0)}\in \ZZ_{> 0}$ and inequalities \eqref{eq:conj}.
 
 \end{thm}

\section{Equivalent representations and Racah polynomials \label{sec:equiv}}

The objective of this section is to compute different finite-dimensional representations labeled by elements $\fg \in \fS_5$ and to show that they are equivalent.
The overlap coefficients (\textit{i.e.} the entries of the transition matrix between two such representations) are computed exactly and are 
expressed in terms of the Racah polynomials.

\subsection{Construction of other representations \label{sec:repsym}} 

The section \ref{sec:rep} was devoted in constructing one finite-dimensional representation of $sR(4)$ labeled by the multiplet $\big\{j^{(1)}, j^{(2)},j^{(3)}, j^{(4)} ,j^{(0)}\big\}$.
The isomorphisms of $sR(4)$ proved in proposition \ref{eq:propiso}, composed
with this representation, allows us to construct new representations. 
To define these new representations, let us introduce the notation 
\begin{equation}
 \pi^{\boldsymbol{J}}_\fe(X) \qquad \text{for} \quad X\in sR(4),\   \boldsymbol{J}=(j^{(1)},j^{(2)},j^{(3)},j^{(4)},j^{(0)})\,,
\end{equation}
for the finite representation given in theorem \ref{th:repfin} with the constraints given in section \ref{sec:finite}. 
It is indexed by the identity element $\fe$ of the symmetric group $\fS_5$ and $\boldsymbol{J}$ associated to the eigenvalues of the central elements $(C_1,C_2,C_3,C_4,C_{1234})$.
The positive integer $N$, see \eqref{eq:N}, characterizing the dimension of the finite 
representation, is denoted $N({\boldsymbol{J}})$. 
Let us also define the maps\footnote{With a slight abuse of notation, we use the same notation for these functions and for 
the isomorphisms \eqref{eq:geners}, \eqref{eq:genert}, \eqref{eq:genern} to which they correspond.} $\mathfrak{s}$, $\mathfrak{t}$ and $\fii$  on any quintuplet by
\begin{equation}\label{eq:A5-quintu}
\begin{aligned}
&\mathfrak{s}\,(z_1,z_2,z_3,z_4,z_0) = (z_0,z_4,z_3,z_2,z_1) \,,\\
& \mathfrak{t}\,(z_1,z_2,z_3,z_4,z_0) = (z_2,-z_3-1,-z_1-1,z_4,z_0) \,,\\
& \mathfrak{i}\,(z_1,z_2,z_3,z_4,z_0) = (z_1,z_2,z_3,z_0,z_4)\,.
\end{aligned}
\end{equation}
The maps are extended by morphism, and it can be checked that it leads to a representation of $\fS_5$. 
For a given element $\fg\in \fS_5$, one defines a new representation of $sR(4)$ as follows
\begin{equation}\label{eq:defpiG}
  \pi^{\bJ}_\fg(X) =\pi^{\fg\bJ}_\fe(\fg(X)) \,.
\end{equation}

By direct calculation, one can show that for $\boldsymbol{J}=(j^{(1)},j^{(2)},j^{(3)},j^{(4)},j^{(0)})$, we have
\begin{equation}
\begin{aligned}
 \pi^{\boldsymbol{J}}_\ft(C_i)&=  \pi^{(j^{(2)},-j^{(3)}-1,-j^{(1)}-1,j^{(4)},j^{(0)})}_\fe(\ft(C_i)) =\mu^{(i)}\, \,.
\end{aligned}
\end{equation}
By abuse of notation, we do not indicate in the r.h.s. of the previous relation the identity matrix.
Similar relations hold for  $\pi^{\boldsymbol{J}}_\fs(C_i)$ and $ \pi^{\boldsymbol{J}}_\fii(C_i)$. Therefore, by morphism we get for any $\mathfrak{g}\in\fS_5$ 
\begin{equation}
\pi^{\bJ}_\fg(C_i) =\mu^{(i)}\, \,.
\end{equation}
This relation shows that the representation of the central element $C_i$ is the same in any such representation. 
The key point to get this property is the choice of the functions  \eqref{eq:A5-quintu}.

From the definitions \eqref{eq:A5-quintu} of $\mathfrak{s}$, $\mathfrak{t}$ and $\fii$, and the explicit form \eqref{eq:N} of $N(\boldsymbol{J})$, we see that $N(\mathfrak{g}\boldsymbol{J})=N(\boldsymbol{J})$.
Therefore, the representations $\pi^{\boldsymbol{J}}_\fg$ have all the same dimension for $\mathfrak{g}\in\fS_5$. 
Let us also mention that the r.h.s. of \eqref{eq:psif} and \eqref{eq:rhof} with $\bJ$ replaced by $\mathfrak{g}\bJ$, for any $\fg\in\mathfrak{S}_5$, remains positive so that 
$\psi_{n,p}(\mathfrak{g}\bJ)$ and $\rho_{n,p}(\mathfrak{g}\bJ)$ are also square roots of positive quantities. \\

The next subsections are devoted to prove that these representations are in fact equivalent \textit{i.e.} for $\fg,\fh \in \fS_5$ there exists an invertible matrix $T_{\fh,\fg}(\boldsymbol{J})$, called transition matrix, 
such that 
\begin{equation}
T_{\fh,\fg}(\bJ)\ \pi_\fg^{\bJ}(X)=  \pi_\fh^{\bJ}(X)\ T_{\fh,\fg}(\bJ)\,, \qquad \text{for any } X\in sR(4) \,. \label{eq4}
\end{equation}
As the representations are symmetric, the transition matrix $T_{\fh,\fg}$ can be chosen orthogonal: 
 \begin{equation}\label{eq:ortho}
  T_{\fh,\fg}(\bJ)\cdot T^t_{\fh,\fg}(\bJ)=1 \,.
 \end{equation}
The transition matrices  $T_{\fh,\fg}$ are $\begin{pmatrix} N+2\\2 \end{pmatrix} \times \begin{pmatrix} N+2\\2 \end{pmatrix}$ matrices with the entries denoted $[T_{\fh,\fg}]_{np}^{mq}$ with $n,p,m,q\geq 0$ 
and $n+p,m+q\leq N$.

\subsection{Transition matrix between $\pi_{\fe}^{\boldsymbol{J}}$ and $\pi_{\ft}^{\boldsymbol{J}}$ \label{sec:ot}}

Let us start by computing the transition matrix $T_{\fe,\ft}=T_{\fe,\ft}(\bJ)$.
The equation \eqref{eq4} for $X=C_{123}$ leads to
\begin{equation}
T_{\fe,\ft}\ \pi_\fe^{\ft\bJ}(C_{123})=  \pi_\fe^{\bJ}(C_{123})\ T_{\fe,\ft}\,,
\end{equation}
that is
\begin{equation}
\big( \mu_q^{(123)}(\mathfrak{t}\boldsymbol{J})-\mu_p^{(123)}(\boldsymbol{J}) \big) \, [T_{\fe,\ft}]_{np}^{mq} = 0 \,.
\end{equation}
Since $\mu_q^{(123)}(\mathfrak{t}\boldsymbol{J}) = \mu_q^{(123)}(\boldsymbol{J})$, one gets 
$[T_{\fe,\ft}]_{np}^{mq}=\delta_{pq}\ [T_{\fe,\ft}]_{np}^{mp}$.
We get then for $X=C_{23}$ in \eqref{eq4}:
\begin{equation}
T_{\fe,\ft}\ \pi_\fe^{\ft\bJ}(C_{12})=  \pi_\fe^{\bJ}(C_{23})\ T_{\fe,\ft}\,,
\end{equation}
\textit{i.e.} inserting the representation \eqref{eq:reprC23},
\begin{equation}
\label{eq:reccurracah}
\mu_m^{(12)}(\mathfrak{t}\boldsymbol{J} ) \, [T_{\fe,\ft}]_{np}^{mp} = \psi_{n+1,p}(\boldsymbol{J}) \ [T_{\fe,\ft}]_{n+1,p}^{mp}  \ + \psi^{00}_{n,p}(\boldsymbol{J})\ [T_{\fe,\ft}]_{np}^{mp}
+ \psi_{n,p}(\boldsymbol{J}) \ [T_{\fe,\ft}]_{n-1,p}^{mp} \,,
\end{equation}
where $\mu_m^{(12)}(t\boldsymbol{J}) = (m-j^{(2)}+j^{(3)}+1)(m-j^{(2)}+j^{(3)})$. 
Now, \eqref{eq4} for $X=C_{12}$ leads to
\begin{align}
&- \psi_{m+1,p}(\ft\boldsymbol{J})  \,[T_{\fe,\ft}]_{np}^{m+1,p}  
- \psi_{m,p}(\ft\boldsymbol{J})   \, [T_{\fe,\ft}]_{np}^{m-1,p}\nonumber\\
&+ \big( \mu_p^{(123)}+\mu^{(1)}+\mu^{(2)}+\mu^{(3)}-\mu_m^{(12)}-\psi_{m,p}^{00}  \big)(\mathfrak{t}\boldsymbol{J})[T_{\fe,\ft}]_{np}^{mp}= 
\mu_n^{(12)}(\boldsymbol{J}) [T_{\fe,\ft}]_{np}^{mp}\,.
\label{eq:diff}
\end{align}
Given the expressions \eqref{eq:pinp} and \eqref{eq:psifapp} of the coefficients $\psi$, the equations \eqref{eq:reccurracah} and \eqref{eq:diff} 
become the recurrence and difference equations of the Racah polynomials, see \eqref{eq:recuq} and \eqref{eq:diffq}.
Therefore, the transition matrix can be expressed in terms of a function $\cP_n(m)$  proportional to the Racah polynomial, see \eqref{eq:qR}, as follows  
\begin{align}
\label{eq:racah_poly}
 [T_{\fe,\ft}]_{np}^{mp} &=\, (-1)^m \,\sigma_\mathfrak{t}(\boldsymbol{J};p) \, \cP_n(m;-2j^{(2)}-1,-2j^{(1)}-1,p-N-1,N-p-2j^{(2)}+ 2j^{(3)}+1)\,,
\end{align}
where $\sigma_\mathfrak{t}(\boldsymbol{J};p)$ is an overall coefficient which is determined below.

Consider now the case $X=C_{234}$ in \eqref{eq4}:
\begin{equation}
T_{\fe,\ft}\ \pi_\fe^{\ft\bJ}(C_{124})=  \pi_\fe^{\bJ}(C_{234})\ T_{\fe,\ft}\,.
\end{equation}
The component $\bra{n,p}\eqref{eq4}\ket{m,q}$ for $q=p-1$ leads to 
\begin{equation}
\label{eq4.16A}
-\rho_{m,p}(\ft\boldsymbol{J}) \,[T_{\fe,\ft}]_{np}^{mp} = \varphi_{n+1,p}^{A}(\boldsymbol{J}) \ [T_{\fe,\ft}]_{n+1,p-1}^{m,p-1} 
+ \varphi_{n,p}^{V}(\boldsymbol{J}) \,[T_{\fe,\ft}]_{n,p-1}^{m,p-1} + \varphi_{n,p}^{D}(\boldsymbol{J})\,[T_{\fe,\ft}]_{n-1,p-1}^{m,p-1}  \,,
\end{equation}
where we have used relation \eqref{eq:Ctrou} to express $C_{124}$ in terms of $C_{34}$. 
Inserting the expression \eqref{eq:racah_poly}, taking into account \eqref{eq:B4}--\eqref{eq:B11} 
and the explicit expressions of $\varphi$ and $\rho$ (see equation \eqref{eq:phidiag} and appendix \ref{app:C}), 
the previous equation leads to $\sigma_\ft(\boldsymbol{J};p)=-\sigma_\ft(\boldsymbol{J};p-1)$. Hence, one gets $\sigma_\ft(\boldsymbol{J};p)=(-1)^p\,\sigma_\ft(\boldsymbol{J};0)$.
The orthogonality relation of the transition matrix with relation \eqref{eq:ortho3} then implies $\sigma_\ft(\boldsymbol{J};0) = \pm 1$. The component $\bra{n,p}\eqref{eq4}\ket{m,q}$ for $q=p+1$ leads to the same result, using now \eqref{eq:reln2} instead of \eqref{eq:reln1}, while for $q=p$, one recovers the recurrence relation \eqref{eq:reccurracah}. 

The component $\bra{n,p}\eqref{eq4}\ket{m,q}$ with $p=q-1$ for $X=C_{34}$ leads to 
\begin{equation}
\label{eq4.16B}
\varphi_{m+1,q}^{A}(\ft\boldsymbol{J}) \ [T_{\fe,\ft}]_{n,q-1}^{m+1,q-1} 
+ \big( \varphi_{m,q}^{V}-\rho_{m,q} \big)(\ft\boldsymbol{J}) \,[T_{\fe,\ft}]_{n,q-1}^{m,q-1} + \varphi_{m,q}^{D}(\ft\boldsymbol{J})\,[T_{\fe,\ft}]_{n,q-1}^{m-1,q-1} = \rho_{n,q}(\boldsymbol{J}) \,[T_{\fe,\ft}]_{nq}^{mq} \,.
\end{equation}
Using now \eqref{eq:reln3} and \eqref{eq:reln4}, it is checked that no new constraint shows up.
Finally, the relations associated to the remaining values of $X$ are equivalent to \eqref{eq4.16A} or \eqref{eq4.16B}.
In conclusion, the two representations are equivalent and the transition matrix $T_{\fe,\ft}$ takes the following form:
\begin{align}
\label{eq:overlap_t}
 [T_{\fe,\ft}(\boldsymbol{J})]_{np}^{mq}  &\sim \delta_{pq}\, (-1)^{m+p} \, \cP_n(m;-2j^{(2)}-1,-2j^{(1)}-1,p-N-1,N-p-2j^{(2)}+ 2j^{(3)}+1)\,,
\end{align}
where $\sim$ stands for \textit{equal up to a global sign}.
\begin{rema}
The formula \eqref{eq:overlap_t} can be generalized for an arbitrary $\boldsymbol{J}'=\fg\boldsymbol{J}$ where $\fg\in\fS_5$ (generated by $\fii$, $\ft$ and $\fs$, see \eqref{eq:A5-quintu}). Indeed, following the lines of the previous proof, the point which requires special attention is the derivation of the relation between $\sigma_\ft(\boldsymbol{J}';p)$ and $\sigma_\ft(\boldsymbol{J}';p-1)$ (which give at most a sign). In fact one gets
\begin{equation}
\sigma_\ft(\boldsymbol{J}';p) = \eta\sigma_\ft(\boldsymbol{J}';p-1)
\end{equation}
where
\begin{equation}
\label{eq:eta}
\eta = -\sgn(\boldsymbol{J}'_3)\,\sgn(N+2\boldsymbol{J}'_3-2\boldsymbol{J}'_2+2) \,,
\end{equation}
($\boldsymbol{J}'_i$ denotes the $i$th component of the quintuplet $\boldsymbol{J}'$).
Hence, one obtains
\begin{align}
\label{eq:overlap_t2}
 [T_{\fe,\ft}(\boldsymbol{J}')]_{np}^{mq}  &\sim \delta_{pq}\, (-1)^{m}\,\eta^p \, \cP_n(m;-2\boldsymbol{J}'_2-1,-2\boldsymbol{J}'_1-1,p-N-1,N-p-2\boldsymbol{J}'_2+ 2\boldsymbol{J}'_3+1)\,.
\end{align}
For example, one gets $\eta=+1$ for $\boldsymbol{J}'=\ft\boldsymbol{J}$ and $\fr^2\boldsymbol{J}$, $\eta=-1$ for $\boldsymbol{J}'=\boldsymbol{J}$, $\ft^2\boldsymbol{J}$, $\fr\boldsymbol{J}$, $\fr^3\boldsymbol{J}$ and $\fr^4\boldsymbol{J}$.
\end{rema}

\subsection{Transition matrix between $\pi_{\fe}^{\boldsymbol{J}}$ and $\pi_{\fs}^{\boldsymbol{J}}$}

We now turn to the calculation of the transition matrix $T_{\fe,\fs}= T_{\fe,\fs}(\bJ)$, where the dependence w.r.t. $\boldsymbol{J}$ is dropped in this subsection.
Taking equation \eqref{eq4} with $\mathfrak{h}=\mathfrak{e}$, $\mathfrak{g}=\mathfrak{s}$ and  $X=C_{12}, C_{123}$, we get two equations
\begin{equation}
\begin{aligned}
&\Big(\mu_q^{(12)}\big(\boldsymbol{J}\big) - \mu_n^{(12)}\big(\boldsymbol{J}\big) \Big)
\,[T_{\fe,\fs}]_{np}^{mq}=0\,,\\
&\Big(\mu_m^{(123)}\big(\boldsymbol{J}\big) - \mu_p^{(123)}\big(\boldsymbol{J}\big) \Big)
\,[T_{\fe,\fs}]_{np}^{mq}=0\,.
\end{aligned}
\label{eq:cons1}
\end{equation}
We have used the definitions
\begin{equation}
[\pi^{\bJ}_\fe(C_{12})]\!\ind{n}{m}{p}{q} = \delta_{nm}\,\delta_{pq}\,\mu_m^{(12)}\big(\boldsymbol{J}\big)\quad \text{and}\qquad
[\pi^{\bJ}_\fe(C_{123})]\!\ind{n}{m}{p}{q}\!\!(\boldsymbol{J}) = \delta_{nm}\,\delta_{pq}\,\mu_p^{(123)}\big(\boldsymbol{J}\big)\,,
\end{equation}
and the properties
\begin{equation}
\mu_m^{(12)}\big(\mathfrak{s}\boldsymbol{J}\big)=\mu_m^{(123)}\big(\boldsymbol{J}\big)\quad \text{and}\qquad
\mu_m^{(123)}\big(\mathfrak{s}\boldsymbol{J}\big)=\mu_m^{(12)}\big(\boldsymbol{J}\big).
\end{equation}
Constraints \eqref{eq:cons1} lead to
\begin{equation}
[T_{\fe,\fs}]_{np}^{mq}=\delta_{p,m}\,\delta_{q,n}\,\sigma_\mathfrak{s}(\boldsymbol{J};n,p)\,,
\end{equation}
where $\sigma_\mathfrak{s}(\boldsymbol{J};n,p)$ is a coefficient to be determined.

Performing the same type of calculation with $X=C_{23}$ and $X=C_{34}$, we obtain
\begin{equation}
\sigma_\mathfrak{s}(\boldsymbol{J};n,p)=\sigma_\mathfrak{s}(\boldsymbol{J};q,p) \quad\text{for}\ |n-q|\leq 1 \quad\text{and}\quad
\sigma_\mathfrak{s}(\boldsymbol{J};n,p)=\sigma_\mathfrak{s}(\boldsymbol{J};n,m) \quad\text{for}\ |p-m|\leq 1\,.
\end{equation}
This proves that the coefficient $\sigma_\mathfrak{s}(\boldsymbol{J};n,p)=\sigma_\mathfrak{s}(\boldsymbol{J})$ does not depend on $n$ and $p$.
Furthermore, demanding that the orthogonality relation \eqref{eq:ortho} is satisfied imposes that $\sigma_\fs(\boldsymbol{J})^2=1$.
Then, we checked that the other constraints \eqref{eq4} are automatically satisfied for any $X\in sR(4)$.
This shows that the transition matrix is the permutation operator up to a sign: 
\begin{equation}\label{eq:Tes}
 T_{\fe,\fs}\sim\,P
\end{equation}
with $[P]_{np}^{mq} =\delta_{p,m}\,\delta_{q,n}$.

\subsection{Transition matrix between $\pi_{\fe}^{\boldsymbol{J}}$ and $\pi_{\fii}^{\boldsymbol{J}}$}

The calculation of the transition matrix $T_{\fe,\fii}= T_{\fe,\fii}(\bJ)$ follows the same lines. 
Taking equation \eqref{eq4} with $\mathfrak{h}=\mathfrak{e}$, $\mathfrak{g}=\mathfrak{i}$ and  $X=C_{12}, C_{123}$ implies 
that $T_{\fe,\fii}(\bJ)$ is diagonal:
\begin{equation}
[T_{\fe,\fii}]_{np}^{mq}=[T_{\fe,\fii}]_{np}^{np}\,\delta_{n,m}\,\delta_{p,q} \,.
\end{equation}
The equations for $X=C_{23}, C_{13}$ are then trivially satisfied. Eq. \eqref{eq4} leads for the remaining generators (those containing an index 4) to the following relations
\begin{equation}
[T_{\fe,\fii}]_{n+1,p}^{n+1,p}=[T_{\fe,\fii}]_{np}^{np} \qquad \text{and} \qquad [T_{\fe,\fii}]_{n,p+1}^{n,p+1}=-[T_{\fe,\fii}]_{np}^{np} \,.
\end{equation}
Therefore, one gets
\begin{equation}
[T_{\fe,\fii}]_{np}^{mq} =\sigma_{\fii}(\boldsymbol{J}) \, (-1)^p \, \delta_{n,m}\,\delta_{p,q} \,.
\end{equation}
Finally, demanding that $T_{\fe,\fii}$ is an orthogonal matrix, we get
\begin{equation}
[T_{\fe,\fii}]_{np}^{mq} \sim (-1)^p \, \delta_{n,m}\,\delta_{p,q}\,.
\end{equation}

\subsection{Equivalence of the representations and transition matrix \label{subsect:4.4}}

In this section, we shall show that the previous computations are enough to determine the transition matrix $T_{\fg,\fh}$ for any $\fh,\fg\in \fS_5$.
If the transition matrices  $T_{\fg,\fh}$ exist, they satisfy the properties given in the following proposition.
\begin{prop}\label{prop:compo-T} For any $\fg,\fh,\fu\in \fS_5$, one gets
\begin{equation} \label{eq:pru1}
T_{\fg,\fh}(\bJ)\sim T_{\fg,\fu}(\bJ)T_{\fu,\fh}(\bJ)\,,\qquad   T_{\fg,\fh}(\bJ)\sim T_{\fh,\fg}(\bJ)^t \qquad\text{and}\qquad  T_{\fu\fh,\fg\fh}(\bJ) \sim T_{\fu,\fg}(\fh\bJ)\,.
\end{equation} 
We recall that $\sim$ stands for \textit{``equal up to a global sign''}.
\end{prop}
\proof From the definition of the transition matrices, one gets for any $X\in sR(4)$, 
\begin{equation}
T_{\fg,\fh}(\bJ)^{-1}  T_{\fg,\fu}(\bJ)T_{\fu,\fh}(\bJ) \pi_\fe^{\fh\bJ}(\fh X) =  \pi_\fe^{\fh \bJ}(\fh X)T_{\fg,\fu}(\bJ)T_{\fu,\fh}(\bJ)T_{\fg,\fh}(\bJ)^{-1} \,.
\end{equation}
Since the representation $\pi_\fe^{\bJ}$ is irreducible, $T_{\fg,\fh}(\bJ)^{-1}  T_{\fg,\fu}(\bJ)T_{\fu,\fh}(\bJ)$ is proportional to the identity by the Schur lemma. 
The orthogonality relation of the transition matrices implies that the factor of proportionality is a sign, which proves the first relation of the proposition.
Both other relations are proven by similar arguments.
\endproof

Up to now, we have proved the existence  of the transition matrices $T_{\fg,\fh}(\bJ)$, with $\fg=\fe$ and $\fh = \ft\,,\,\fs\,,\,\fii$ solely. 
Then, at this point,
the above proposition can be applied only for these matrices. 
The following proposition ensures the existence (and thus the validity of proposition 
\ref{prop:compo-T}) for all transition matrices $T_{\fg,\fh}(\bJ)$, with $\fg\,,\,\fh\in\fS_5$.
\begin{prop}
 All the representations $\pi_{\fg}^{\bJ}$, $g\in\fS_5$, defined by \eqref{eq:defpiG} are equivalent. 
\end{prop}
\proof The equivalence between $\pi_{\fg}^{\bJ}$ and $\pi_{\fh}^{\bJ}$ is proven by showing that there exists a transition matrix $T_{\fg,\fh}$.
Relations \eqref{eq4} for $T_{\fh,\ft\fh}(\bJ)$ (for $\mathfrak{h}\in \fS_5$ and $X\in sR(4)$) reads
\begin{equation}
T_{\fh,\ft\fh}(\bJ)\ \pi_\fe^{\ft\fh\bJ}(\ft\fh(X))=\pi_\fe^{\fh\bJ}(\fh(X))\ T_{\fh,\ft\fh}(\bJ)\,,
\end{equation}
whereas the one for $T_{\fe,\ft}(\fh\bJ)$ (with $X$ replaced by $\fh(X)$) is recast as  
\begin{equation}
T_{\fe,\ft}(\fh\bJ)\ \pi_\fe^{\ft\fh\bJ}(\ft\fh(X))=\pi_\fe^{\fh\bJ}(\fh(X))\ T_{\fe,\ft}(\fh\bJ)\,.
\end{equation}
We deduce that $T_{\fh,\ft\fh}(\bJ)$ and $T_{\fe,\ft}(\fh\bJ)$ satisfy the same set of relations. In section \ref{sec:ot}, we proved that this set has a unique (up to a global sign) solution given by \eqref{eq:overlap_t}.
Therefore one gets, for any $\fh\in \fS_5$,
\begin{equation}
 T_{\fh,\ft\fh}(\bJ) \sim T_{\fe,\ft}(\fh\bJ)\,.
\end{equation}
 The following relations are proven similarly
\begin{equation}
 T_{\fh,\fs\fh}(\bJ) \sim T_{\fe,\fs}(\fh\bJ) ,\quad  T_{\fh,\fii\fh}(\bJ) \sim T_{\fe,\fii}(\fh\bJ)\,.
\end{equation}
Finally for any $\mathfrak{g} \in \fS_5$, one knows that $\mathfrak{g}$ can be written as a product of $\mathfrak{s}$, $\mathfrak{t}$ and $\fii$. Therefore, 
any $T_{\fg,\fh}$ can be expressed as a product of $T_{\fe,\ft}$, $T_{\fe,\fs}$ and $T_{\fe,\fii}$ by using the first relation in \eqref{eq:pru1}, hence the proposition.\endproof

From the previous proof, one gets a way to compute all the transition matrices in terms of $T_{\fe,\ft}$, $T_{\fe,\fs}$ and $T_{\fe,\fii}$. 
For example, one gets  
\begin{equation}
 T_{\fe,\fs\ft\fs}(\bJ)\sim T_{\fe,\ft\fs}(\bJ)T_{\ft\fs,\fs\ft\fs}(\bJ)\sim T_{\fe,\fs}(\bJ)T_{\fs,\ft\fs}(\bJ)T_{\ft\fs,\fs\ft\fs}(\bJ)
 \sim T_{\fe,\fs}(\bJ)T_{\fe,\ft}(\fs\bJ)T_{\fe,\fs}(\ft\fs\bJ)\;.
\end{equation}
Using the third relation in \eqref{eq:pru1} with $\fu=\fe$ and $\fh=\fg^{-1}$, and the second one, we deduce that 
\begin{equation}
 T_{\fe,\fg^{-1}}(\bJ) \sim T_{\fe,\fg}(\fg^{-1}\bJ) ^t.  \label{eq:invt}
\end{equation}

\subsection{Overlap coefficients and the icosidodecahedron \label{sec:OCI}}

This subsection provides a geometrical description of the transition matrices.

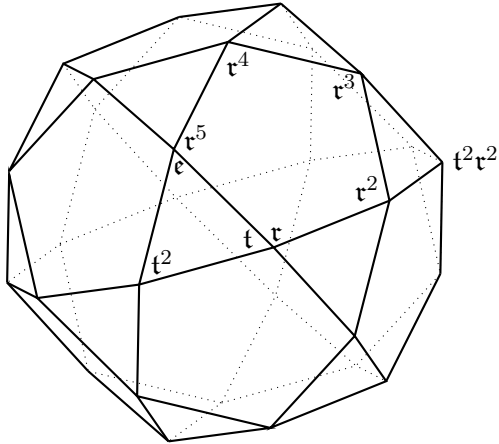
\begin{figure}[htbp]
\begin{center}
\begin{tikzpicture}[scale=1.5]
\node at (-0.4, 0.49) {$\fe$};
\node at (-0.25, 0.75) {$\fr^5$};

\node at (-0.55, -0.35) {$\ft^2$};

\node at (0.2, -0.15) {$\ft$};
\node at (0.45, -0.1) {$\fr$};

\node at (0.12, 1.4) {$\fr^4$};

\node at (1.25, 0.3)  {$\fr^2$}; 

\node at (1.05, 1.2)   {$\fr^3$};

\node at (2.2, 0.57)  {$\ft^2\fr^2$};

\draw [thick] (-1.908111787, -.5319769439)--(-1.890410045, .4242268192)-- (-1.643196449, -.6682790668) -- (-1.908111787, -.5319769439);

\draw [dotted]  (-1.908111787, -.5319769439)--(-1.196979681, -1.314983596)--(-1.444193277, -.1924777098) -- (-1.908111787, -.5319769439);

\draw [thick] (-1.908111787, -.5319769439)--(-1.196979681, -1.314983596);

\draw [dotted] (-1.196979681, -1.314983596)--(-.4925605311, -1.935210443)--(-0.0286420207, -1.595711209) --(-1.196979681, -1.314983596);

\draw [thick] (-1.196979681, -1.314983596)--(-.4925605311, -1.935210443);

\draw  [thick] (-.4925605311, -1.935210443)--(-.7683376602, -1.535525063)--(.3999999999, -1.816252676) -- (-.4925605311, -1.935210443);

\draw  [thick] (-.7683376602, -1.535525063)--(-1.643196449, -.6682790668)--(-.7506359181, -.5493212998) --  (-.7683376602, -1.535525063);

\draw  [dotted] (.7683376602, 1.535525063)--(.4925605311, 1.935210443)--(-.3999999999, 1.816252676) --(.7683376602, 1.535525063);

\draw [thick] (.4925605311, 1.935210443)--(-.3999999999, 1.816252676) ;

\draw  [dotted] (1.643196449, .6682790668)--(.7683376602, 1.535525063)--(.7506359181, .5493212998)--(1.643196449, .6682790668);

\draw [dotted] (1.908111787, .5319769439)--(1.643196449, .6682790668)--(1.890410045, -.4542268192)-- (1.908111787, .5319769439);

\draw [thick] (1.890410045, -.4542268192)-- (1.908111787, .5319769439);

\draw [ thick]  (1.908111787, .5319769439)--(1.444193277, .1924777098)--(1.196979681, 1.314983596) -- (1.908111787, .5319769439);

\draw [thick]  (.4925605311, 1.935210443)--(1.196979681, 1.314983596)--(0.0286420207, 1.595711209) -- (.4925605311, 1.935210443);

\draw [dotted]  (-1.139774127, 1.003548119)--(-.3999999999, 1.816252676)--(-1.415551256, 1.403233499) -- (-1.139774127, 1.003548119);

\draw [thick] (-.3999999999, 1.816252676)--(-1.415551256, 1.403233499);

\draw [dotted] (-1.444193277, -.1924777098)--(-.4286420207, .2205414672)--(-1.139774127, 1.003548119)-- (-1.444193277, -.1924777098);

\draw [dotted]  (.7506359181, .5493212998)--(-.4286420207, .2205414672)--(.4462167681, -.6467045291) --(.7506359181, .5493212998);

\draw  [dotted] (-0.286420207e-1, -1.595711209)--(1.150635918, -1.266931376)--(.4462167681, -.6467045291) --(-0.286420207e-1, -1.595711209);

\draw  [dotted] (1.415551256, -1.403233499)--(1.890410045, -.4542268192)--(1.150635918, -1.266931376) -- (1.415551256, -1.403233499);

\draw [thick] (1.415551256, -1.403233499)--(1.890410045, -.4542268192);

\draw  [thick] (1.139774127, -1.003548119)--(1.415551256, -1.403233499)--(.3999999999, -1.816252676) -- (1.139774127, -1.003548119);

\draw  [thick] (1.444193277, .1924777098)--(1.139774127, -1.003548119)--(.4286420207, -.2205414672) --  (1.444193277, .1924777098);

\draw [thick]  (.4286420207, -.2205414672)--(-.7506359181, -.5493212998)--(-.4462167681, .6467045291) -- (.4286420207, -.2205414672);

\draw [ thick]  (0.286420207e-1, 1.595711209)--(-.4462167681, .6467045291)--(-1.150635918, 1.266931376) -- (0.286420207e-1, 1.595711209);

\draw [ thick] (-1.890410045, .4542268192) -- (-1.415551256, 1.403233499)--(-1.150635918, 1.266931376) --(-1.890410045, .4542268192);

\end{tikzpicture}

\caption{The labels of the vertices of the icosidodecahedron in terms of elements of $\fS_5$. We recall that $\ft$ and $\fr$ are defined by \eqref{eq:genert} and \eqref{eq:gener}.
  \label{fig2}}
\end{center}
\end{figure}

For a given isomorphism $\mathfrak{g} \in \fS_5$ of $sR(4)$ defined in proposition \ref{eq:propiso}, the matrices  $\pi^{\bJ}_\fg(\mathfrak{g}^{-1}(C_{12}))$ and $\pi^{\bJ}_\fg(\mathfrak{g}^{-1}(C_{123}))$ are diagonal.
Indeed one gets 
\begin{equation}\label{eq:CC}
\begin{aligned}
\pi^{\bJ}_\fg(\mathfrak{g}^{-1}(C_{12}))=\pi^{\fg\bJ}_\fe(C_{12}),\qquad \pi^{\bJ}_\fg(\mathfrak{g}^{-1}(C_{123}))=\pi^{\fg\bJ}_\fe(C_{123}),
\end{aligned}
\end{equation}
which are diagonal.
We associate the map $\mathfrak{g} \in \fS_5$ with the vertex $(\mathfrak{g}^{-1}(C_{12}),\mathfrak{g}^{-1}(C_{123}))$ of the icosidodecahedron, as illustrated in figure \ref{fig2}.
Let us notice that different elements of $\fS_5$ can be associated to the same vertex.
 
The transition matrix $T_{\fh,\fg}$ corresponds to the change of basis from the one where 
$\mathfrak{h}^{-1}(C_{12})$ and $\mathfrak{h}^{-1}(C_{123})$ are diagonal 
to the one where $\mathfrak{g}^{-1}(C_{12})$ and $\mathfrak{g}^{-1}(C_{123})$ are diagonal. 
A path on the icosidodecahedron between the vertices labeled by $\fh$ and $\fg$ can then be associated to this transition matrix  $T_{\fh\fg}$.
Indeed, let us choose a path between the vertices $\fh=\fg_0$ to $\fg=\fg_\ell$, given by the sequence of vertices labeled $\fg_0$, $\fg_1,\dots,\fg_{\ell-1}$ and $\fg_\ell$.
Finally, using \eqref{eq:pru1}, the transition matrix $T_{\fg_0,\fg_\ell}$ can be written as follows 
\begin{equation}\label{eq:Tpath}
 T_{\fg_0,\fg_\ell}\sim T_{\fg_0,\fg_1}T_{\fg_1,\fg_2} \dots T_{\fg_{\ell-1},\fg_\ell}\ .
\end{equation}
Let us emphasize that the value \eqref{eq:Tpath} of $T_{\fg_0,\fg_\ell}$ does not depend on the chosen path.

Each transition matrix $T_{\fg_i,\fg_{i+1}}$  in the r.h.s. of \eqref{eq:Tpath} corresponds just to one edge of the icosidodecahedron and can be rewritten as, using \eqref{eq:pru1},
\begin{equation}
 T_{\fg_i,\fg_{i+1}}(\bJ) \sim T_{\fe,\fg_{i+1}\fg_i^{-1}}(\fg_i\bJ)\ .
\end{equation}
The transition matrix on the r.h.s. of the previous relation corresponds to one of the four edges connected to the vertex $\fe$ and can be expressed in terms of monovariate Racah polynomials.
Therefore expression \eqref{eq:Tpath} allows one to give an expression of any transition matrix in terms of monovariate Racah polynomials.

The different expressions of $T_{\fh,\fg}$ depending on the choice of the path provide relation between the monovariate Racah polynomial.
We exploit this result in the next section to give new proof of known relations.

\section{Properties of the monovariate Racah polynomials \label{sec:pmr}}

In this section, we study the relations satisfied by the monovariate Racah polynomials, as explained in Section \ref{sec:OCI}, 
for cycles on the icosidodecahedron. The basic examples are given by the triangle and the pentagon faces of the icosidodecahedron, 
which correspond to the identities $\mathfrak{t}^3=\fe$ and $\fr^5=\fe$ (see relation \eqref{eq:gener}), respectively. 
These two cycles can be illustrated in red on the following figures:
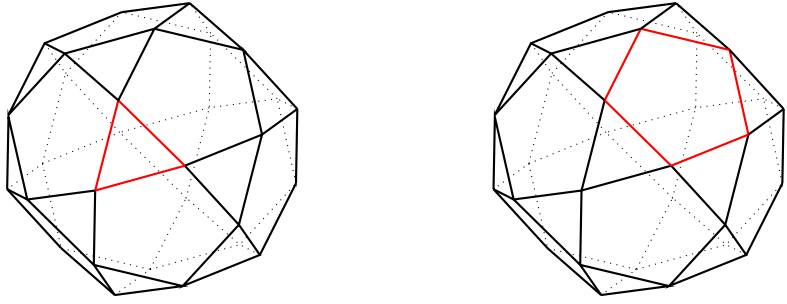
\begin{center}
\begin{tabular}{cc}
\begin{tikzpicture}[scale=1]
\draw [thick] (-1.908111787, -.5319769439)--(-1.890410045, .4242268192)-- (-1.643196449, -.6682790668) -- (-1.908111787, -.5319769439);

\draw [dotted]  (-1.908111787, -.5319769439)--(-1.196979681, -1.314983596)--(-1.444193277, -.1924777098) -- (-1.908111787, -.5319769439);

\draw [thick] (-1.908111787, -.5319769439)--(-1.196979681, -1.314983596);

\draw [dotted] (-1.196979681, -1.314983596)--(-.4925605311, -1.935210443)--(-0.0286420207, -1.595711209) --(-1.196979681, -1.314983596);

\draw [thick] (-1.196979681, -1.314983596)--(-.4925605311, -1.935210443);

\draw  [thick] (-.4925605311, -1.935210443)--(-.7683376602, -1.535525063)--(.3999999999, -1.816252676) -- (-.4925605311, -1.935210443);

\draw  [thick] (-.7683376602, -1.535525063)--(-1.643196449, -.6682790668)--(-.7506359181, -.5493212998) --  (-.7683376602, -1.535525063);

\draw  [dotted] (.7683376602, 1.535525063)--(.4925605311, 1.935210443)--(-.3999999999, 1.816252676) --(.7683376602, 1.535525063);

\draw [thick] (.4925605311, 1.935210443)--(-.3999999999, 1.816252676) ;

\draw  [dotted] (1.643196449, .6682790668)--(.7683376602, 1.535525063)--(.7506359181, .5493212998)--(1.643196449, .6682790668);

\draw [dotted] (1.908111787, .5319769439)--(1.643196449, .6682790668)--(1.890410045, -.4542268192)-- (1.908111787, .5319769439);

\draw [thick] (1.890410045, -.4542268192)-- (1.908111787, .5319769439);

\draw [ thick]  (1.908111787, .5319769439)--(1.444193277, .1924777098)--(1.196979681, 1.314983596) -- (1.908111787, .5319769439);

\draw [thick]  (.4925605311, 1.935210443)--(1.196979681, 1.314983596)--(0.0286420207, 1.595711209) -- (.4925605311, 1.935210443);

\draw [dotted]  (-1.139774127, 1.003548119)--(-.3999999999, 1.816252676)--(-1.415551256, 1.403233499) -- (-1.139774127, 1.003548119);

\draw [thick] (-.3999999999, 1.816252676)--(-1.415551256, 1.403233499);

\draw [dotted] (-1.444193277, -.1924777098)--(-.4286420207, .2205414672)--(-1.139774127, 1.003548119)-- (-1.444193277, -.1924777098);

\draw [dotted]  (.7506359181, .5493212998)--(-.4286420207, .2205414672)--(.4462167681, -.6467045291) --(.7506359181, .5493212998);

\draw  [dotted] (-0.286420207e-1, -1.595711209)--(1.150635918, -1.266931376)--(.4462167681, -.6467045291) --(-0.286420207e-1, -1.595711209);

\draw  [dotted] (1.415551256, -1.403233499)--(1.890410045, -.4542268192)--(1.150635918, -1.266931376) -- (1.415551256, -1.403233499);

\draw [thick] (1.415551256, -1.403233499)--(1.890410045, -.4542268192);

\draw  [thick] (1.139774127, -1.003548119)--(1.415551256, -1.403233499)--(.3999999999, -1.816252676) -- (1.139774127, -1.003548119);

\draw  [thick] (1.444193277, .1924777098)--(1.139774127, -1.003548119)--(.4286420207, -.2205414672) --  (1.444193277, .1924777098);

\draw [thick, color=red]  (.4286420207, -.2205414672)--(-.7506359181, -.5493212998);
\draw [thick, color=red] (-.7506359181, -.5493212998)--(-.4462167681, .6467045291) ;
\draw [thick, color=red]  (-.4462167681, .6467045291) -- (.4286420207, -.2205414672);

\draw [ thick]  (0.286420207e-1, 1.595711209)--(-.4462167681, .6467045291)--(-1.150635918, 1.266931376) -- (0.286420207e-1, 1.595711209);

\draw [ thick] (-1.890410045, .4542268192) -- (-1.415551256, 1.403233499)--(-1.150635918, 1.266931376) --(-1.890410045, .4542268192);

\end{tikzpicture}
&\hspace{2cm}
\begin{tikzpicture}[scale=1]
\draw [thick] (-1.908111787, -.5319769439)--(-1.890410045, .4242268192)-- (-1.643196449, -.6682790668) -- (-1.908111787, -.5319769439);

\draw [dotted]  (-1.908111787, -.5319769439)--(-1.196979681, -1.314983596)--(-1.444193277, -.1924777098) -- (-1.908111787, -.5319769439);

\draw [thick] (-1.908111787, -.5319769439)--(-1.196979681, -1.314983596);

\draw [dotted] (-1.196979681, -1.314983596)--(-.4925605311, -1.935210443)--(-0.0286420207, -1.595711209) --(-1.196979681, -1.314983596);

\draw [thick] (-1.196979681, -1.314983596)--(-.4925605311, -1.935210443);

\draw  [thick] (-.4925605311, -1.935210443)--(-.7683376602, -1.535525063)--(.3999999999, -1.816252676) -- (-.4925605311, -1.935210443);

\draw  [thick] (-.7683376602, -1.535525063)--(-1.643196449, -.6682790668)--(-.7506359181, -.5493212998) --  (-.7683376602, -1.535525063);

\draw  [dotted] (.7683376602, 1.535525063)--(.4925605311, 1.935210443)--(-.3999999999, 1.816252676) --(.7683376602, 1.535525063);

\draw [thick] (.4925605311, 1.935210443)--(-.3999999999, 1.816252676) ;

\draw  [dotted] (1.643196449, .6682790668)--(.7683376602, 1.535525063)--(.7506359181, .5493212998)--(1.643196449, .6682790668);

\draw [dotted] (1.908111787, .5319769439)--(1.643196449, .6682790668)--(1.890410045, -.4542268192)-- (1.908111787, .5319769439);

\draw [thick] (1.890410045, -.4542268192)-- (1.908111787, .5319769439);

\draw [ thick]  (1.908111787, .5319769439)--(1.444193277, .1924777098);
\draw [ thick, color=red]  (1.444193277, .1924777098)--(1.196979681, 1.314983596);
\draw [ thick]  (1.196979681, 1.314983596) -- (1.908111787, .5319769439);

\draw [thick]  (.4925605311, 1.935210443)--(1.196979681, 1.314983596);
\draw [thick, color=red]  (1.196979681, 1.314983596)--(0.0286420207, 1.595711209);
\draw [thick]  (0.0286420207, 1.595711209) -- (.4925605311, 1.935210443);

\draw [dotted]  (-1.139774127, 1.003548119)--(-.3999999999, 1.816252676)--(-1.415551256, 1.403233499) -- (-1.139774127, 1.003548119);

\draw [thick] (-.3999999999, 1.816252676)--(-1.415551256, 1.403233499);

\draw [dotted] (-1.444193277, -.1924777098)--(-.4286420207, .2205414672)--(-1.139774127, 1.003548119)-- (-1.444193277, -.1924777098);

\draw [dotted]  (.7506359181, .5493212998)--(-.4286420207, .2205414672)--(.4462167681, -.6467045291) --(.7506359181, .5493212998);

\draw  [dotted] (-0.286420207e-1, -1.595711209)--(1.150635918, -1.266931376)--(.4462167681, -.6467045291) --(-0.286420207e-1, -1.595711209);

\draw  [dotted] (1.415551256, -1.403233499)--(1.890410045, -.4542268192)--(1.150635918, -1.266931376) -- (1.415551256, -1.403233499);

\draw [thick] (1.415551256, -1.403233499)--(1.890410045, -.4542268192);

\draw  [thick] (1.139774127, -1.003548119)--(1.415551256, -1.403233499)--(.3999999999, -1.816252676) -- (1.139774127, -1.003548119);

\draw  [thick] (1.444193277, .1924777098)--(1.139774127, -1.003548119);
\draw  [thick] (1.139774127, -1.003548119)--(.4286420207, -.2205414672) ;
\draw  [thick, color=red] (.4286420207, -.2205414672) --  (1.444193277, .1924777098);

\draw [thick]  (.4286420207, -.2205414672)--(-.7506359181, -.5493212998);
\draw [thick] (-.7506359181, -.5493212998)--(-.4462167681, .6467045291) ;
\draw [thick, color=red]  (-.4462167681, .6467045291) -- (.4286420207, -.2205414672);

\draw [ thick, color=red]  (0.286420207e-1, 1.595711209)--(-.4462167681, .6467045291);
\draw [ thick]  (-.4462167681, .6467045291)--(-1.150635918, 1.266931376) ;
\draw [ thick]  (-1.150635918, 1.266931376) -- (0.286420207e-1, 1.595711209);

\draw [ thick] (-1.890410045, .4542268192) -- (-1.415551256, 1.403233499)--(-1.150635918, 1.266931376) --(-1.890410045, .4542268192);

\end{tikzpicture}
\end{tabular}
\end{center}
As we will see, they are intimately related to the Racah relation and the Biedenharn--Eliott relation.
Other cycles can be considered, such as a cycle involving a triangle and a pentagon together. However, it can be deduced by gluing together the basic ones and 
using unitarity: we illustrate this on a relation involving six Racah polynomials, see section \ref{sec:otre}.

\subsection{Racah relation}

We look at the interpretation of the relation $\mathfrak t^3=\fe$.
Using the results of section \ref{subsect:4.4}, we deduce that 
\begin{equation}
T_{\fe,\ft^3}(\bJ) \sim T_{\fe,\ft}(\bJ) T_{\fe,\ft}(\ft\bJ) T_{\fe,\ft}(\ft^2\bJ) \,,
\end{equation}
that is
\begin{equation}\label{eq:t3}
\sum_{\atopn{k_1,\ell_1\geq0}{k_1+\ell_1\leq N}}\sum_{\atopn{k_2,\ell_2\geq0}{k_2+\ell_2\leq N}}  [T_{\fe,\ft}(\bJ)]_{np}^{k_1\ell_1} \, [T_{\fe,\ft}(\ft\bJ)]_{k_1\ell_1}^{k_2\ell_2} \, [T_{\fe,\ft}(\ft^2\bJ)]_{k_2\ell_2}^{mq}
\sim \delta_{n,m}\,\delta_{p,q} \,.
\end{equation}
It corresponds to the cycle around the triangle $(C_{12},C_{123})$, $(C_{13},C_{123})$ and $(C_{23},C_{123})$.
Then using expression \eqref{eq:overlap_t} for $T_{\fe,\ft}(\bJ)$ and the action \eqref{eq:A5-quintu} of $\ft$ on $\bJ$, 
for parameters $m,n,p\geq0$ such that $m+p\leq N$ and $n+p\leq N$, this expression becomes
\begin{equation}\begin{aligned}
\sum_{k_1,k_2=0}^{N-p} &(-1)^{k_1+k_2+m} \, \cP_n(k_1;-2j^{(2)}-1,-2j^{(1)}-1,p-N-1,N-p-2j^{(2)}+ 2j^{(3)}+1) \, \\
& \times\ \cP_{k_1}(k_2;2j^{(3)}+1,-2j^{(2)}-1,p-N-1,N-p-2j^{(1)}+ 2j^{(3)}+1)   \\
& \times\  \cP_{k_2}(m;-2j^{(1)}-1,2j^{(3)}+1,p-N-1,N-p-2j^{(1)}- 2j^{(2)}-1)\\
\sim \delta_{n,m}\,.
\end{aligned}
\end{equation}
When the parameters $j^{(k)}$ are nonnegative half-integers, we showed that this last equation is a consequence of the Racah and orthogonality  relations
for the $6j$-symbols \cite{Messiah} .

\subsection{Biedenharn--Elliot relation}

We recall that the element $\mathfrak{r}$ of $\fS_{5}$ corresponds to the rotation of the icosidodecahedron with an angle $2\pi/5$ around the axis passing through 
its center and the center of the pentagon $(C_{12},C_{123})$, $(C_{23},C_{123})$, $(C_{234},C_{23})$, $(C_{34},C_{234})$, $(C_{12},C_{34})$. 
It satisfies $\mathfrak{r}^5=\mathfrak{e}$ and one defines a function $\fr$ on $\bJ$ as
\begin{equation}
\label{eq:rJ}
\fr\bJ=(j^{(4)}, j^{(0)}, -j^{(1)}-1, j^{(2)}, -j^{(3)}-1). 
\end{equation}
We now turn to the calculation of the transition matrix $T_{\fe,\fr}(\bJ)$, before investigating the Biedenharn--Elliot relation. 

\paragraph{Transition matrix $T_{\fe,\fr}(\bJ)$.} The computation of $T_{\fe,\fr}(\bJ)$ follows the same lines as the computation of $T_{\fe,\ft}(\bJ)$. 
The equation \eqref{eq4} for $X=C_{123}$ leads to
\begin{equation}
\big( \mu_m^{(123)}(\boldsymbol{J})-\mu_p^{(123)}(\boldsymbol{J}) \big) \, [T_{\fe,\fr}(\bJ)]_{np}^{mq} = 0 \,,
\end{equation}
hence $[T_{\fe,\fr}(\bJ)]_{np}^{mq}=\delta_{mp}\ [T_{\fe,\fr}(\bJ)]_{np}^{pq}$. 

Considering $X=C_{23}$ and $X=C_{12}$ in \eqref{eq4}, and using the relations $\rho_{p,q}(\mathfrak{r}\boldsymbol{J}) = \psi_{q,p}(\mathfrak{t}\boldsymbol{J})$ and $\rho^{00}_{p,q}(\mathfrak{r}\boldsymbol{J})=\big( \mu_p^{(123)}+\mu^{(1)}+\mu^{(2)}+\mu^{(3)}-\mu_q^{(12)}-\psi_{q,p}^{00}  \big)(\mathfrak{t}\boldsymbol{J})$, one finds that $(-1)^{q}\ [T_{\fe,\fr}]_{np}^{pq}$ satisfies the same recurrence and difference equations than $[T_{\fe,\ft}]_{np}^{qp}$, see \eqref{eq:reccurracah} and \eqref{eq:diff}. 

Equation \eqref{eq4} for $X=C_{234}$ leads to, for the component $\bra{n,p}\eqref{eq4}\ket{m,q}$ with $m=p-1$:
\begin{equation}
\label{eqC234r1}
\psi_{p,q}(\fr\boldsymbol{J}) \,[T_{\fe,\fr}]_{np}^{pq} = \varphi_{n+1,p}^{A}(\boldsymbol{J}) \ [T_{\fe,\fr}]_{n+1,p-1}^{p-1,q} 
+ \varphi_{n,p}^{V}(\boldsymbol{J}) \,[T_{\fe,\fr}]_{n,p-1}^{p-1,q} + \varphi_{n,p}^{D}(\boldsymbol{J})\,[T_{\fe,\fr}]_{n-1,p-1}^{p-1,q}  \,.
\end{equation}
A similar relation is obtained for $m=p+1$. 
Using the relation $\psi_{p,q}(\mathfrak{r}\boldsymbol{J}) = \rho_{q,p}(\mathfrak{t}\boldsymbol{J})$ shows that $(-1)^{p+q} [T_{\fe,\fr}]_{np}^{pq}$ also satisfies the equation \eqref{eq4.16A} obeyed by $[T_{\fe,\ft}]_{np}^{qp}$ (together with a similar one associated to $m=p+1$). The case $m=p$ is treated along the same lines, using now $\psi^{00}_{p,q}(\mathfrak{r}\boldsymbol{J})=\big( \mu_q^{(12)}+\mu^{(3)}+\mu^{(4)}+\mu^{(0)}-\mu_p^{(123)}-\rho_{q,p}^{00}  \big)(\mathfrak{t}\boldsymbol{J})$.

In the same way, equation \eqref{eq4} for $X=C_{34}$ shows that $(-1)^{p+q} [T_{\fe,\fr}]_{np}^{pq}$ satisfies the  equation \eqref{eq4.16B} obeyed by $[T_{\fe,\ft}]_{np}^{qp}$, using now the relations $\varphi_{p,q}^{D}(\mathfrak{r}\boldsymbol{J}) = \varphi_{q,p}^{D}(\mathfrak{t}\boldsymbol{J})$, 
$\varphi_{p,q}^{A}(\mathfrak{r}\boldsymbol{J}) = \varphi_{q,p}^{A}(\mathfrak{t}\boldsymbol{J})$ and 
$\varphi_{p,q}^{H}(\mathfrak{r}\boldsymbol{J}) = (\rho_{q,p}-\varphi_{q,p}^{V})(\mathfrak{t}\boldsymbol{J})$. 

The other generators do not lead to new conditions. It follows that $(-1)^{p+q} [T_{\fe,\fr}]_{np}^{pq}$ can be identified with $[T_{\fe,\ft}]_{np}^{qp}$, that is the transition matrix $T_{\fe,\fr}(\bJ)$ takes the following form:
\begin{align}
\label{eq:overlap_r}
[T_{\fe,\fr}(\bJ)]_{np}^{mq} &\sim \delta_{mp} \, \cP_n(q;-2j^{(2)}-1,-2j^{(1)}-1,p-N-1,N-p-2j^{(2)}+ 2j^{(3)}+1) \,.
\end{align}
\begin{rema}
In the same way, formula \eqref{eq:overlap_r} can be generalized for a generic quintuplet $\bJ'=\fg\bJ$ and one gets 
\begin{align}
\label{eq:overlap_r2}
[T_{\fe,\fr}(\bJ')]_{np}^{mq} &\sim \delta_{mp} (-\eta)^m\, \cP_n(q;-2\boldsymbol{J}'_2-1,-2\boldsymbol{J}'_1-1,p-N-1,N-p-2\boldsymbol{J}'_2+ 2\boldsymbol{J}'_3+1) \,,
\end{align}
where $\eta$ is given by \eqref{eq:eta}.
\end{rema}

\paragraph{Biedenharn--Elliot relation.} 
Now the relation $\mathfrak{r}^5=\fe$ leads to the following equation:
\begin{align}
\delta_{n,m}\delta_{p,q} &= [T_{\fe,\fr^5}(\mathfrak{r}\bJ)]_{np}^{mq} \nonumber \\
&\sim \sum_{\atopn{k_1,k_2,k_3,k_4}{\ell_1,\ell_2,\ell_3,\ell_4}}
[T_{\fe,\fr}(\bJ)]_{np}^{k_1\ell_1} \, [T_{\fe,\fr}(\fr\bJ)]_{k_1\ell_1}^{k_2\ell_2} \, [T_{\fe,\fr}(\fr^2\bJ)]_{k_2\ell_2}^{k_3\ell_3} \, [T_{\fe,\fr}(\fr^3\bJ)]_{k_3\ell_3}^{k_4\ell_4} \, [T_{\fe,\fr}(\fr^4\bJ)]_{k_4\ell_4}^{mq} \,,
\end{align}
where the sum is performed on integers $k_j,\ell_j\geq0$ such that $k_j+\ell_j\leq N$. 
Using \eqref{eq:overlap_r} and the action of $\fr$ on $\bJ$, see \eqref{eq:rJ}, one finally obtains
\begin{equation}
\begin{aligned}
& \delta_{n,m}\delta_{p,q} \sim \sum_{a,b,c} (-1)^b\,\cP_n(a;-2j^{(2)}-1,-2j^{(1)}-1,p-N-1,N-p-2j^{(2)}+ 2j^{(3)}+1) \\
&\qquad\times\cP_p(b;-2j^{(0)}-1,-2j^{(4)}-1,a-N-1,N-a-2j^{(0)}- 2j^{(1)}-1)  \\ 
& \qquad\times \cP_a(c;2j^{(3)}+1,-2j^{(2)}-1,b-N-1,N-b-2j^{(4)}+ 2j^{(3)}+1) \\
&\qquad\times \cP_b(m;-2j^{(1)}-1,-2j^{(0)}-1,c-N-1,N-c-2j^{(2)}- 2j^{(1)}-1) \\
&\qquad\times \cP_c(q;-2j^{(4)}-1,2j^{(3)}+1,m-N-1,N-m-2j^{(0)}- 2j^{(4)}-1) \, ,
\end{aligned}
\end{equation}
where now $a,b,c\geq0$ satisfy $a\leq N-p$, $a+b\leq N$, $c+b\leq N$, $c\leq N-m$.
This last equation is similar to the Biedenharn--Elliott relation for the $6j$-symbols when the parameters $j^{(k)}$ are nonnegative half-integers \cite{Bieden, Elliott}.

\subsection{Additional relations for the Racah polynomials \label{sec:otre}}

Other relations of the same type can be obtained using the algebra $\fS_5$. 
For instance, defining $\fu=\fs\ft\fs\ft^2\fs$, we have again $\fu^2=\fe$. 
This leads to the relation 
\begin{equation}\label{eq:rel_6a}
\II \sim T_{\fe,\fu^2}(\bJ)\ \sim\ T_{\fe,\fu}(\bJ)\, T_{\fe,\fu}(\fu\bJ)\,,
\end{equation}
where $\II$ is the identity matrix.
The transition matrix  $T_{\fe,\fu}(\bJ)$ is computed using \eqref{eq:pru1}, 
\eqref{eq:Tes} and the expression $\fu=\fs\ft\fs\ft^2\fs$:
\begin{equation}\label{eq:overlap_u_v0}
T_{\fe,\fu}(\bJ)\sim \,P\,T_{\fe,\ft^2}(\fs\bJ)\,P\,T_{\fe,\ft}(\fs\ft^2\fs\bJ)\,P\,,
\end{equation}
where $P$ is the permutation matrix introduced in \eqref{eq:Tes}.
From \eqref{eq:pru1}, we also get
\begin{equation}\label{eq:overlap_t2}
T_{\fe,\ft^2}(\bJ)\ \sim\ T_{\fe,\ft}(\bJ)\, T_{\fe,\ft}(\ft\bJ)\,.
\end{equation}
Gathering \eqref{eq:overlap_u_v0} and \eqref{eq:overlap_t2}, we obtain
\begin{equation}\label{eq:overlap_u}
T_{\fe,\fu}(\bJ)\sim \,P\,T_{\fe,\ft}(\fs\bJ)\,T_{\fe,\ft}(\ft\fs\bJ)\,P\,T_{\fe,\ft}(\fs\ft^2\fs\bJ)\,P\,.
\end{equation}
Then, \eqref{eq:rel_6a} can be rewritten as 
\begin{equation}\label{eq:rel_6b}
\II \,\sim\,P\,T_{\fe,\ft}(\fs\bJ)\,T_{\fe,\ft}(\ft\fs\bJ)\,P\,T_{\fe,\ft}(\fs\ft^2\fs\bJ)
\,T_{\fe,\ft}(\fs\fu\bJ)\,T_{\fe,\ft}(\ft\fs\fu\bJ)\,P\,T_{\fe,\ft}(\fs\ft^2\fs\fu\bJ)\,P\,.
\end{equation}
Relation \eqref{eq:rel_6b}, once expressed in terms of Racah polynomials using expression 
 \eqref{eq:racah_poly}, leads to an identity implying the product of 6 Racah polynomials.
As mentioned at the beginning of the section, this relation can be deduced from the unitarity, Racah and 
Biedenharn--Elliott relations. 
Indeed, we can rewrite the r.h.s of relation \eqref{eq:rel_6b} as 
\begin{equation}
\begin{aligned}
&P\,T_{\fe,\ft}(\fs\bJ)\,T_{\fe,\ft}(\ft\fs\bJ)\,P\,\Big(T_{\fe,\ft}(\bJ')
\,T_{\fe,\ft}(\ft\bJ')\,T_{\fe,\ft}(\ft^2\bJ')\Big)\,P\,T_{\fe,\ft}(\ft^2\fs\bJ)\,P
\\
&\sim\ P\,T_{\fe,\ft}(\fs\bJ)\,T_{\fe,\ft}(\ft\fs\bJ)\,P^2\,T_{\fe,\ft}(\ft^2\fs\bJ)\,P
\ \sim\ P\,T_{\fe,\ft}(\fs\bJ)\,T_{\fe,\ft}(\ft\fs\bJ)\,T_{\fe,\ft}(\ft^2\fs\bJ)\,P\ \sim\ P^2\ \sim\ \II\,,
\end{aligned}
\end{equation}
where $\boldsymbol{J'}=st^2s\boldsymbol{J}$ and we used twice the Racah relation \eqref{eq:t3} to get the second line.

\section{Bivariate Racah polynomials \label{sec:bivar}}

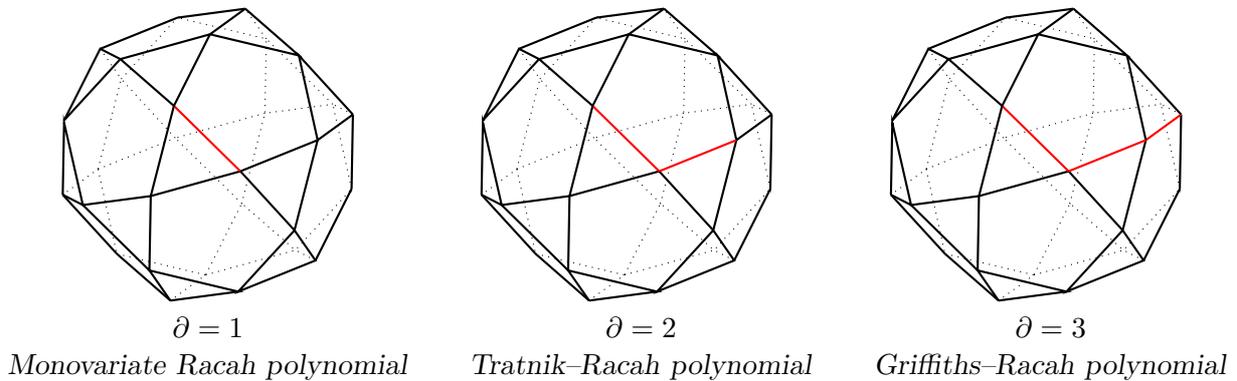
\begin{figure}[htbp]
\begin{center}
\begin{tabular}{ccc}
\begin{tikzpicture}[scale=1]
\draw [thick] (-1.908111787, -.5319769439)--(-1.890410045, .4242268192)-- (-1.643196449, -.6682790668) -- (-1.908111787, -.5319769439);

\draw [dotted]  (-1.908111787, -.5319769439)--(-1.196979681, -1.314983596)--(-1.444193277, -.1924777098) -- (-1.908111787, -.5319769439);

\draw [thick] (-1.908111787, -.5319769439)--(-1.196979681, -1.314983596);

\draw [dotted] (-1.196979681, -1.314983596)--(-.4925605311, -1.935210443)--(-0.0286420207, -1.595711209) --(-1.196979681, -1.314983596);

\draw [thick] (-1.196979681, -1.314983596)--(-.4925605311, -1.935210443);

\draw  [thick] (-.4925605311, -1.935210443)--(-.7683376602, -1.535525063)--(.3999999999, -1.816252676) -- (-.4925605311, -1.935210443);

\draw  [thick] (-.7683376602, -1.535525063)--(-1.643196449, -.6682790668)--(-.7506359181, -.5493212998) --  (-.7683376602, -1.535525063);

\draw  [dotted] (.7683376602, 1.535525063)--(.4925605311, 1.935210443)--(-.3999999999, 1.816252676) --(.7683376602, 1.535525063);

\draw [thick] (.4925605311, 1.935210443)--(-.3999999999, 1.816252676) ;

\draw  [dotted] (1.643196449, .6682790668)--(.7683376602, 1.535525063)--(.7506359181, .5493212998)--(1.643196449, .6682790668);

\draw [dotted] (1.908111787, .5319769439)--(1.643196449, .6682790668)--(1.890410045, -.4542268192)-- (1.908111787, .5319769439);

\draw [thick] (1.890410045, -.4542268192)-- (1.908111787, .5319769439);

\draw [ thick]  (1.908111787, .5319769439)--(1.444193277, .1924777098)--(1.196979681, 1.314983596) -- (1.908111787, .5319769439);

\draw [thick]  (.4925605311, 1.935210443)--(1.196979681, 1.314983596)--(0.0286420207, 1.595711209) -- (.4925605311, 1.935210443);

\draw [dotted]  (-1.139774127, 1.003548119)--(-.3999999999, 1.816252676)--(-1.415551256, 1.403233499) -- (-1.139774127, 1.003548119);

\draw [thick] (-.3999999999, 1.816252676)--(-1.415551256, 1.403233499);

\draw [dotted] (-1.444193277, -.1924777098)--(-.4286420207, .2205414672)--(-1.139774127, 1.003548119)-- (-1.444193277, -.1924777098);

\draw [dotted]  (.7506359181, .5493212998)--(-.4286420207, .2205414672)--(.4462167681, -.6467045291) --(.7506359181, .5493212998);

\draw  [dotted] (-0.286420207e-1, -1.595711209)--(1.150635918, -1.266931376)--(.4462167681, -.6467045291) --(-0.286420207e-1, -1.595711209);

\draw  [dotted] (1.415551256, -1.403233499)--(1.890410045, -.4542268192)--(1.150635918, -1.266931376) -- (1.415551256, -1.403233499);

\draw [thick] (1.415551256, -1.403233499)--(1.890410045, -.4542268192);

\draw  [thick] (1.139774127, -1.003548119)--(1.415551256, -1.403233499)--(.3999999999, -1.816252676) -- (1.139774127, -1.003548119);

\draw  [thick] (1.444193277, .1924777098)--(1.139774127, -1.003548119)--(.4286420207, -.2205414672) --  (1.444193277, .1924777098);

\draw [thick]  (.4286420207, -.2205414672)--(-.7506359181, -.5493212998);
\draw [thick] (-.7506359181, -.5493212998)--(-.4462167681, .6467045291) ;
\draw [thick, color=red]  (-.4462167681, .6467045291) -- (.4286420207, -.2205414672);

\draw [ thick]  (0.286420207e-1, 1.595711209)--(-.4462167681, .6467045291)--(-1.150635918, 1.266931376) -- (0.286420207e-1, 1.595711209);

\draw [ thick] (-1.890410045, .4542268192) -- (-1.415551256, 1.403233499)--(-1.150635918, 1.266931376) --(-1.890410045, .4542268192);

\end{tikzpicture}
&
\begin{tikzpicture}[scale=1]
\draw [thick] (-1.908111787, -.5319769439)--(-1.890410045, .4242268192)-- (-1.643196449, -.6682790668) -- (-1.908111787, -.5319769439);

\draw [dotted]  (-1.908111787, -.5319769439)--(-1.196979681, -1.314983596)--(-1.444193277, -.1924777098) -- (-1.908111787, -.5319769439);

\draw [thick] (-1.908111787, -.5319769439)--(-1.196979681, -1.314983596);

\draw [dotted] (-1.196979681, -1.314983596)--(-.4925605311, -1.935210443)--(-0.0286420207, -1.595711209) --(-1.196979681, -1.314983596);

\draw [thick] (-1.196979681, -1.314983596)--(-.4925605311, -1.935210443);

\draw  [thick] (-.4925605311, -1.935210443)--(-.7683376602, -1.535525063)--(.3999999999, -1.816252676) -- (-.4925605311, -1.935210443);

\draw  [thick] (-.7683376602, -1.535525063)--(-1.643196449, -.6682790668)--(-.7506359181, -.5493212998) --  (-.7683376602, -1.535525063);

\draw  [dotted] (.7683376602, 1.535525063)--(.4925605311, 1.935210443)--(-.3999999999, 1.816252676) --(.7683376602, 1.535525063);

\draw [thick] (.4925605311, 1.935210443)--(-.3999999999, 1.816252676) ;

\draw  [dotted] (1.643196449, .6682790668)--(.7683376602, 1.535525063)--(.7506359181, .5493212998)--(1.643196449, .6682790668);

\draw [dotted] (1.908111787, .5319769439)--(1.643196449, .6682790668)--(1.890410045, -.4542268192)-- (1.908111787, .5319769439);

\draw [thick] (1.890410045, -.4542268192)-- (1.908111787, .5319769439);

\draw [ thick]  (1.908111787, .5319769439)--(1.444193277, .1924777098)--(1.196979681, 1.314983596) -- (1.908111787, .5319769439);

\draw [thick]  (.4925605311, 1.935210443)--(1.196979681, 1.314983596)--(0.0286420207, 1.595711209) -- (.4925605311, 1.935210443);

\draw [dotted]  (-1.139774127, 1.003548119)--(-.3999999999, 1.816252676)--(-1.415551256, 1.403233499) -- (-1.139774127, 1.003548119);

\draw [thick] (-.3999999999, 1.816252676)--(-1.415551256, 1.403233499);

\draw [dotted] (-1.444193277, -.1924777098)--(-.4286420207, .2205414672)--(-1.139774127, 1.003548119)-- (-1.444193277, -.1924777098);

\draw [dotted]  (.7506359181, .5493212998)--(-.4286420207, .2205414672)--(.4462167681, -.6467045291) --(.7506359181, .5493212998);

\draw  [dotted] (-0.286420207e-1, -1.595711209)--(1.150635918, -1.266931376)--(.4462167681, -.6467045291) --(-0.286420207e-1, -1.595711209);

\draw  [dotted] (1.415551256, -1.403233499)--(1.890410045, -.4542268192)--(1.150635918, -1.266931376) -- (1.415551256, -1.403233499);

\draw [thick] (1.415551256, -1.403233499)--(1.890410045, -.4542268192);

\draw  [thick] (1.139774127, -1.003548119)--(1.415551256, -1.403233499)--(.3999999999, -1.816252676) -- (1.139774127, -1.003548119);

\draw  [thick] (1.444193277, .1924777098)--(1.139774127, -1.003548119);
\draw  [thick] (1.139774127, -1.003548119)--(.4286420207, -.2205414672) ;
\draw  [thick, color=red] (.4286420207, -.2205414672) --  (1.444193277, .1924777098);

\draw [thick]  (.4286420207, -.2205414672)--(-.7506359181, -.5493212998);
\draw [thick] (-.7506359181, -.5493212998)--(-.4462167681, .6467045291) ;
\draw [thick, color=red]  (-.4462167681, .6467045291) -- (.4286420207, -.2205414672);

\draw [ thick]  (0.286420207e-1, 1.595711209)--(-.4462167681, .6467045291)--(-1.150635918, 1.266931376) -- (0.286420207e-1, 1.595711209);

\draw [ thick] (-1.890410045, .4542268192) -- (-1.415551256, 1.403233499)--(-1.150635918, 1.266931376) --(-1.890410045, .4542268192);

\end{tikzpicture}
&
\begin{tikzpicture}[scale=1]
\draw [thick] (-1.908111787, -.5319769439)--(-1.890410045, .4242268192)-- (-1.643196449, -.6682790668) -- (-1.908111787, -.5319769439);

\draw [dotted]  (-1.908111787, -.5319769439)--(-1.196979681, -1.314983596)--(-1.444193277, -.1924777098) -- (-1.908111787, -.5319769439);

\draw [thick] (-1.908111787, -.5319769439)--(-1.196979681, -1.314983596);

\draw [dotted] (-1.196979681, -1.314983596)--(-.4925605311, -1.935210443)--(-0.0286420207, -1.595711209) --(-1.196979681, -1.314983596);

\draw [thick] (-1.196979681, -1.314983596)--(-.4925605311, -1.935210443);

\draw  [thick] (-.4925605311, -1.935210443)--(-.7683376602, -1.535525063)--(.3999999999, -1.816252676) -- (-.4925605311, -1.935210443);

\draw  [thick] (-.7683376602, -1.535525063)--(-1.643196449, -.6682790668)--(-.7506359181, -.5493212998) --  (-.7683376602, -1.535525063);

\draw  [dotted] (.7683376602, 1.535525063)--(.4925605311, 1.935210443)--(-.3999999999, 1.816252676) --(.7683376602, 1.535525063);

\draw [thick] (.4925605311, 1.935210443)--(-.3999999999, 1.816252676) ;

\draw  [dotted] (1.643196449, .6682790668)--(.7683376602, 1.535525063)--(.7506359181, .5493212998)--(1.643196449, .6682790668);

\draw [dotted] (1.908111787, .5319769439)--(1.643196449, .6682790668)--(1.890410045, -.4542268192)-- (1.908111787, .5319769439);

\draw [thick] (1.890410045, -.4542268192)-- (1.908111787, .5319769439);

\draw [ thick, color=red]  (1.908111787, .5319769439)--(1.444193277, .1924777098);
\draw [ thick]  (1.444193277, .1924777098)--(1.196979681, 1.314983596);
\draw [ thick]  (1.196979681, 1.314983596) -- (1.908111787, .5319769439);

\draw [thick]  (.4925605311, 1.935210443)--(1.196979681, 1.314983596)--(0.0286420207, 1.595711209) -- (.4925605311, 1.935210443);

\draw [dotted]  (-1.139774127, 1.003548119)--(-.3999999999, 1.816252676)--(-1.415551256, 1.403233499) -- (-1.139774127, 1.003548119);

\draw [thick] (-.3999999999, 1.816252676)--(-1.415551256, 1.403233499);

\draw [dotted] (-1.444193277, -.1924777098)--(-.4286420207, .2205414672)--(-1.139774127, 1.003548119)-- (-1.444193277, -.1924777098);

\draw [dotted]  (.7506359181, .5493212998)--(-.4286420207, .2205414672)--(.4462167681, -.6467045291) --(.7506359181, .5493212998);

\draw  [dotted] (-0.286420207e-1, -1.595711209)--(1.150635918, -1.266931376)--(.4462167681, -.6467045291) --(-0.286420207e-1, -1.595711209);

\draw  [dotted] (1.415551256, -1.403233499)--(1.890410045, -.4542268192)--(1.150635918, -1.266931376) -- (1.415551256, -1.403233499);

\draw [thick] (1.415551256, -1.403233499)--(1.890410045, -.4542268192);

\draw  [thick] (1.139774127, -1.003548119)--(1.415551256, -1.403233499)--(.3999999999, -1.816252676) -- (1.139774127, -1.003548119);

\draw  [thick] (1.444193277, .1924777098)--(1.139774127, -1.003548119);
\draw  [thick] (1.139774127, -1.003548119)--(.4286420207, -.2205414672) ;
\draw  [thick, color=red] (.4286420207, -.2205414672) --  (1.444193277, .1924777098);

\draw [thick]  (.4286420207, -.2205414672)--(-.7506359181, -.5493212998);
\draw [thick] (-.7506359181, -.5493212998)--(-.4462167681, .6467045291) ;
\draw [thick, color=red]  (-.4462167681, .6467045291) -- (.4286420207, -.2205414672);

\draw [ thick]  (0.286420207e-1, 1.595711209)--(-.4462167681, .6467045291)--(-1.150635918, 1.266931376) -- (0.286420207e-1, 1.595711209);

\draw [ thick] (-1.890410045, .4542268192) -- (-1.415551256, 1.403233499)--(-1.150635918, 1.266931376) --(-1.890410045, .4542268192);

\end{tikzpicture}
\\
\textsl{$\partial=1$}\quad &\quad \textsl{$\partial=2$}\quad &\quad \textsl{$\partial=3$}
\\
\textsl{Monovariate Racah polynomial}\quad &\quad \textsl{Tratnik--Racah polynomial}\quad &\quad \textsl{Griffiths--Racah polynomial}
\end{tabular}
\end{center}
\caption{Paths of different length on the folded icosidodecahedron. \label{fig4}}
\end{figure}

In the previous section, we have shown how the cycles on the icosidodecahedron allow us to obtain properties of monovariate Racah polynomials.
Here, we study paths on the icosidodecahedron with different lengths $\partial$. We will consider the paths illustrated in red in figure \ref{fig4}.

We shall show that the different transition matrices associated to these paths (see section \ref{sec:OCI}) for $\partial=2$ and $\partial=3$ 
are associated to multivariate polynomials of Racah type.

The path of length one is associated to the usual monovariate Racah polynomial. Indeed, gathering the computations done 
in sections \ref{sec:equiv} and \ref{sec:pmr}, the properties of the monovariate Racah polynomials can be recast in the following way.
Equation \eqref{eq4} for $\fh=\fe$, $\fg=\fr$ and $X=C_{12}$ or $C_{23}$
becomes
\begin{eqnarray}
T_{\fe,\fr}(\bJ) \pi_{\fr}(C_{23})= \pi_{\fe}(C_{23})T_{\fe,\fr}(\bJ)\,, \qquad
T_{\fe,\fr}(\bJ) \pi_{\fr}(C_{12})= \pi_{\fe}(C_{12})T_{\fe,\fr}(\bJ)\,.
\end{eqnarray}
These relations correspond to the recurrence and difference relations of the monovariate Racah polynomials as used previously.
The orthogonality relation of the transition matrix 
\begin{eqnarray}
T_{\fe,\fr}(\bJ) T_{\fe,\fr}(\bJ)^t  = \II  
\end{eqnarray}
is associated to the unitarity relation of the monovariate Racah polynomials.

\subsection{Tratnik polynomials\label{sec:Tratnik}}

Let us define the function, for $n_1,n_2,m_1,m_2\geq 0$, $n_1+n_2\leq N$, $m_1+m_2\leq N$ and $n_2+m_1\leq N$,
\begin{eqnarray}
 \cT_{n_1,n_2}(m_1,m_2)\!\!\!&=&\!\!\!\cP_{n_1}(m_1;-2j^{(2)}-1,-2j^{(1)}-1 ,n_2-N-1,N-n_2-2j^{(2)}+2j^{(3)}+1 )\nonumber\\
 &\times& \!\!\! \cP_{n_2}(m_2;-2j^{(0)}-1,-2j^{(4)}-1 ,m_1-N-1,N-m_1-2j^{(0)}-2j^{(1)}-1 ),\qquad \label{def:trat}
\end{eqnarray}
where $\cT_{n_1,n_2}(m_1,m_2)$ stands for $\cT_{n_1,n_2}(m_1,m_2;N;\bJ)$. We recall that $\bJ=(j^{(1)},j^{(2)},j^{(3)},j^{(4)},j^{(0)} )$ and $N=j^{(1)}+j^{(2)}-j^{(3)}+j^{(4)}+j^{(0)}$. We recall that the function $\cP$, defined in \eqref{eq:qR}, is proportional to the Racah polynomial.
The rest of this section is devoted to provide different properties of this function as its recurrence and 
difference relations as well as its connection with a transition matrix. We recover the results on Tratnik polynomials \cite{Tra,GI,BV,DIVV}.

\paragraph{Link with the transition matrix $T_{\fe,\fr^2}(\bJ)$.}
Let us focus on the transition matrix associated to the path of length $\partial=2$, as shown on figure \ref{fig4}. 
As explained in section \ref{sec:OCI}, it corresponds to the change of basis between the representations where $(C_{12},C_{123})$ or $(C_{23},C_{234})$ are diagonal. 
Using expression \eqref{eq:gener} of the automorphism $\fr$, one remarks that $(\fr^{-2}(C_{12}),\fr^{-2}(C_{123}))=(C_{23},C_{234})$.
Therefore this path corresponds to the transition matrix $T_{\fe,\fr^2}(\bJ)$. Then, using expression \eqref{eq:pru1} and \eqref{eq:Tpath}, one gets
\begin{equation}
 T_{\fe,\fr^2}(\bJ)\sim T_{\fe,\fr}(\bJ)T_{\fe,\fr}(\fr \bJ)\ .
\end{equation}
Finally, using the explicit expression of $T_{\fe,\fr}$ given in \eqref{eq:overlap_r} and of $\fr \bJ$ given in \eqref{eq:rJ}, one obtains
\begin{equation}
 [T_{\fe,\fr^2}(\bJ)]_{n_1n_2}^{m_1m_2}\sim \cT_{n_1,n_2}(m_1,m_2) \ \Theta(n_2+m_1\leq N),
\end{equation}
where $\Theta$ is the test function which is equal to one if the condition in its argument is true and vanishes otherwise.

\paragraph{Link with the Tratnik polynomials.} 

Given the expression \eqref{eq:qR} of the functions $\cP$, the function $\cT$ reads
\begin{eqnarray}
 \cT_{n_1,n_2}(m_1,m_2)\!\!\!&=&\!\!\!\cA \ r_{n_1}(m_1;-2j^{(2)}-1,-2j^{(1)}-1 ,n_2-N-1,N-n_2-2j^{(2)}+2j^{(3)}+1 )\nonumber\\
 &\times& \!\!\! r_{n_2}(m_2;-2j^{(0)}-1,-2j^{(4)}-1 ,m_1-N-1,N-m_1-2j^{(0)}-2j^{(1)}-1 ),\qquad \label{def:trat0}
\end{eqnarray}
where $\cA$ is a normalization factor given by
\begin{equation}
 \cA = (-1)^{n_1+n_2}\,\sqrt{ \omega_{n_1,n_2}(\bJ) \, \omega_{m_1,n_2}(\fs\fr\bJ) \, \omega_{n_2,m_1}(\fr\bJ)\, \omega_{m_2,m_1}(\fs\fr^2\bJ) } 
\end{equation}
and 
\begin{align}
 \omega_{n_1,n_2}(\bJ)=&\frac{ (2n_1-2j^{(1)}-2j^{(2)}-1)\,  (-2j^{(1)}+n_1)_{N-n_1-n_2}  }   
 {(-2j^{(1)}-2j^{(2)}-1+n_1)_{N-n_1-n_2+1}} \\
 &\times\frac{(-2j^{(2)})_{n_1}\, (N-n_2-2j^{(1)}-2j^{(2)}+2j^{(3)}+1)_{n_1}\, (n_2-N)_{n_1} }
 {n_1!\, (N-n_2-2j^{(1)}-2j^{(2)})_{n_1}\, (n_2-N-2j^{(3)}-1)_{n_1}}.\nonumber
\end{align}
We recall that $\fs\fr\bJ=(-j^{(3)}-1,j^{(2)},-j^{(1)}-1,j^{(0)},j^{(4)})$, $\fr\bJ=(j^{(4)},j^{(0)},-j^{(1)}-1,j^{(2)},-j^{(3)}-1)$ and $\fs\fr^2\bJ=(j^{(1)},j^{(0)},-j^{(4)}-1,-j^{(3)}-1,j^{(2)})$.

Relation \eqref{def:trat0} shows that, up to a global normalization, $\cT_{n_1,n_2}(m_1,m_2)$ defines a polynomial in $m_1$ and $m_2$.
To recognize the usual definition of the Tratnik polynomials, let us relabel $n_1,m_1,n_2,m_2$ as  $x_1,k_1,N-x_2,N-k_1-k_2$.
With this relabeling, the inequalities satisfied by $n$ and $m$ read:
\begin{equation}
0\leq x_1 \leq x_2 \leq N \ , \qquad k_1,k_2\geq 0, \qquad k_1+k_2 \leq N \ .
\end{equation}
Then, using \eqref{eq:dua} to transform the first function $r$ and \eqref{eq:wh}, \eqref{eq:dua} followed by \eqref{eq:wh} to transform the second one,
relation \eqref{def:trat0} can be rewritten as follows
\begin{align}
 &\cT_{x_1,N-x_2}(k_1,N-k_1-k_2)= \widetilde{\cA}\ r_{k_1}(x_1;-2j^{(2)}-1,2j^{(3)}+1 ,-x_2-1,x_2-2j^{(1)}-2j^{(2)}-1 )\nonumber\\
 \times&r_{k_2}(x_2-k_1;2k_1-2j^{(2)}+2j^{(3)}+1,-2j^{(4)}-1  ,k_1-N-1,k_1+N-2j^{(1)}-2j^{(2)}+2j^{(3)}+1)\qquad\qquad \label{def:trat2}
\end{align}
with
\begin{equation}
 \widetilde{\cA}=\frac{(k_1-N+2j^{(1)}+1)_{N-x_2}\, (-2j^{(4)})_{N-x_2}\, (k_1+2j^{(0)}-N+1)_{k_2}\, (2k_1-2j^{(2)}+2j^{(3)}+2)_{k_2}  }
 {(N-k_1-2j^{(0)}-2j^{(1)}-2j^{(4)}-1 )_{N-x_2}\, (-2j^{(0)})_{N-x_2} \, (k_1-N+2j^{(1)}+1)_{k_2}\, (-2j^{(4)})_{k_2} }  \, 
\cA.\qquad
\end{equation}
Note that in order to get \eqref{def:trat2}, we have used the relation $N=j^{(1)}+j^{(2)}-j^{(3)}+j^{(4)}+j^{(0)}$.
Finally, by setting $j^{(1)}=-\frac{1}{2}(\beta_0+1)$, $j^{(2)}=\frac{1}{2}(\beta_0-\beta_1)$, $j^{(3)}=\frac{1}{2}(\beta_2-\beta_1-2)$, $j^{(4)}=\frac{1}{2}(\beta_2-\beta_3)$, 
the above expression becomes 
\begin{align}
& \hspace{-4mm}\cT_{x_1,N-x_2}(k_1,N-k_1-k_2)=\widetilde{\cA}\ r_{k_1}(x_1;\beta_1-\beta_0-1,\beta_2-\beta_1-1 ,-x_2-1,x_2+\beta_1 )\nonumber\\
  &\hspace{3.2cm}\times r_{k_2}(x_2-k_1;2k_1+\beta_2-\beta_0-1,\beta_3-\beta_2-1 ,k_1-N-1,k_1+N+\beta_2). \label{def:trat3} 
\end{align}
In this last expression, we recognize the expression $R_2$ of the Tratnik polynomials \cite{Tra,GI,BV,DIVV} (see for instance 
relation\footnote{Beware that our definition for $r$ differs from the one of \cite{Tra} by a normalization factor.} (3.10) with $p=2$ in \cite{Tra}):
\begin{equation}
\cT_{x_1,N-x_2}(k_1,N-k_1-k_2)=(-1)^{N+x_1-x_2}\,\sqrt{\cW(x_1,x_2)\,\cK(k_1,k_2)}\, R_2(k_1,k_2;x_1,x_2 )
\end{equation} 
with
\begin{align}
\cW(x_1,x_2) &= \frac{(2x_1-2j^{(1)}-2j^{(2)}-1)(2x_2-2j^{(1)}-2j^{(2)}+2j^{(3)}+1)}{x_1!\,(x_2-x_1)!\,(N-x_2)!}
\nonumber\\
&\quad\times\ \frac{\Gamma\big( 2j^{(1)}+1-x_1 \big)\,\Gamma\big( 2j^{(1)}+2j^{(2)}+1-x_1-x_2 \big)\,\Gamma\big( x_2-x_1+2j^{(3)}+2 \big)}{\Gamma\big( 2j^{(2)}+1-x_1 \big)\,\Gamma\big( 2j^{(1)}+2j^{(2)}+2-x_1 \big)} \nonumber\\
&\quad\times\ \frac{\Gamma\big( x_1+x_2-2j^{(1)}-2j^{(2)}+2j^{(3)}+1 \big)\,\Gamma\big( x_2-N+2j^{(0)}+1 \big)}{\Gamma\big( N+x_2-2j^{(1)}-2j^{(2)}+2j^{(3)}+2 \big)\,\Gamma\big( x_2-N+2j^{(4)}+1 \big)}
\,,\label{eq:weight}\\
\cK(k_1,k_2) &= \frac{(N-k_1-k_2)!}{k_1!\,k_2!}\, \big(2k_1-2j^{(2)}+2j^{(3)}+1\big)\big(2k_1+2k_2-2j^{(2)}+2j^{(3)}-2j^{(4)}+1\big)
\nonumber\\
&\quad\times\ \frac{\Gamma\big(k_1-2j^{(2)} +2j^{(3)}+1\big)\,\Gamma\big(2j^{(2)}+1-k_1 \big)\,\Gamma\big( 2j^{(4)}+1-k_2\big)}
{\Gamma\big( k_1+2j^{(3)}+2 \big)\,\Gamma\big(k_1+k_2+2j^{(1)}-N+1 \big)\,\Gamma\big(k_1+k_2+2j^{(0)}-N+1 \big)}
\nonumber\\
&\quad\times\ \frac{\Gamma\big( 2k_1+k_2-2j^{(2)}+2j^{(3)}-2j^{(4)}+1\big)}
{\Gamma\big( 2k_1+k_2-2j^{(2)}+2j^{(3)}+2\big)\,\Gamma\big(k_1+k_2+2j^{(1)}+2j^{(0)}-N+1 \big)}\,.
\end{align}
The functions $\cW$ and $\cK$ are written in such a way that all individual factors are well-defined when the parameters $j^{(k)}$ are positive half-integers.

\paragraph{Unitarity relations.}
Since our approach relies on orthogonal transition matrices, the unitarity relations are immediate. Indeed, 
  the entries of the relation
$T_{\fe,\fr^2}(\bJ) T_{\fe,\fr^2}(\bJ)^t  = \II $ provide the following relations
\begin{equation}
\sum_{\atopn{m_1,m_2\geq0}{m_1+m_2\leq N}} \cT_{n_1,n_2}(m_1,m_2)\, \cT_{n'_1,n'_2}(m_1,m_2) = \delta_{n_1,n_1'}\delta_{n_2,n_2'}\,,
\end{equation}
which is interpreted as the unitarity relation for $\cT$. In the same way, starting from $T_{\fe,\fr^2}(\bJ)^t T_{\fe,\fr^2}(\bJ)  = \II $, we obtain
\begin{equation}\label{eq:trat-unit}
\sum_{\atopn{n_1,n_2\geq0}{n_1+n_2\leq N}} \cT_{n_1,n_2}(m_1,m_2)\, \cT_{n_1,n_2}(m'_1,m'_2) = \delta_{m_1,m_1'}\delta_{m_2,m_2'}\,.
\end{equation}
This relation can be recast as 
\begin{equation}\label{eq:trat-unit2}
\sum_{0\leq x_1\leq x_2\leq N} \cW(x_1,x_2)\,R_2(k_1,k_2;x_1,x_2)\, R_2(k'_1,k'_2;x_1,x_2) = \cK(k_1,k_2)^{-1}\,\delta_{k_1,k_1'}\delta_{k_2,k_2'}\,.
\end{equation}
One recognizes in $\cW$ the weight computed in \cite{Tra} (up to the use of the identity $\Gamma(z)\Gamma(1-z)=\frac{\pi}{\sin(\pi z)}$ which is valid when $z$ is not an integer).

\paragraph{Difference relations.}

Equation \eqref{eq4} 
for $\fh=\fe$, $\fg=\fr^2$ and $X=C_{12}$ or $C_{123}$
becomes
\begin{align}
T_{\fe,\fr^2}(\bJ) \pi^{\fr^2\bJ}_{\fe}(C_{234})= \pi^{\bJ}_{\fe}(C_{12})T_{\fe,\fr^2}(\bJ)\,, \qquad
T_{\fe,\fr^2}(\bJ) \pi^{\fr^2\bJ}_{\fe}(C_{34})= \pi^{\bJ}_{\fe}(C_{123})T_{\fe,\fr^2}(\bJ)\,,
\end{align}
where we have used  that $\fr^2(C_{12})=C_{234}$ and $\fr^2(C_{123})=C_{34}$.
We recall that the representation  $\pi_\fe$ is given in theorem \ref{th:repfin}. Then, the previous relations read 
\begin{align}
&\varphi^{V}_{m_1,m_2+1}(\fr^2\bJ)\,\cT_{n_1,n_2}(m_1,m_2+1)+\varphi^{V}_{m_1,m_2}(\fr^2\bJ)\,\cT_{n_1,n_2}(m_1,m_2-1)\nonumber\\
+\,&\varphi^{H}_{m_1+1,m_2}(\fr^2\bJ)\,\cT_{n_1,n_2}(m_1+1,m_2)+\varphi^{H}_{m_1,m_2}(\fr^2\bJ)\,\cT_{n_1,n_2}(m_1-1,m_2)\nonumber\\
+\,&\varphi^{D}_{m_1+1,m_2+1}(\fr^2\bJ)\,\cT_{n_1,n_2}(m_1+1,m_2+1)+\varphi^{D}_{m_1,m_2}(\fr^2\bJ)\,\cT_{n_1,n_2}(m_1-1,m_2-1)\nonumber\\
+\,&\varphi^{A}_{m_1+1,m_2}(\fr^2\bJ)\,\cT_{n_1,n_2}(m_1+1,m_2-1)+\varphi^{A}_{m_1,m_2+1}(\fr^2\bJ)\,\cT_{n_1,n_2}(m_1-1,m_2+1)\nonumber\\
+\,&\varphi^{00}_{m_1m_2}(\fr^2\bJ)\,\cT_{n_1,n_2}(m_1,m_2)\ =\ (n_1-j^{(1)}-j^{(2)}-1)(n_1-j^{(1)}-j^{(2)})\,\cT_{n_1,n_2}(m_1,m_2)
\end{align}
and
\begin{align}
&\rho_{m_1,m_2+1}(\fr^2\bJ)\,\cT_{n_1,n_2}(m_1,m_2+1)+\rho_{m_1,m_2}(\fr^2\bJ)\,\cT_{n_1,n_2}(m_1,m_2-1)\nonumber\\
+&\rho^{00}_{m_1m_2}(\fr^2\bJ)\,\cT_{n_1,n_2}(m_1,m_2)\ =\ (n_2-j^{(0)}-j^{(4)}-1)(n_2-j^{(0)}-j^{(4)})\,\cT_{n_1,n_2}(m_1,m_2)\,.
\end{align}
These two equations can be interpreted as difference equations for $\cT_{n_1,n_2}(m_1,m_2)$.

\paragraph{Recurrence relations.}

Equation \eqref{eq4} 
for $\fh=\fe$, $\fg=\fr^2$ and $X=C_{23}$ or $C_{234}$
becomes
\begin{eqnarray}
T_{\fe,\fr^2}(\bJ) \pi^{\fr^2\bJ}_{\fe}(C_{12})= \pi^{\bJ}_{\fe}(C_{23})T_{\fe,\fr^2}(\bJ)\,, \qquad
T_{\fe,\fr^2}(\bJ) \pi^{\fr^2\bJ}_{\fe}(C_{123})= \pi^{\bJ}_{\fe}(C_{234})T_{\fe,\fr^2}(\bJ)\,,
\end{eqnarray}
where we have used $\fr^2(C_{23})=C_{12}$ and $\fr^2(C_{234})=C_{123}$.
Using the representation $\pi_\fe$, the previous relations lead to 
\begin{align}
 &(m_1-j^{(2)}+j^{(3)}) (m_1-j^{(2)}+j^{(3)}+1)\cT_{n_1,n_2}(m_1,m_2)=\psi^{00}_{n_1n_2}(\bJ)\,\cT_{n_1,n_2}(m_1,m_2)\\
 &\qquad+ \psi_{n_1+1,n_2}(\bJ)\,\cT_{n_1+1,n_2}(m_1,m_2)+ \psi_{n_1,n_2}(\bJ)\,\cT_{n_1-1,n_2}(m_1,m_2),\nonumber
\end{align}
and
\begin{align}
 &(m_2-j^{(0)}-j^{(1)}-1) (m_2-j^{(0)}-j^{(1)})\cT_{n_1,n_2}(m_1,m_2)\ =\ \varphi^{00}_{n_1n_2}(\bJ)\,\cT_{n_1,n_2}(m_1,m_2)\label{eq:recuT2}\\
&\qquad +\varphi^{V}_{n_1,n_2+1}(\bJ)\,\cT_{n_1,n_2+1}(m_1,m_2)+\varphi^{V}_{n_1,n_2}(\bJ)\,\cT_{n_1,n_2-1}(m_1,m_2)\nonumber\\
&\qquad +\varphi^{H}_{n_1+1,n_2}(\bJ)\,\cT_{n_1+1,n_2}(m_1,m_2)+\varphi^{H}_{n_1,n_2}(\bJ)\,\cT_{n_1-1,n_2}(m_1,m_2)\nonumber\\
&\qquad +\varphi^{D}_{n_1+1,n_2+1}(\bJ)\,\cT_{n_1+1,n_2+1}(m_1,m_2)+\varphi^{D}_{n_1,n_2}(\bJ)\,\cT_{n_1-1,n_2-1}(m_1,m_2)\nonumber\\
&\qquad  +\varphi^{A}_{n_1+1,n_2}(\bJ)\,\cT_{n_1+1,n_2-1}(m_1,m_2)+\varphi^{A}_{n_1,n_2+1}(\bJ)\,\cT_{n_1-1,n_2+1}(m_1,m_2).\nonumber
\end{align}

\subsection{Griffiths-like polynomials\label{sec:Griffiths}}

Let us now define the function, for $n_1,n_2,m_1,m_2\geq 0$, $n_1+n_2\leq N$ and $m_1+m_2\leq N$,
\begin{equation}\label{eq:griff0}
\begin{aligned}
 \cG_{n_1,n_2}(m_1,m_2)=\sum_{ a=0}^{M} (-1)^{a+m_2}
& \cP_{n_1}(a; -2j^{(2)}-1,-2j^{(1)}-1 ,n_2-N-1,N-n_2-2j^{(2)}+2j^{(3)}+1 )\\
 \times &\cP_{n_2}(m_2; -2j^{(0)}-1,-2j^{(4)}-1 ,a-N-1,N-a-2j^{(0)}-2j^{(1)}-1 )\\
 \times& \cP_{m_1}(a; -2j^{(2)}-1,-2j^{(4)}-1 ,m_2-N-1,N-m_2-2j^{(2)}+2j^{(3)}+1 )\,,
\end{aligned}
\end{equation}
where $M=\min(N-n_2,N-m_2)$ and $\cG_{n_1,n_2}(m_1,m_2)$ stands for $\cG_{n_1,n_2}(m_1,m_2; N,\bJ)$. 
We have used \eqref{def:trat} to get the last equality in \eqref{eq:griff0}.
As in the previous section, we provide different properties of this function as its recurrence and 
difference relations as well as its connection with a transition matrix. These relations, which derive directly from our construction, 
are obtained for the first time, to the best of our knowledge.  

\paragraph{Expression in term of Tratnik polynomials.} We can rewrite the above expression as
\begin{equation}\label{eq:griffb}
\begin{aligned}
 \cG_{n_1,n_2}(m_1,m_2)=&\\
 \sum_{ a=0}^{M} (-1)^{a+m_2}\,
 \cT_{n_1,n_2}(a,m_2)&
\,\cP_{m_1}(a; -2j^{(2)}-1,-2j^{(4)}-1 ,m_2-N-1,N-m_2-2j^{(2)}+2j^{(3)}+1 )\,.
\end{aligned}
\end{equation}
In the same way, using the duality relation \eqref{eq:dua} on $\cP_{n_2}(m_2)$, we also have
\begin{equation}\label{eq:griffc}
\begin{aligned}
 \cG_{n_1,n_2}(m_1,m_2)=&\\
 \sum_{ a=0}^{M} (-1)^{a+n_2}\,&
\cP_{n_1}(a; -2j^{(2)}-1,-2j^{(1)}-1 ,n_2-N-1,N-n_2-2j^{(2)}+2j^{(3)}+1 )\\
&\times\ \cT_{m_1,m_2}(a,n_2)\Big|_{j^{(1)}\,\leftrightarrow\, j^{(4)}}
\,.
\end{aligned}
\end{equation}
Remark that equations \eqref{eq:griffc} and\eqref{eq:griffb} show that $ \cG_{n_1,n_2}(m_1,m_2)$ is invariant under the transformation
$j^{(1)}\,\leftrightarrow\, j^{(4)}$ and $m_k\leftrightarrow n_k$, $k=1,2$.
\paragraph{Link with the transition matrix $T_{\fe,\ft^2\fr^2}(\bJ)$.}

Let us focus on the transition matrix associated to the path of length $\partial=3$, as shown on figure \ref{fig4}.  
As explained in section \ref{sec:OCI}, it corresponds to the change of basis between the representations where $(C_{12},C_{123})$ or $(C_{24},C_{234})$ are diagonal. 
Using expressions \eqref{eq:genert} and \eqref{eq:gener}  of the automorphisms $\ft$ and $\fr$, one finds  $(\fr^{-2}\ft^{-2}(C_{12}),\fr^{-2}\ft^{-2}(C_{123}))=(C_{24},C_{234})$.
Therefore this path corresponds to the transition matrix $T_{\fe,\ft^2\fr^2}(\bJ)$. Then, using expression \eqref{eq:pru1} and \eqref{eq:Tpath}, one gets
\begin{equation}
 T_{\fe,\ft^2\fr^2}(\bJ)\sim T_{\fe,\fr}(\bJ)T_{\fe,\fr}(\fr \bJ)T_{\fe,\ft^{2}}(\fr^2 \bJ)\ .
\end{equation}
From \eqref{eq:invt}, one has 
\begin{equation}
  T_{\fe,\ft^{2}}(\bJ) \sim T_{\fe,\ft^{-1}}(\bJ) \sim T_{\fe,\ft}(\ft^{-1}\bJ)^t\,.
\end{equation}
Finally, using the explicit expression of $T_{\fe,\ft}$ and $T_{\fe,\fr}$ given in \eqref{eq:overlap_t} and \eqref{eq:overlap_r}, one obtains
\begin{equation}
 [T_{\fe,\ft^2\fr^2}(\bJ)]_{n_1n_2}^{m_1m_2}\sim \cG_{n_1,n_2}(m_1,m_2) .
\end{equation}

As in the previous case, using the definition of the function $\cP$, one obtains in terms of the Racah polynomials $r_n$:
\begin{equation}
\begin{aligned}
 \cG_{n_1,n_2}(m_1,m_2)&=\sum_{ a=0}^{M} (-1)^{a}\,\cB(n_1,n_2,m_1,m_2,a)\\
& \times  r_{n_1}(a; -2j^{(2)}-1,-2j^{(1)}-1 ,n_2-N-1,N-n_2-2j^{(2)}+2j^{(3)}+1 )\\
& \times r_{n_2}(m_2; -2j^{(0)}-1,-2j^{(4)}-1 ,a-N-1,N-a-2j^{(0)}-2j^{(1)}-1 )\\
& \times r_{m_1}(a; -2j^{(2)}-1,-2j^{(4)}-1 ,m_2-N-1,N-m_2-2j^{(2)}+2j^{(3)}+1 )\,,
\end{aligned}
 \label{def:griff0}
\end{equation}
where we recall that $M=\min(N-n_2,N-m_2)$ and $\cB$ is a normalization factor given by
\begin{equation}
 \cB=(-1)^{n_1+n_2+m_1+m_2} \sqrt{ \omega_{n_1,n_2}(\bJ) \, \omega_{a,n_2}(\fs\fr\bJ) \, \omega_{n_2,a}(\fr\bJ) \, \omega_{m_2,a}(\fs\fr^2\bJ) \, \omega_{m_1,m_2}(\ft^2\fr^2\bJ) \, \omega_{a,m_2}(\fs\fr\ft^2\fr^2\bJ)} \,,
\end{equation}
with $\ft^2\fr^2\bJ=(j^{(4)},j^{(2)},j^{(3)},j^{(0)},j^{(1)})$ and $\fs\fr\ft^2\fr^2\bJ=(-j^{(3)}-1,j^{(2)},-j^{(4)}-1,j^{(1)},j^{(0)})$.

Using the results obtained for the Tratnik polynomials, the above expression can be recast as
\begin{equation}
\begin{aligned}
 \cG_{x_1,N-x_2}(m_1,m_2)&=\sum_{ a=0}^{M} \sqrt{\cW(x_1,x_2)\,\cK(a,N-a-m_2) 
 \, \omega_{m_1,m_2}(\ft^2\fr^2\bJ) \, \omega_{a,m_2}(\fs\fr\ft^2\fr^2\bJ)}\\
& \qquad\qquad\qquad \times  (-1)^{N+x_1-x_2+m_1+m_2+a} \, R_2(a,N-a-m_2;x_1,x_2)\, r_{m_1}(a )\,,
\end{aligned}
 \label{def:griff1}
\end{equation}
where $r_{m_1}(a)$ stands for $r_{m_1}(a; -2j^{(2)}-1,-2j^{(4)}-1 ,m_2-N-1,N-m_2-2j^{(2)}+2j^{(3)}+1 )$ and 
now $M=\min(x_2,N-m_2)$. Expression \eqref{def:griff1} shows that, up to a global normalization $\sqrt{\cW(x_1,x_2)}$, 
$\cG_{x_1,N-x_2}(m_1,m_2)$ is a polynomial in $x_1$ and $x_2$.
We now show that this polynomial is orthogonal and obeys recurrence and difference relations.

\paragraph{Unitarity relations.}
Once again, the entries of the relation
$T_{\fe,\ft^2\fr^2}(\bJ) T_{\fe,\ft^2\fr^2}(\bJ)^t  = \II $ provide the unitarity relation for $\cG$
\begin{equation}
\sum_{\atopn{m_1,m_2\geq0}{m_1+m_2\leq N}} \cG_{n_1,n_2}(m_1,m_2)\, \cG_{n'_1,n'_2}(m_1,m_2) = \delta_{n_1,n_1'}\delta_{n_2,n_2'}\,.
\end{equation}
In the same way, starting from $T_{\fe,\fr^2}(\bJ)^t T_{\fe,\fr^2}(\bJ)  = \II $, we obtain
\begin{equation}
\sum_{\atopn{n_1,n_2\geq0}{n_1+n_2\leq N}} \cG_{n_1,n_2}(m_1,m_2)\, \cG_{n_1,n_2}(m'_1,m'_2) = \delta_{m_1,m_1'}\delta_{m_2,m_2'}\,.
\end{equation}

\paragraph{Difference relations.}

Equation \eqref{eq4} 
for $\fh=\fe$, $\fg=\ft^2\fr^2$ and $X=C_{12}$ or $C_{123}$
becomes
\begin{eqnarray}
T_{\fe,\ft^2\fr^2}(\bJ) \pi^{\ft^2\fr^2\bJ}_{\fe}(C_{134})= \pi^{\bJ}_{\fe}(C_{12})T_{\fe,\ft^2\fr^2}(\bJ)\,, \qquad
T_{\fe,\ft^2\fr^2}(\bJ) \pi^{\ft^2\fr^2\bJ}_{\fe}(C_{14})= \pi^{\bJ}_{\fe}(C_{123})T_{\fe,\ft^2\fr^2}(\bJ)\,,
\end{eqnarray}
where the relations $\ft^2\fr^2(C_{12})=C_{134}$ and $\ft^2\fr^2(C_{123})=C_{14}$ have been used. Let us recall 
that $\ft^2\fr^2\bJ=(j^{(4)},j^{(2)},j^{(3)},j^{(0)},j^{(1)})$ and that the representation  $\pi_\fe$ is given in theorem \ref{th:repfin}. Then, the previous relations read 
\begin{align}
&-\widetilde{\varphi}^{V}_{m_1,m_2+1}(\ft^2\fr^2\bJ)\,\cG_{n_1,n_2}(m_1,m_2+1) -\widetilde{\varphi}^{V}_{m_1,m_2}(\ft^2\fr^2\bJ)\,\cG_{n_1,n_2}(m_1,m_2-1)\nonumber\\
&-\varphi^{H}_{m_1+1,m_2}(\ft^2\fr^2\bJ)\,\cG_{n_1,n_2}(m_1+1,m_2)-\varphi^{H}_{m_1,m_2}(\ft^2\fr^2\bJ)\,\cG_{n_1,n_2}(m_1-1,m_2)\nonumber\\
&-\varphi^{D}_{m_1+1,m_2+1}(\ft^2\fr^2\bJ)\,\cG_{n_1,n_2}(m_1+1,m_2+1)-\varphi^{D}_{m_1,m_2}(\ft^2\fr^2\bJ)\,\cG_{n_1,n_2}(m_1-1,m_2-1)\nonumber\\
&-\varphi^{A}_{m_1+1,m_2}(\ft^2\fr^2\bJ)\,\cG_{n_1,n_2}(m_1+1,m_2-1)-\varphi^{A}_{m_1,m_2+1}(\ft^2\fr^2\bJ)\,\cG_{n_1,n_2}(m_1-1,m_2+1)\nonumber\\
&-\varphi^{0\overline{0}}_{m_1m_2}(\ft^2\fr^2\bJ)\,\cG_{n_1,n_2}(m_1,m_2)\ =\ (n_1-j^{(1)}-j^{(2)}-1)(n_1-j^{(1)}-j^{(2)})\,\cG_{n_1,n_2}(m_1,m_2)\,,
\end{align}
and
\begin{align}
&-\varphi^{V}_{m_1,m_2+1}(\ft^2\fr^2\bJ)\,\cG_{n_1,n_2}(m_1,m_2+1)-\varphi^{V}_{m_1,m_2}(\ft^2\fr^2\bJ)\,\cG_{n_1,n_2}(m_1,m_2-1)\nonumber\\
&-\widetilde{\varphi}^{H}_{m_1+1,m_2}(\ft^2\fr^2\bJ)\,\cG_{n_1,n_2}(m_1+1,m_2)-\widetilde{\varphi}^{H}_{m_1,m_2}(\ft^2\fr^2\bJ)\,\cG_{n_1,n_2}(m_1-1,m_2)\nonumber\\
&-\varphi^{D}_{m_1+1,m_2+1}(\ft^2\fr^2\bJ)\,\cG_{n_1,n_2}(m_1+1,m_2+1)-\varphi^{D}_{m_1,m_2}(\ft^2\fr^2\bJ)\,\cG_{n_1,n_2}(m_1-1,m_2-1)\nonumber\\
&-\varphi^{A}_{m_1+1,m_2}(\ft^2\fr^2\bJ)\,\cG_{n_1,n_2}(m_1+1,m_2-1)-\varphi^{A}_{m_1,m_2+1}(\ft^2\fr^2\bJ)\,\cG_{n_1,n_2}(m_1-1,m_2+1)\nonumber\\
&-\varphi^{\overline{0}0}_{m_1m_2}(\ft^2\fr^2\bJ)\,\cG_{n_1,n_2}(m_1,m_2)\ =\ (n_2-j^{(0)}-j^{(4)}-1)(n_2-j^{(0)}-j^{(4)})\,\cG_{n_1,n_2}(m_1,m_2)\,,
\end{align}
where the explicit expressions of $\varphi$ can be found in appendix \ref{app:C}.
These two equations can be interpreted as difference equations for $\cG_{n_1,n_2}(m_1,m_2)$.

\paragraph{Recurrence relations.}

Equation \eqref{eq4} for $\fh=\fe$, $\fg=\ft^2\fr^2$ and $X=C_{24}$ or $C_{234}$ becomes
\begin{eqnarray}
T_{\fe,\ft^2\fr^2}(\bJ) \pi^{\ft^2\fr^2\bJ}_{\fe}(C_{12})= \pi^{\bJ}_{\fe}(C_{24})T_{\fe,\ft^2\fr^2}(\bJ)\,, \qquad
T_{\fe,\ft^2\fr^2}(\bJ) \pi^{\ft^2\fr^2\bJ}_{\fe}(C_{123})= \pi^{\bJ}_{\fe}(C_{234})T_{\fe,\ft^2\fr^2}(\bJ)\,, \label{eq:46t}
\end{eqnarray}
where we used again $\ft^2\fr^2(C_{24})=C_{12}$ and $\ft^2\fr^2(C_{234})=C_{123}$.
The first equation in \eqref{eq:46t} leads to 
\begin{align}
  &(m_1-j^{(2)}-j^{(4)}-1) (m_1-j^{(2)}-j^{(4)})\cG_{n_1,n_2}(m_1,m_2)\ =\ \varphi^{\overline{0}\overline{0}}_{n_1n_2}(\bJ)\,\cG_{n_1,n_2}(m_1,m_2)\\
&\qquad +\widetilde{\varphi}^{V}_{n_1,n_2+1}(\bJ)\,\cG_{n_1,n_2+1}(m_1,m_2)+\widetilde{\varphi}^{V}_{n_1,n_2}(\bJ)\,\cG_{n_1,n_2-1}(m_1,m_2)\nonumber\\
&\qquad +\widetilde{\varphi}^{H}_{n_1+1,n_2}(\bJ)\,\cG_{n_1+1,n_2}(m_1,m_2)+\widetilde{\varphi}^{H}_{n_1,n_2}(\bJ)\,\cG_{n_1-1,n_2}(m_1,m_2)\nonumber\\
&\qquad +\varphi^{D}_{n_1+1,n_2+1}(\bJ)\,\cG_{n_1+1,n_2+1}(m_1,m_2)+\varphi^{D}_{n_1,n_2}(\bJ)\,\cG_{n_1-1,n_2-1}(m_1,m_2)\nonumber\\
&\qquad +\varphi^{A}_{n_1+1,n_2}(\bJ)\,\cG_{n_1+1,n_2-1}(m_1,m_2)+\varphi^{A}_{n_1,n_2+1}(\bJ)\,\cG_{n_1-1,n_2+1}(m_1,m_2),\nonumber
\end{align}
with $\varphi^{\overline{0}\overline{0}}$ defined by \eqref{eq:min3}.
The second equation in \eqref{eq:46t} provides an equation similar to \eqref{eq:recuT2}:
\begin{align}
 &(m_2-j^{(0)}-j^{(1)}-1) (m_2-j^{(0)}-j^{(1)})\cG_{n_1,n_2}(m_1,m_2)\ =\ \varphi^{00}_{n_1n_2}(\bJ)\,\cG_{n_1,n_2}(m_1,m_2)\label{eq:recuG2}\\
&\qquad +\varphi^{V}_{n_1,n_2+1}(\bJ)\,\cG_{n_1,n_2+1}(m_1,m_2)+\varphi^{V}_{n_1,n_2}(\bJ)\,\cG_{n_1,n_2-1}(m_1,m_2)\nonumber\\
&\qquad +\varphi^{H}_{n_1+1,n_2}(\bJ)\,\cG_{n_1+1,n_2}(m_1,m_2)+\varphi^{H}_{n_1,n_2}(\bJ)\,\cG_{n_1-1,n_2}(m_1,m_2)\nonumber\\
&\qquad +\varphi^{D}_{n_1+1,n_2+1}(\bJ)\,\cG_{n_1+1,n_2+1}(m_1,m_2)+\varphi^{D}_{n_1,n_2}(\bJ)\,\cG_{n_1-1,n_2-1}(m_1,m_2)\nonumber\\
&\qquad  +\varphi^{A}_{n_1+1,n_2}(\bJ)\,\cG_{n_1+1,n_2-1}(m_1,m_2)+\varphi^{A}_{n_1,n_2+1}(\bJ)\,\cG_{n_1-1,n_2+1}(m_1,m_2).\nonumber
\end{align}

By contrast to the Tratnik polynomials, we obtain here 9-point relations for both recurrence and difference relations.
In the case of the multivariate Krawtchouk polynomial, the generic Griffiths polynomials has the same kind of complication for the recurrence and difference relations \cite{Griff,HoRa}. 
This motivates the denomination of Griffiths function for $\cG_{n_1,n_2}(m_1,m_2)$.

\section{Conclusion \label{sec:conclu}}

In this paper, the algebra $sR(4)$ and its representations are studied in detail. The knowledge of these representations allows us to provide 
important relations for bivariate functions. Numerous generalizations are possible and will provide important results in different contexts,
such as representation theory, super-integrable models and algebraic combinatorics. We detail some of them in the following.
 
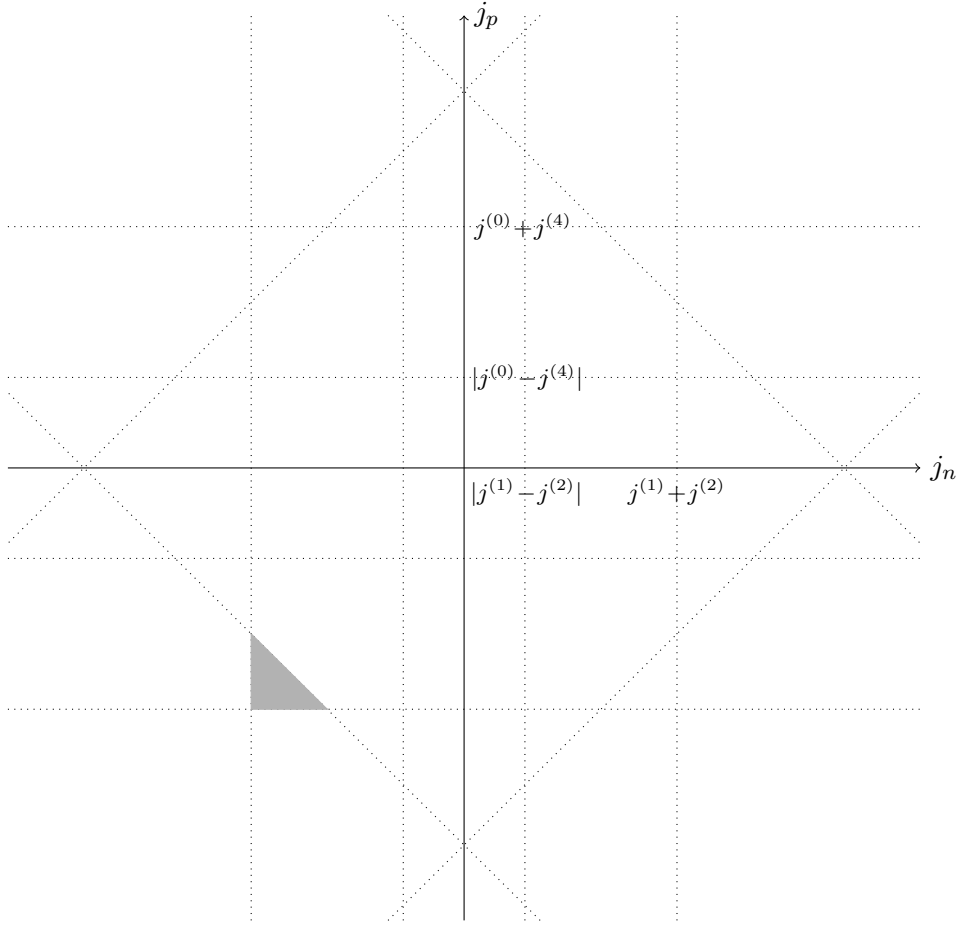
\begin{figure}[htb]
\begin{center}
\begin{tikzpicture}[scale=0.4]
\draw [->]  (-15,0)--(15,0) ;
\draw (15,0) node[right] {$j_n$};
\draw (0,15) node[right] {$j_p$};
\draw [->]  (0,-15)--(0,15);

\draw [dotted]  (2,-15)--(2,15);
\draw [dotted]  (7,-15)--(7,15);
\draw [dotted]  (-2,-15)--(-2,15);
\draw [dotted]  (-7,-15)--(-7,15);

\draw [dotted]  (-15,3)--(15,3) ;
\draw [dotted]  (-15,8)--(15,8) ;
\draw [dotted]  (-15,-3)--(15,-3) ;
\draw [dotted]  (-15,-8)--(15,-8) ;

 \draw [dotted] ( -2.5,15)--(15, -2.5);
\draw [dotted] (2.5,15)--(-15, -2.5);
 \draw [dotted] (-2.5,-15)--(15, 2.5);
 \draw [dotted] (2.5,-15)--(-15, 2.5);
 
 \draw [black!30,fill=black!30] (-7,-5.5)--(-4.5, -8) --(-7,-8) -- (-7,-5.5);
 
 \node[below] at (2.05,0) {\footnotesize$|j^{(1)}\!-\!j^{(2)}|$};
 \node[below] at (7,0) {\footnotesize$j^{(1)}\!+\!j^{(2)}$};
 \node[right] at (-0.1,3) {\footnotesize$|j^{(0)}\!-\!j^{(4)}|$};
 \node[right] at (0,8) {\footnotesize$j^{(0)}\!+\!j^{(4)}$};
 
\end{tikzpicture}

\caption{Regions where the functions \eqref{eq:Y} and \eqref{eq:Y2} keep the same sign. The grey region corresponds to the case studied in this paper.
  \label{region}}
\end{center}
\end{figure}

\paragraph{Other regions.} We have restricted ourselves to the case where the vectors $|n,p\rangle$ of the finite representation studied in section \ref{sec:finite} are 
 such that $n,p\geq 0$ and $n+p\leq N$. 
This choice implies that the functions defined by \eqref{eq:Y} and \eqref{eq:Y2} remain positive. 
In figure \ref{region}, we show schematically where these functions vanish on dashed lines. Therefore, these lines define regions where these functions keep the same sign. 
The grey region corresponds to the region studied in this paper.
We see from this figure that other regions are possible, which lead to some other constraints between $n,p$. A comprehensive investigation of these different cases remains to be done.
Let us remark that for the multivariate Hahn polynomials, 
such a detailed study of the different possible regions has been investigated in \cite{IX,IX2}. 

\paragraph{Classification of the representations.} In the present paper, we have focused on the case 
where the joint spectrum of $C_{12}$ and $C_{123}$ is nondegenerate. It may be interesting to consider the cases when this spectrum is degenerate and 
also the cases when $C_{12}$ or $C_{123}$ are nondiagonalizable. 
We expect that a complete representation theory of $sR(4)$ is possible as it has been done for $R(3)$ in \cite{Huang,HB}. 
We also remarked that the nondiagonal entries of the representation matrices factorize because of the vanishing of the Casimir elements in the special Racah algebra $sR(4)$. 
In the case of $R(4)$ algebra, the values of these Casimir elements provide additional parameters in the representation theory. 
It would be also interesting to study this case in detail.

\paragraph{Higher rank Racah algebra.}

Up to now, we have focused on $R(3)$ and $R(4)$ algebras, but the general case of $R(\sN)$ algebras can be studied following the same lines.
The definition of the special Racah algebra $sR(\sN)$ can be found in \cite{CGPV}. 
We have noticed in the present paper that the defining
 relation  \eqref{eq:comm4} of $R(4)$  \cite{DGVV1, DIVV, BV, CGPV} can be replaced by
  relation  \eqref{eq:comm4bis}. This equivalence is valid  
in the general case of a $R(\sN)$ algebra ($\sN\geq4$). For $\sN\geq5$, supplementary relations are needed, and we believe that they can be replaced by
\begin{equation}\label{eq:comm4bis-N}
\begin{aligned}
&(C_{aj}-C_j-C_a)\,[C_{ik},C_{k\ell}] +(C_{ai}-C_i-C_a)[C_{k\ell},C_{jk}]  +(C_{a\ell}-C_\ell-C_a) [C_{ij},C_{jk}] \\
& \qquad+(C_{ak}-C_k-C_a)[C_{i\ell},C_{j\ell}] 
=0\,.
\end{aligned}
\end{equation}
This replacement should simplify the study of the representation theory of the $sR(\sN)$ algebra.

\paragraph{Automorphism group of $R(\sN)$.}
We proved in this paper that the permutation group $\fS_5$ is an automorphism group of the $sR(4)$ algebra. In the same way, we believe that 
the permutation group $\fS_{\sN+1}$ is an automorphism group for $sR(\sN)$. It can be constructed in the following way.
The permutation group $\fS_5$ acts naturally on the set  $\{0,1,2,3,4\}$: for $g\in \fS_5$, we note $g_{\{i\}}$ the image of $i$. 
For any subset $I\subset\{1,2,3,4\}$ we define $g_I$ by iteration on the cardinality of the subset $I$ as:
\begin{equation}\label{eq:gI}
g_I = \big(g_{\{i\}}\cup g_{I\setminus\{i\}}\big)\setminus \big(g_{\{i\}} \cap g_{I\setminus\{i\}}\big)\,,\quad \forall i\in I.
\end{equation}
This definition is consistent with $g_{\{0\}}$, when identifying $\{0\}$ with $\{1,2,3,4\}$. Then, the relation $C_{g_{I}}=\mathfrak{g}(C_I)$ is compatible 
with the definition of the automorphisms $\fg$ defined by proposition \ref{eq:propiso}. Let us remark that this construction also works for $R(3)$, and leads to an  automorphism group $\fS_4$.
Such a procedure should be also possible in the general case of $R(\sN)$ algebras, leading to an automorphism group $\fS_{\sN+1}$ for $sR(\sN)$.

\paragraph{Connection graph for $sR(\sN)$.} We have shown that the connection graph between abelian subalgebras of $sR(4)$ is associated to an icosidodecahedron, 
which belongs to the class of  Archimedean solids. 
The construction of a connection graph for $sR(\sN)$ is also a point which is worth studying.
 In particular, one can wonder whether it is also an Archimedean polyhedron in $\sN-1$ dimensions. It would give a geometrical interpretation of algebraic results.
 For instance, we have understood the Racah and Biedenharn--Elliott relations as cycles around triangular and pentagonal faces of the icosidodecahedron.

\paragraph{Multivariate Racah polynomials.} Mimicking our approach to the $sR(\sN)$ algebra, the study of its representations should lead to multivariate Racah polynomials
with $N-2$ variables. Depending on the precise structure of the connection graph, there will be different types of multivariate polynomials associated to the 
inequivalent paths of this graph. Starting from equivalent symmetric representations, the transition matrices should provide the orthogonal, recurrence and difference relations for these polynomials.

\paragraph{$q$-deformation.} One may think to a deformed version of the previous results. Indeed, 
a deformation of the Racah algebra is obtained in different contexts, leading to the Askey--Wilson algebras $AW(3)$ \cite{Zhedanov,GZ3},
see \cite{CFGPRV} for a recent review.
Several works have been done to define higher rank algebras $AW(\sN)$ \cite{PW,BCV,Dec}. Recently, the defining relations of $AW(4)$ 
have been introduced in \cite{cooke},
 using the Skein algebra and an
approach based on the centralizer of four copies of $U_q(su(2))$.

The identification of an automorphism group for $AW(\sN)$ algebras would simplify the study of their algebraic structure, as done in \cite{Ter2,CFGPRV}
for $AW(3)$. We expect that it is associated with the braid group on $\sN+1$ strands.
The representation theory of $AW(4)$ may be done similarly to the construction proposed in this paper. One can wonder whether the folded-icosidodecahedron plays also a role.
In the deformed case, the transition matrices will be associated to the Askey--Wilson polynomials and must provide the recurrence and the difference relations for
some multivariate Askey--Wilson polynomials, which generalize the results given in \cite{Iliev,GI2}.

\subsection*{Acknowledgements: } The authors thank warmly J. Gaboriaud, L. Poulain~d'Andecy and L. Vinet for their interest in this work and enlightening discussions.
N. Cramp\'e thanks LAPTh for its hospitality and is supported by the international research project AAPT of the CNRS.
This work is partially supported by Universit\'e Savoie Mont Blanc and Conseil Savoie Mont Blanc grant APOINT.

\begin{appendix}

\section{Presentation of the $R(4)$ algebra\label{app:algR4}}

As mentioned in section \ref{sec:sR4}, the algebra $R(4)$ is generated by the elements $C_{j}$, $j=1,2,3,4$, $C_{1234}$,
$C_{12}$, $C_{23}$, $C_{34}$, $C_{123}$ and $C_{234}$. We will call this basis the \textsl{contiguous basis}.
The remaining elements can be reconstructed using 
the relation \eqref{eq:defCI}:
\begin{equation}\label{eq:Ctrou}
\begin{split}
&C_{13}=-C_{23}+C_{123}-C_{12}+C_1+C_2+C_3\\
&C_{24} = C_{234} - C_{23} - C_{34} + C_2 + C_3 + C_4\,,\\
&C_{14} = C_{1234} - C_{123} - C_{234} + C_{23} + C_1 + C_4\,,\\
&C_{124} =C_{1234} - C_{123} - C_{34} + C_{12} + C_{3} + C_{4}\,,\\
&C_{134} = C_{1234} - C_{234} - C_{12} + C_{34} + C_1 + C_2\,.
\end{split}
\end{equation}
Plugging these expressions in the relations \eqref{eq:comm2}-\eqref{eq:comm3}, we obtain different types of relations:
\paragraph{Commutativity relations ($R(2)$-type relations).} In addition to the centrality of the generators $C_{j}$, $j=1,2,3,4$, $C_{1234}$, 
we get 5 relations:
\begin{equation}
[C_{12}\,,\,C_{34}]=0\,,\quad [C_{12}\,,\,C_{123}]=0\,,\quad [C_{23}\,,\,C_{123}]=0\,,\quad [C_{23}\,,\,C_{234}]=0\,,\quad 
[C_{34}\,,\,C_{234}]=0\,.
\end{equation}
We call them $R(2)$-type relations because $R(2)$ is abelian.

\paragraph{$R(3)$-type relations.} We get 10 of them:
\begin{align}
\label{eq:Rac1} 
&\frac{1}{2}\big[C_{12}, [C_{12},C_{23}] \big] = C_{12}^2 + \{C_{12},C_{23}\} - (C_{1}+C_{2}+C_{3}+C_{123})C_{12} - (C_{1}-C_{2})(C_{3}-C_{123}) \,, \\
\label{eq:Rac2} 
&\frac{1}{2}\big[C_{23}, [C_{23},C_{12}] \big] = C_{23}^2 + \{C_{12},C_{23}\} - (C_{1}+C_{2}+C_{3}+C_{123})C_{23} - (C_{1}-C_{123})(C_{3}-C_{2}) \,,\\
\label{eq:Rac23-34} 
&\frac{1}{2}\big[C_{23}, [C_{23},C_{34}] \big] = C_{23}^2 + \{C_{23},C_{34}\} - (C_{2}+C_{3}+C_{4}+C_{234})C_{23} - (C_{2}-C_{3})(C_{4}-C_{234}) \,,\\
&\frac{1}{2}\big[C_{34}, [C_{34},C_{23}] \big] = C_{34}^2 + \{C_{23},C_{34}\} - (C_{2}+C_{3}+C_{4}+C_{234})C_{34} - (C_{2}-C_{234})(C_{4}-C_{3}) \,,\\
\label{eq:Rac12-234} 
&\frac{1}{2}\big[C_{12}, [C_{12},C_{234}] \big] = C_{12}^2 + \{C_{12},C_{234}\} - (C_{1}+C_{2}+C_{34}+C_{1234})C_{12} - (C_{1}-C_{2})(C_{34}-C_{1234}) \,,\\
&\frac{1}{2}\big[C_{234}, [C_{234},C_{12}] \big] = C_{234}^2 + \{C_{12},C_{234}\} - (C_{1}+C_{2}+C_{34}+C_{1234})C_{234} - (C_{1}-C_{1234})(C_{34}-C_{2}) \,,\\
\label{eq:Rac34-123} 
&\frac{1}{2}\big[C_{34}, [C_{34},C_{123}] \big] = C_{34}^2 + \{C_{123},C_{34}\} - (C_{12}+C_{3}+C_{4}+C_{1234})C_{34} - (C_{12}-C_{1234})(C_{4}-C_{3}) \,,\\
\label{eq:Rac123-34} 
&\frac{1}{2}\big[C_{123}, [C_{123},C_{34}] \big] = C_{123}^2 + \{C_{123},C_{34}\} - (C_{12}+C_{3}+C_{4}+C_{1234})C_{123} - (C_{12}-C_{3})(C_{4}-C_{1234}) \,, \\
&\frac{1}{2}\big[C_{234}, [C_{234},C_{123}] \big] = C_{234}^2 + \{C_{123},C_{234}\} - (C_{1}+C_{23}+C_{4}+C_{1234})C_{234} - (C_{1}-C_{1234})(C_{4}-C_{23}) \,,\\
\intertext{}
\label{eq:Rac123-234} 
&\frac{1}{2}\big[C_{123}, [C_{123},C_{234}] \big] = C_{123}^2 + \{C_{123},C_{234}\} - (C_{1}+C_{23}+C_{4}+C_{1234})C_{123} - (C_{1}-C_{23})(C_{4}-C_{1234}) \,.
\end{align}
This kind of relations already appears in the $R(3)$ algebra.

\paragraph{$R(4)$-type relations.} These types of relations do not exist in the $R(3)$ algebra. They correspond to the relations  
\eqref{eq:comm3bis} and \eqref{eq:comm4bis}. 
In the contiguous basis, the set of relations   \eqref{eq:comm3bis} is equivalent to the five following relations
\begin{eqnarray}
  [C_{12},C_{23}] + [C_{23},C_{34}] - [C_{123},C_{34}] - [C_{12},C_{234}] + [C_{123},C_{234}] = 0,
\end{eqnarray}
\begin{eqnarray}
\frac12\big[C_{34}\,,\,[C_{12},C_{23}]\big] &=&  
C_{12}\big(C_{23}+C_{34}-C_{234}-C_3\big) +C_{23}\big(C_{34}-C_{1234}\big) 
-C_{34}\,C_{2}
\nonumber\\
&&+C_{123}\big(C_{234}-C_{34}-C_{2}\big)-C_{234}\,C_3+(C_2+C_3)\,C_{1234}+C_2\,C_3\,,
\quad\\
\frac12 \big[C_{23}\,,\,[C_{123},C_{34}]\big]  &=& 
C_{12}\big(-C_{23}+C_{234}-C_4\big) +C_{23}\big(C_{123}-C_{4}\big) 
+C_{34}\big(C_{23}-C_{1}\big)
\nonumber\\
&&+C_{123}\big(C_{34}-C_{234}-C_{3}\big)-C_{234}\,C_3+(C_1+C_3)\,C_{4}+C_1\,C_3\,,
\\
\frac12 \big[C_{234}\,,\,[C_{23},C_{12}]\big]  &=& 
C_{12}\big(C_{234}-C_4\big) +C_{23}\big(C_{12}-C_{34}+C_{234}-C_{1}\big) 
-C_{34}\,C_{1}
\nonumber\\
&&+C_{123}\big(C_{34}-C_{234}-C_{2}\big)-C_{234}\,C_2+(C_1+C_2)\,C_{4}+C_1\,C_2\,,
\\
\frac12 \big[C_{234}\,,\,[C_{123},C_{34}]\big]  &=& 
C_{12}\big(-C_{23}+C_{234}+C_4\big) 
+C_{34}\big(C_{23} -C_{123} -C_{234} +C_{1234}\big)
\nonumber\\
&&+C_{23}\,C_{1234} +C_{123}\big(-C_{234}+C_{2}\big)+C_{234}\,C_4
\nonumber\\
&&-(C_2+C_{1234})\,C_{4}-C_2\,C_{1234}\,.
\end{eqnarray}

In the contiguous basis, the four relations \eqref{eq:comm4bis} read
\begin{eqnarray}
&&
\begin{aligned}
&\big(C_{12}+C_{234}-C_{2}-C_{1234}\big)\,[C_{12},C_{23}]
+\big(C_{12}-C_2+C_1\big)\,[C_{23},C_{34}]
\\
&\qquad\quad +\big(C_{123}-C_{23}-C_{12}+C_2\big)\,[C_{12},C_{234}]
+\big(C_{12}-C_2-C_1\big)\,[C_{34},C_{123}] 
= 0,
\end{aligned}
\\[1ex]
&&\begin{aligned}
&\big(-C_{234}+C_{34}+C_2\big)\,[C_{12},C_{23}]
+\big(C_{12}-C_1+C_2\big)\,[C_{23},C_{34}] 
\\
&\qquad\quad +\big(C_{23}-C_{2}-C_3\big)\,[C_{12},C_{234}]
+2C_2\,[C_{34},C_{123}]
= 0,
\end{aligned}
\\[1ex]
&&\begin{aligned}
&\big(-C_{34}-C_3+C_4\big)\,[C_{12},C_{23}] +\big(C_{123}-C_{12}-C_3\big)\,[C_{23},C_{34}] 
\\
&\qquad\quad +2C_3\,[C_{12},C_{234}] +\big(C_{23}-C_{2}-C_3\big)\,[C_{34},C_{123}]
= 0,
\end{aligned}
\\[1ex]
&&\begin{aligned}
&\big(-C_{34}+C_3-C_4\big)\,[C_{12},C_{23}]
+\big(-C_{123}-C_{34}+C_3+C_{1234}\big)\,[C_{23},C_{34}] 
\\
&\qquad\quad +\big(C_{34}-C_{4}-C_3\big)\,[C_{12},C_{234}]
+\big(C_{234}-C_{23}-C_{34}+C_3\big)\,[C_{34},C_{123}]
= 0.
\end{aligned}
\end{eqnarray}

\paragraph{Casimir elements.} As mentioned in section \ref{sec:sR4}, 
in addition to the central elements $C_{j}$, $j=1,2,3,4$ and $C_{1234}$, 
one shows that the combinations $w_{123}$, $w_{124}$, $w_{134}$, $w_{234}$ and $x_{1234}$, as given in 
\eqref{eq:w123} and \eqref{eq:x1234} (with the expressions \eqref{eq:Ctrou}), are Casimir elements.  

\paragraph{Other presentations.} There exist different generating sets for the $R(4)$ algebra. For instance, using relation \eqref{eq:defCI}, the elements with three or more indices (\textit{i.e.} $C_{123}$, $C_{124}$, $C_{134}$, $C_{234}$ and $C_{1234}$) could be expressed only in terms of the elements 
of type $C_{ij}$ and $C_i$, which would be an alternative basis of $R(4)$.

\section{Finite-dimensional real symmetric representation of $sR(4)$ \label{app:C}}

In this appendix, we provide the explicit expressions of the matrix entries for the representation of $sR(4)$ in Theorem \ref{th:repfin}.

\begin{align}
\psi_{n,p} &= \Big[{\frac{n(2j^{(1)}+1-n)(2j^{(2)}+1-n)(2j^{(1)}+2j^{(2)}+2-n)(N-n-p+2j^{(3)}+2)(N-n-p+1) }{(2j^{(1)}+2j^{(2)} +2-2n)^2(2j^{(1)}+2j^{(2)}+1-2n)(2j^{(1)}+2j^{(2)}+3-2n) } } \nonumber \\
 &\qquad\times { (p-n-N+2j^{(1)}+2j^{(2)}+1)(n-p-N +2j^{(0)}+2j^{(4)}) }\Big]^{\frac12}\,, \label{eq:psifapp}
\end{align}

\begin{align}
\rho_{n,p} &= \Big[{\frac{p(2j^{(4)}+1-p)(2j^{(0)}+1-p)(2j^{(4)}+2j^{(0)}+2-p)(N-n-p+2j^{(3)}+2)(N-n-p+1) }{(2j^{(0)}+2j^{(4)} +2-2p)^2(2j^{(0)}+2j^{(4)}+1-2p)(2j^{(0)}+2j^{(4)}+3-2p) } } \nonumber \\
 &\qquad\times { (p-n-N+2j^{(1)}+2j^{(2)})(n-p-N +2j^{(0)}+2j^{(4)}+1) }\Big]^{\frac12}\,, \label{eq:rhofapp} 
\end{align}

\begin{align}
\varphi^{D}_{n,p} &= \Big[{
\frac{n(2j^{(1)}+1-n)(2j^{(2)}+1-n)(2j^{(1)}+2j^{(2)}+2-n) }
     {(2j^{(1)}+2j^{(2)} +2-2n)^2(2j^{(1)}+2j^{(2)}+1-2n)(2j^{(1)}+2j^{(2)}+3-2n) } } \nonumber \\
&\qquad\times {
\frac{p(2j^{(4)}+1-p)(2j^{(0)}+1-p)(2j^{(4)}+2j^{(0)}+2-p) }
     {(2j^{(0)}+2j^{(4)} +2-2p)^2(2j^{(0)}+2j^{(4)}+1-2p)(2j^{(0)}+2j^{(4)}+3-2p) } } \nonumber \\
 &\qquad\times { (N-n-p+2j^{(3)}+2)(N-n-p+2j^{(3)}+3)(N-n-p+1)(N-n-p+2) }\Big]^{\frac12} \nonumber \\
 &\quad\times \sgn\big((n-p-2j_1-2j_2+N+2j_3)(p-n+2j_1+2j_2-N+1)\big) \,,\label{eq:phiDapp} 
\end{align}

\begin{align}
\varphi^{A}_{n,p} &= \Big[{
\frac{n(2j^{(1)}+1-n)(2j^{(2)}+1-n)(2j^{(1)}+2j^{(2)}+2-n)(n-p-2j^{(1)}-2j^{(2)}+2j^{(3)}+N) }
     {(2j^{(1)}+2j^{(2)} +2-2n)^2(2j^{(1)}+2j^{(2)}+1-2n)(2j^{(1)}+2j^{(2)}+3-2n) } } \nonumber \\
&\qquad\times {
\frac{p(2j^{(4)}+1-p)(2j^{(0)}+1-p)(2j^{(4)}+2j^{(0)}+2-p)(p-n-2j^{(0)}-2j^{(4)}+2j^{(3)}+N) }
     {(2j^{(0)}+2j^{(4)} +2-2p)^2(2j^{(0)}+2j^{(4)}+1-2p)(2j^{(0)}+2j^{(4)}+3-2p) } } \nonumber \\
 &\qquad\times { (p-n-N+2j^{(1)}+2j^{(2)}+1)(n-p-N +2j^{(0)}+2j^{(4)}+1) } \Big]^{\frac12}\nonumber \\
 &\quad\times   \sgn(n+p-N-2j_3-2) \,,\label{eq:phiAapp} 
\end{align}

\begin{align}
\varphi^{H}_{n,p} &= \frac{(p-j^{(4)}-j^{(0)}-1)(p-j^{(4)}-j^{(0)})-j^{(4)}(j^{(4)}+1)+j^{(0)}(j^{(0)}+1)}{2(p-j^{(4)}-j^{(0)}-1)(p-j^{(4)}-j^{(0)})} \nonumber \\
&\times \Big[{\frac{n(2j^{(1)}+1-n)(2j^{(2)}+1-n)(2j^{(1)}+2j^{(2)}+2-n)(N-n-p+2j^{(3)}+2)(N-n-p+1) }{(2j^{(1)}+2j^{(2)} +2-2n)^2(2j^{(1)}+2j^{(2)}+1-2n)(2j^{(1)}+2j^{(2)}+3-2n) } } \nonumber \\
 &\qquad\times { (p-n-N+2j^{(1)}+2j^{(2)}+1)(n-p-N +2j^{(0)}+2j^{(4)}) }\Big]^{\frac12}\,, \label{eq:phiHapp}
\end{align}

\begin{align}
\varphi^{V}_{n,p} &= \frac{(n-j^{(1)}-j^{(2)}-1)(n-j^{(1)}-j^{(2)})+j^{(1)}(j^{(1)}+1)-j^{(2)}(j^{(2)}+1)}{2(n-j^{(1)}-j^{(2)}-1)(n-j^{(1)}-j^{(2)})} \nonumber \\
&\times \Big[{\frac{p(2j^{(4)}+1-p)(2j^{(0)}+1-p)(2j^{(4)}+2j^{(0)}+2-p)(N-n-p+2j^{(3)}+2)(N-n-p+1) }{(2j^{(0)}+2j^{(4)} +2-2p)^2(2j^{(0)}+2j^{(4)}+1-2p)(2j^{(0)}+2j^{(4)}+3-2p) } } \nonumber \\
 &\qquad\times { (p-n-N+2j^{(1)}+2j^{(2)})(n-p-N +2j^{(0)}+2j^{(4)}+1) }\Big]^{\frac12}\,, \label{eq:phiVapp} 
\end{align}
where $\sgn$ is the sign function.
One defines also 
\begin{align}
 \widetilde{\varphi}^{H}_{n,p} &= \varphi^{H}_{n,p} - \psi_{n,p}\,, \label{eq:varphit}\\
 \widetilde{\varphi}^{V}_{n,p} &= \varphi^{V}_{n,p} - \rho_{n,p}\,.
\end{align}
Let us recall that 
\begin{align}
 \varphi^{00}_{np}&=\frac{(j^{(1)}(j^{(1)}+1)-j^{(2)}(j^{(2)}+1))( \rho^{00}_{np}-j^{(0)}(j^{(0)}+1) )}{2(n-j^{(1)}-j^{(2)}-1)(n-j^{(1)}-j^{(2)})}\\
 &+\frac{j^{(1)}(j^{(1)}+1)+j^{(2)}(j^{(2)}+1)+j^{(0)}(j^{(0)}+1)-(n-j^{(1)}-j^{(2)}-1)(n-j^{(1)}-j^{(2)})+  \rho^{00}_{np} }{2} \nonumber\\
 &=\frac{(j^{(0)}(j^{(0)}+1)-j^{(4)}(j^{(4)}+1))( \psi^{00}_{np}-j^{(1)}(j^{(1)}+1) )}{2(p-j^{(0)}-j^{(4)}-1)(p-j^{(0)}-j^{(4)})}\\
 &+\frac{j^{(1)}(j^{(1)}+1)+j^{(4)}(j^{(4)}+1)+j^{(0)}(j^{(0)}+1)-(p-j^{(0)}-j^{(4)}-1)(p-j^{(0)}-j^{(4)})+  \psi^{00}_{np} }{2}\,. \nonumber
\end{align}
We define 
\begin{align}
  \varphi^{0\overline{0}}_{np}&=\frac{(j^{(1)}(j^{(1)}+1)-j^{(2)}(j^{(2)}+1))( \rho^{00}_{np}-j^{(0)}(j^{(0)}+1) )}{2(n-j^{(1)}-j^{(2)}-1)(n-j^{(1)}-j^{(2)})} \label{eq:min1}\\
 &-\frac{j^{(1)}(j^{(1)}+1)+j^{(2)}(j^{(2)}+1)+j^{(0)}(j^{(0)}+1)-(n-j^{(1)}-j^{(2)}-1)(n-j^{(1)}-j^{(2)})+  \rho^{00}_{np} }{2}\,, \nonumber\\
  \varphi^{\overline{0}0}_{np}&=  \frac{(j^{(0)}(j^{(0)}+1)-j^{(4)}(j^{(4)}+1))( \psi^{00}_{np}-j^{(1)}(j^{(1)}+1) )}{2(p-j^{(0)}-j^{(4)}-1)(p-j^{(0)}-j^{(4)})} \label{eq:min2}\\
 &-\frac{j^{(1)}(j^{(1)}+1)+j^{(4)}(j^{(4)}+1)+j^{(0)}(j^{(0)}+1)-(p-j^{(0)}-j^{(4)}-1)(p-j^{(0)}-j^{(4)})+  \psi^{00}_{np} }{2} \,,\nonumber  \\
  \varphi^{\overline{0}\overline{0}}_{np}&=\varphi^{00}_{np}-\rho^{00}_{np}-\psi^{00}_{np}+j^{(2)}(j^{(2)}+1)+j^{(3)}(j^{(3)}+1)+j^{(4)}(j^{(4)}+1)\,. \label{eq:min3}
\end{align}

\section{Racah polynomials \label{app:r}}

The Racah polynomials are defined through the following hypergeometric functions \cite{Koek}, for a given positive integer $\cN$ and $0\leq n,m \leq \cN$,
\begin{equation}
r_n(m;\alpha,\beta,-\cN-1,\delta) = {_4F_3}\Big[ \genfrac{}{}{0pt}{0}{-n,n+\alpha+\beta+1,-m,m-\cN+\delta}{\alpha+1,\beta+\delta+1,-\cN} \Big\vert 1 \Big] \,,
\end{equation}
which corresponds to the case $\gamma=-\cN-1$ in \cite{Koek}.
They satisfy the recurrence relation
\begin{equation} \label{eq:recuR}
 \lambda(m)\, r_n(m)=A_{n}\  r_{n+1}(m) -(A_n+C_n)\,r_n(m)+ C_n \ r_{n-1}(m)\,,
\end{equation}
with $\lambda(m)=m(m-\cN+\delta)$,
and the difference equation
\begin{equation}\label{eq:diffR}
 n(n+\alpha+\beta+1)\, r_n(m)=B_m \ r_{n}(m+1) -(B_m+D_m )\, r_n(m)+ D_m\ r_n(m-1)\,,
\end{equation}
where $r_n(m)$ stands for $r_n(m;\alpha,\beta,-\cN-1,\delta)$ and 
\begin{align}
\label{eq:B4}
 &A_n=\frac{(n+\alpha+1)(n+\beta+\delta+1)(n-\cN)(n+\alpha+\beta+1)}{(2n+\alpha+\beta+1)(2n+\alpha+\beta+2)},\\
\label{eq:B5}
 &C_n=\frac{n(n+\alpha+\beta+\cN+1)(n+\alpha-\delta)(n+\beta)}{(2n+\alpha+\beta)(2n+\alpha+\beta+1)},\\
\label{eq:B6}
 &B_m=\frac{(m+\alpha+1)(m+\beta+\delta+1)(m-\cN)(m-\cN+\delta)}{(2m-\cN+\delta)(2m-\cN+\delta+1)},\\
\label{eq:B7}
 &D_m=\frac{m(m-\cN-1-\alpha+\delta)(m-\cN-1-\beta)(m+\delta)}{(2m-\cN-1+\delta)(2m-\cN+\delta)}.
\end{align}
The parameters $\alpha,\ \beta$ and $\delta$ are chosen in such a way that the previous coefficients satisfy 
$A_nC_{n+1}>0$ and $B_mD_{m+1}>0$, for any $0\leq n,m < \cN$.

Using the invariance of the hypergeometric function w.r.t. the permutations of its parameters, one can show that
\begin{eqnarray}
 r_n(x;\alpha,\beta,-\cN-1,\delta)
 &=&r_n(x;\beta+\delta,\alpha-\delta,-\cN-1,\delta) \label{eq:3}\\
 &=& r_x(n;\alpha, \delta-\alpha-\cN-1, -\cN-1,\alpha+\beta+\cN+1)\,.\label{eq:dua}
\end{eqnarray}
The last equality provides the proof of the duality between $x$ and $n$.
There is also the Whipple relation \cite{Wipple} for the ${_4F_3}$ hypergeometric function, which reads as follows  
\begin{equation}
  r_n(x;\alpha,\beta,-\cN-1,\delta)=\frac{(\alpha-\delta+1)_n(\beta+1)_n}{(\beta+\delta+1)_n(\alpha+1)_n}\ r_n(\cN-x;\beta,\alpha,-\cN-1,-\delta),\label{eq:wh}
\end{equation}
where $(y)_n=y(y+1)\dots(y+n-1)$ is the Pochhammer symbol (by convention $(y)_0=1$). This previous equality is the symmetry of the function $r_n$ w.r.t. the transformation of the variable $x \rightarrow \cN-x$. 

The functions $r$ satisfies also the following property
\begin{equation}\label{eq:reln1}
\widetilde{\lambda}(m)\ r_n(m)=E_{n}\ \widetilde{r}_{n+1}(m) +F_n\  \widetilde{r}_{n}(m)  +G_n\  \widetilde{r}_{n-1}(m) \,,
\end{equation}
where $\widetilde{r}_{n}(m)=r_n(m;\alpha,\beta,-\cN-2,\delta+1)$ and 
\begin{align}
\label{eq:B9}
 &\widetilde{\lambda}(m)= \frac{(m+1+\delta)(\cN-m+1)(\beta+\delta+1)}{\cN+1}\,,\\
\label{eq:B10}
& E_n= \frac{n+2+\beta+\delta}{n-\cN}A_n\ , \qquad  G_n=  \frac{n-1+\alpha-\delta}{n+\alpha+\beta+\cN+1}C_n\ ,\\
\label{eq:B11}
& F_n=-E_n-G_n+(\delta+1)(\beta+\delta+1)\, .
\end{align}
Relation \eqref{eq:reln1} is proven for small $n$ by direct computation. Assuming that it holds up to a given $n$,
 we transform the l.h.s. of \eqref{eq:reln1} for $n+1$ as follows
\begin{eqnarray}
\widetilde{\lambda}({m})\ r_{n+1}(m) &=&\frac{\widetilde{\lambda}({m}) }{A_n}\Big((\lambda(m) +A_n+C_n)\ r_n(m) - C_n\ r_{n-1}(m) \Big)\\
&=&\frac{1 }{A_n}\Big( (\lambda(m) +A_n+C_n)(E_{n}\ \widetilde{r}_{n+1}(m) +F_n\  \widetilde{r}_{n}(m)  +G_n\  \widetilde{r}_{n-1}(m))
\quad \\
&&\qquad -C_n (E_{n-1}\ \widetilde{r}_{n}(m) +F_{n-1}\  \widetilde{r}_{n-1}(m)  +G_{n-1}\  \widetilde{r}_{n-2}(m)) \Big)\,.\nonumber
\end{eqnarray}
We have used the recurrence relation \eqref{eq:recuR}, then the recursion hypothesis.
The terms $\lambda({m})\widetilde r_j(m)$ ($j=n+1,n,n-1$) are replaced thanks to the recurrence relation \eqref{eq:recuR} and, using the explicit expressions of $A_n$, $C_n$, $E_n$ and $F_n$,
this expression reproduces the r.h.s. of \eqref{eq:reln1} for $n+1$. It proves \eqref{eq:reln1} for any $n$.

In the same way, introducing $\widetilde{r}\,'_{n}(m)=r_n(m;\alpha,\beta,-\cN,\delta-1)$, one has the property
\begin{equation}\label{eq:reln2}
\widetilde{\lambda}\,'(m)\ r_n(m)=E'_{n}\ \widetilde{r}\,'_{n+1}(m) +F'_n\  \widetilde{r}\,'_{n}(m)  +G\,'_n\  \widetilde{r}\,'_{n-1}(m) \,,
\end{equation}
where 
\begin{align}
&\widetilde{\lambda}\,'(m)= \frac{\cN(m+\beta+\delta)(\cN-m+\beta)}{\beta+\delta} \ ,\\ 
& E'_n= \frac{n-\cN+1}{n+1+\beta+\delta}A_n\ , \quad  G'_n=  \frac{n+\alpha+\beta+\cN}{n+\alpha-\delta}C_n\ ,\quad
 F'_n=-E'_n-G'_n+\cN(\beta+\cN)\, .
\end{align}

The functions $r$ satisfies also the following properties
\begin{equation}\label{eq:reln3}
\widetilde{\mu}(n)\ r_n(m)=H_{m}\ \widetilde{r}_{n}(m+1) +I_m\  \widetilde{r}_{n}(m)  +J_m\  \widetilde{r}_{n}(m-1) \,,
\end{equation}
where 
\begin{align}
\label{eq:B12}
 &\widetilde{\mu}(n)= \frac{(n+2+\alpha+\beta+\cN)(\cN-n+1)(\beta+\delta+1)}{\cN+1}\,,\\
\label{eq:B13}
& H_m= \frac{m+2+\beta+\delta}{m-\cN}B_m\ , \qquad  J_m=  \frac{m-\cN-\beta-1}{m+\delta}D_m\ ,\\
\label{eq:B14}
& I_m=-H_m-J_m+(\alpha+\beta+2+\cN)(\beta+\delta+1)\,,
\end{align}
and
\begin{equation}\label{eq:reln4}
\widetilde{\mu}\,'(n)\ r_n(m)=H'_{m}\ \widetilde{r}\,'_{n}(m+1) +I'_m\  \widetilde{r}\,'_{n}(m)  +J'_m\  \widetilde{r}\,'_{n}(m-1) \,,
\end{equation}
where 
\begin{align}
\label{eq:B15}
 &\widetilde{\mu}\,'(n)= \frac{\cN(n+\beta+\delta)(n+\alpha-\delta+1)}{\beta+\delta}\,,\\
\label{eq:B16}
& H'_m= \frac{\cN-m-1}{m+\beta+\delta+1}B_m\ , \qquad  J'_m=  \frac{\delta+m-1}{\cN-m+\beta+1}D_m\ ,\\
& I'_m=-H'_m-J'_m+(\alpha-\delta+1)\cN\, .
\end{align}

Let us define 
\begin{equation}\label{eq:qR}
 \cP_n(m) = \sqrt{\frac{(\alpha-\delta+1)_\cN(\beta+1)_\cN}{(\alpha+\beta+2)_\cN(-\delta)_\cN}} \, \prod_{i=0}^{n-1}\frac{A_i}{\sqrt{A_iC_{i+1}}}\prod_{j=0}^{m-1}\frac{B_j}{\sqrt{B_jD_{j+1}}}\,
 r_n(m;\alpha,\beta,-\cN-1,\delta)\,,
\end{equation}
where $\cP_n(m)$ stands for $\cP_n(m;\alpha,\beta,-\cN-1,\delta)$.
Let us emphasize that we do not simplify the normalization coefficients in \eqref{eq:qR} since we do not fix the sign of $A_i$ or $B_j$. 
These polynomials satisfy the recurrence relation
\begin{equation} \label{eq:recuq}
 m(m-\cN+\delta)\ \cP_n(m)=\sqrt{A_{n}C_{n+1}}\  \cP_{n+1}(m) -(A_n+C_n)\ \cP_n(m)+ \sqrt{A_{n-1}C_n} \ \cP_{n-1}(m)\,,
\end{equation}
and the difference equation
\begin{equation}\label{eq:diffq}
 n(n+\alpha+\beta+1)\ \cP_n(m)=\sqrt{B_mD_{m+1}} \ \cP_{n}(m+1) -(B_m+D_m )\ \cP_n(m)+ \sqrt{B_{m-1}D_m}\ \cP_n(m-1)\,.
\end{equation}
The orthogonality relation reads
\begin{equation}\label{eq:ortho3}
 \sum_{m=0}^{\cN} \cP_{n}(m)\,\cP_{n'}(m) = \delta_{n,n'}\,.
\end{equation}

\end{appendix}

\end{document}